\documentclass[draft]{article}

\usepackage{amssymb}
\usepackage{amsmath}
\usepackage{amsthm}
\usepackage{graphics}

\newcommand{\diag}{\mathop{\mathrm {diag}}\nolimits}

\newcommand{\Hom}{\mathop{\mathrm {Hom}}\nolimits}

\newcommand{\sI}{\sqrt{-1}}

\newcommand{\cF}{\mathcal{F}}

\newcommand{\cJ}{\mathcal{J}}

\newcommand{\cP}{\mathcal{P}}
\newcommand{\cQ}{\mathcal{Q}}

\newcommand{\cS}{\mathcal{S}}
\newcommand{\cT}{\mathcal{T}}

\newcommand{\cW}{\mathcal{W}}

\newcommand{\rC}{\mathrm{C}}

\newcommand\rG {\mathrm{G}}

\newcommand{\rI}{\mathrm{I}}

\newcommand{\rM}{\mathrm{M}}

\newcommand{\rO}{\mathrm{O}}
\newcommand{\rP}{\mathrm{P}}

\newcommand{\rR}{\mathrm{R}}
\newcommand{\rS}{\mathrm{S}}
\newcommand{\rU}{\mathrm{U}}
\newcommand{\rW}{\mathrm{W}}
\newcommand{\ra}{\mathrm{a}}
\newcommand{\rb}{\mathrm{b}}
\newcommand{\rc}{\mathrm{c}}

\newcommand{\rf}{\mathrm{f}}
\newcommand{\rg}{\mathrm{g}}
\newcommand{\rk}{\mathrm{k}}

\newcommand{\rp}{\mathrm{p}}

\newcommand{\rr}{\mathrm{r}}

\newcommand{\ru}{\mathrm{u}}

\newcommand{\bC}{\mathbb{C}}

\newcommand{\bN}{\mathbb{N}}
\newcommand{\bQ}{\mathbb{Q}}
\newcommand{\bR}{\mathbb{R}}
\newcommand{\bZ}{\mathbb{Z}}

\newcommand{\mW}{\mathbf{W}}

\newcommand{\me}{\mathbf{e}}

\newcommand{\gge}{\mathfrak{e}}
\newcommand{\g }{\mathfrak{g}}

\newcommand{\gl}{\mathfrak{l}}

\newcommand{\go}{\mathfrak{o}}

\newcommand{\gs}{\mathfrak{s}}
\newcommand{\gu}{\mathfrak{u}}

\newtheorem{thm}{Theorem}[section]
\newtheorem{lem}[thm]{Lemma}
\newtheorem{prop}[thm]{Proposition}
\newtheorem{cor}[thm]{Corollary}

\newtheorem{rem}[thm]{Remark}

\numberwithin{equation}{section}
\theoremstyle{remark}

\begin{document}
\allowdisplaybreaks


\title{Calculus of archimedean Rankin--Selberg integrals 
with recurrence relations}

\author{Taku Ishii\footnote{Faculty of Science and Technology, 
Seikei University, 3-3-1 Kichijoji-kitamachi, 
Musashino, Tokyo, 180-8633, Japan. 
E-mail: \texttt{ishii@st.seikei.ac.jp}}, 
and Tadashi Miyazaki\footnote{Department of Mathematics, 
College of Liberal Arts and Sciences, Kitasato University, 
1-15-1 Kitasato, Minamiku, Sagamihara, Kanagawa, 252-0373, Japan. 
E-mail: \texttt{miyaza@kitasato-u.ac.jp}}}

\maketitle


\begin{abstract}
Let $n$ and $n'$ be positive integers such that 
$n-n'\in \{0,1\}$. 
Let $F$ be either $\bR$ or $\bC$. 
Let $K_n$ and $K_{n'}$ be maximal compact subgroups of 
$\mathrm{GL}(n,F)$ and $\mathrm{GL}(n',F)$, respectively. 
We give the explicit descriptions 
of archimedean Rankin--Selberg integrals at the 
minimal $K_n$- and $K_{n'}$-types for pairs of 
principal series representations of 
$\mathrm{GL}(n,F)$ and $\mathrm{GL}(n',F)$, 
using their recurrence relations. 
Our results for $F=\bC$ can be applied to 
the arithmetic study of critical values of 
automorphic $L$-functions. 
\end{abstract}

\section{Introduction}
\label{sec:intro}

The theory of automorphic $L$-functions 
via integral representations has its origin in the work of 
Hecke \cite{Hecke_001} for $\mathrm{GL}(2)$, 
and the works of Rankin \cite{Rankin_001}, 
Selberg \cite{Selberg001} for $\mathrm{GL}(2)\times \mathrm{GL}(2)$. 
As a direct outgrowth of their works, 
the theory of Rankin--Selberg integrals for 
$\mathrm{GL}(n)\times \mathrm{GL}(n')$ were developed by 
Jacquet, Piatetski-Shapiro, and Shalika \cite{J_P_S_S_001}. 
Our interest here is the archimedean local theory of 
their Rankin--Selberg integrals. 

Let $F$ be either $\bR$ or $\bC$. 
We fix a maximal compact subgroup $K_n$ of $\mathrm{GL}(n,F)$. 
Let $\Pi $ and $\Pi'$ be irreducible generic 
Casselman--Wallach representations of $\mathrm{GL}(n,F)$ and 
$\mathrm{GL}(n',F)$, respectively. 
We denote by $L(s,\Pi \times \Pi')$ the archimedean $L$-factor 
for $\Pi\times \Pi'$. 
The theory of 
archimedean Rankin--Selberg integrals for $\Pi\times \Pi'$ 
was developed by Jacquet and Shalika. 
In \cite{Jacuqet_Shalika_001}, they 
showed that any archimedean Rankin--Selberg integral 
for $\Pi \times \Pi'$ is extended 
to $\bC$ as a holomorphic multiple of 
$L(s,\Pi \times {\Pi'})$, is bounded at infinity in vertical strips, 
and satisfies the local functional equation. 
Furthermore, in  \cite{Jacquet_001}, Jacquet refined 
the proofs of the above results, and showed further that 
$L(s,\Pi \times {\Pi'})$ can be expressed 
as a linear combination of 
archimedean Rankin--Selberg integrals for $\Pi \times {\Pi'}$ 
if $n-n'\in \{0,1\}$. 
Their results are sufficient for the proofs of 
important analytic properties of automorphic $L$-functions such as 
the analytic continuations, the functional equations and 
the converse theorems. 
However, in the studies of arithmetic properties of automorphic $L$-functions, 
the precise knowledges of archimedean Rankin--Selberg integrals 
at the special $K_n$- and $K_{n'}$-types 
are required. 
For example, Sun's non-vanishing result \cite{Sun_001} 
at the minimal $K_n$- and $K_{n-1}$-types 
is vital to the arithmetic study of critical values
of automorphic $L$-functions for $\mathrm{GL}(n)\times \mathrm{GL}(n-1)$ 
by the cohomological method. 
The goal of this paper is to give explicit descriptions of 
archimedean Rankin--Selberg integrals at the minimal 
$K_n$- and $K_{n'}$-types for pairs of 
principal series representations of 
$\mathrm{GL}(n,F)$ and $\mathrm{GL}(n',F)$ 
with $n-n'\in \{0,1\}$. 
We generalize Stade's results \cite{Stade_002}, \cite{Stade_001} 
(see \cite{Ishii_Stade_002} for a simplified proof) 
for the spherical case to general case.

Let us explain our main result 
for the $\mathrm{GL}(n)\times \mathrm{GL}(n-1)$-case. 
Assume that $\Pi$ and $\Pi'$ are irreducible 
principal series representations of $\mathrm{GL}(n,F)$ and 
$\mathrm{GL}(n-1,F)$, respectively. 
Let $\psi $ be the standard additive character of $F$. 
We regard $\mathrm{GL}(n-1,F)$ as a subgroup of 
$\mathrm{GL}(n,F)$ via the embedding 
\begin{align*}
&\iota_n\colon \mathrm{GL}(n-1,F)\ni g\mapsto 
\begin{pmatrix}
g& \\ 
 &1
\end{pmatrix}\in \mathrm{GL}(n,F).
\end{align*}
For $W\in \cW (\Pi ,\psi )$ and 
$W'\in \cW (\Pi',\psi^{-1})$, 
we define the archimedean Rankin--Selberg integral 
$Z(s,W,W')$ by 
\begin{align*}
&Z(s,W,W')=\int_{N_{n-1}\backslash \mathrm{GL}(n-1,F)}
W(\iota_n(g))W'(g)|\det g|_F^{s-1/2}\,dg\qquad 
(\mathrm{Re}(s)\gg 0), 
\end{align*}
where $\cW (\Pi ,\psi )$, $\cW (\Pi',\psi^{-1} )$ 
are the Whittaker models of $\Pi $, $\Pi'$, 
respectively, 
$N_{n-1}$ is the upper triangular unipotent subgroup of 
$\mathrm{GL}(n-1,F)$, 
and $|\cdot |_F$ is the usual norm on $F$. 
Let $(\tau_{\mathrm{min}} ,V_{\mathrm{min}})$ be 
the minimal $K_n$-type of $\Pi$, and we fix 
a $K_n$-embedding $\mW \colon V_{\mathrm{min}}
\to \cW (\Pi ,\psi )$. 
Let $(\tau_{\mathrm{min}}',V_{\mathrm{min}}')$ be 
the minimal $K_{n-1}$-type of $\Pi'$, and we fix 
a $K_{n-1}$-embedding $\mW'\colon V_{\mathrm{min}}'
\to \cW (\Pi',\psi^{-1} )$. 
Here we give $\mW$ and $\mW'$ concretely 
by the Jacquet integrals. We note that 
\[
V_{\mathrm{min}}\otimes_\bC V_{\mathrm{min}}'\ni 
v\otimes v'\mapsto Z(s,\mW (v),\mW'(v'))\in 
\bC_{\mathrm{triv}}
\]
defines an element of 
$\Hom_{K_{n-1}}(V_{\mathrm{min}}\otimes_\bC V_{\mathrm{min}}', 
\bC_{\mathrm{triv}})$, where 
$\bC_{\mathrm{triv}}=\bC$ is the trivial $K_{n-1}$-module. 
In the first main theorem (Theorem \ref{thm:main1}), 
we give the explicit description of 
this $K_{n-1}$-homomorphism. 
More precisely, under the assumption 
$\Hom_{K_{n-1}}(V_{\mathrm{min}}\otimes_\bC V_{\mathrm{min}}', 
\bC_{\mathrm{triv}})\neq \{0\}$, 
we show the equality 
\begin{equation}
\label{eqn:intro_001}
Z(s,\mW(v),\mW'(v'))
=L(s,\Pi \times \Pi')\Psi (v\otimes v')\qquad 
(v\in V_{\mathrm{min}},\ v'\in V_{\mathrm{min}}')
\end{equation}
with some nonzero $\Psi \in 
\Hom_{K_{n-1}}(V_{\mathrm{min}}\otimes_\bC V_{\mathrm{min}}', 
\bC_{\mathrm{triv}})$ independent of $s$, 
and describe $\Psi $ explicitly in terms of 
Gelfand--Tsetlin type bases of 
$V_{\mathrm{min}}$ and $V_{\mathrm{min}}'$. 
In the second main theorem (Theorem \ref{thm:main2}), 
we give a similar description 
for the $\mathrm{GL}(n)\times \mathrm{GL}(n)$-case. 
Since the statement of Theorem \ref{thm:main2} is slightly complicated, 
we leave it to Section 2.

We introduce some applications of our results 
(Theorems \ref{thm:main1} and \ref{thm:main2}) 
for $F=\bC$. In the arithmetic study of 
critical values of automorphic $L$-functions 
for $\mathrm{GL}(n)\times \mathrm{GL}(n')$ 
with $n-n'\in \{0,1\}$ by the cohomological method, 
the archimedean Rankin--Selberg integrals 
at the minimal $K_n$- and $K_{n'}$-types play 
important roles, and the hypothesis of the non-vanishing of them 
at critical points is called the non-vanishing hypothesis 
for $\mathrm{GL}(n,F)\times \mathrm{GL}(n',F)$. 
It is known that a local component at the complex place 
of irreducible regular algebraic cuspidal automorphic 
representation of $\mathrm{GL}(n)$ is 
a cohomological principal series representation 
({\it cf.} \cite[Proposition 2.14]{Raghuram_001}). 
Hence, Theorem \ref{thm:main1} gives 
another proof of the non-vanishing hypothesis 
for $\mathrm{GL}(n,\bC )\times \mathrm{GL}(n-1,\bC )$ 
at all critical points, 
which were originally proved by Sun \cite{Sun_001} 
and were used in Grobner--Harris \cite{Grobner_Harris_001} 
and Raghuram \cite{Raghuram_001}. 
In \cite{Dong_Xue_001}, 
Dong and Xue proved the non-vanishing hypothesis 
for $\mathrm{GL}(n,\bC )\times \mathrm{GL}(n,\bC )$ 
at the central critical point, and they indicate that 
it is hard to generalize their result to all critical points by 
the technique of the translation functor. 
Theorem \ref{thm:main2} proves the non-vanishing hypothesis 
for $\mathrm{GL}(n,\bC )\times \mathrm{GL}(n,\bC )$ at all critical points, 
and allows us to improve the archimedean part of Greni\'{e}'s theorem 
\cite[Theorem 2]{Grenie_001} into more explicit form 
({\it cf.} Remark \ref{rem:cohomological_2}). 
We expect that our explicit results will be applied to deeper study of 
special values of automorphic $L$-functions. 

There are some related results to be mentioned here. 
In the cases of $\mathrm{GL}(n)\times \mathrm{GL}(n-1)$ and 
$\mathrm{GL}(n)\times \mathrm{GL}(n)$, 
we expect that the archimedean Rankin--Selberg integrals 
for appropriate Whittaker functions are equal to the associated $L$-factors. 
This expectation was proved 
by Jacquet--Langlands \cite{Jacquet_Langlands_001} and 
Popa \cite{Popa_001} 
for the $\mathrm{GL}(2)\times \mathrm{GL}(1)$-case; 
by Jacquet \cite{Jacquet_003}, Zhang \cite{Zhang_001} and 
the second author \cite{Miyazaki_003} 
for the $\mathrm{GL}(2)\times \mathrm{GL}(2)$-case; 
by Hirano and the authors \cite{Hirano_Ishii_Miyazaki_pre}  
for the $\mathrm{GL}(3)\times \mathrm{GL}(2)$-case. 
The results of this paper, that is, the formula 
(\ref{eqn:intro_001}) and the analogous formula 
for $\mathrm{GL}(n)\times \mathrm{GL}(n)$ in Theorem \ref{thm:main2} 
can be regarded as additional evidences of this expectation 
for the higher rank cases. 

Let us briefly explain the idea of the proofs of our main theorems. 
The key ingredients are two kinds of the Godement sections 
for a principal series representation $\Pi$ of $\mathrm{GL}(n,F)$. 
One is defined by Jacquet \cite{Jacquet_001} as 
an integral transform of the standard section for some 
principal series representation of $\mathrm{GL}(n-1,F)$, 
and gives a recursive integral representation of 
a Whittaker function for $\Pi$. 
The other seems to be new, 
and is defined as 
an integral transform of the standard section for 
the same representation $\Pi$ of $\mathrm{GL}(n,F)$. 
It gives an integral representation of 
a Whittaker function for $\Pi$, which is related to 
the local theta correspondence in 
Watanabe \cite[\S 2]{Watanabe_001}. 
Using two kinds of the Godement sections, 
we construct the recurrence relations of 
the archimedean Rankin--Selberg integrals 
for pairs of principal series representations of $\mathrm{GL}(n,F)$ 
and $\mathrm{GL}(n',F)$ with $n-n'\in \{0,1\}$. 
Based on the representation theory of $K_n$, 
we write down these recurrence relations at the minimal $K_n$- 
and $K_{n'}$-types, explicitly, and prove the main theorems by induction. 
Here we remark that the explicit recurrence relations for the spherical case 
coincide with those in \cite{Ishii_Stade_002}, 
which follow from explicit formulas of the radial parts of 
spherical Whittaker functions in \cite{Ishii_Stade_001}.

This paper consists of five sections together with an appendix. 
In Section \ref{sec:main}, 
we introduce basic notation and state the main theorems. 
In Section \ref{sec:relation}, 
we define two kinds of the Godement sections 
and give the recurrence relations of 
the archimedean Rankin--Selberg integrals. 
Section \ref{sec:fin_dim} is devoted to 
some preliminary results on the theory of 
finite dimensional representations of $K_n$ and 
$\mathrm{GL}(n,\bC )$. 
In Section \ref{sec:poof_main}, we prove the main theorems using 
the results in Sections \ref{sec:relation} and \ref{sec:fin_dim}. 
As an appendix, we generalize the explicit formulas of the radial parts of 
Whittaker functions in \cite{Ishii_Stade_001} 
using the Godement section. 

This work was supported by JSPS KAKENHI Grant Numbers JP19K03452, JP18K03252.

\section{Main results}
\label{sec:main}

In this section, we introduce basic notation and our main results. 
We describe each objects explicitly as possible, 
although not all of them are necessary to state our main theorems. 
The authors believe that they are of interest and useful 
for further studies.

\subsection{Notation}
\label{subsec:notation}

We denote by $\bZ $, $\bQ $, $\bR $ and $\bC $ the ring of rational integers, 
the rational number field, 
the real number field and the complex number field, respectively. 
Let $\bR_{+}^\times $ be the multiplicative group of positive real numbers. 
Let $\bN_0$ be the set of non-negative integers. 
The real part, the imaginary part and 
the complex conjugate of a complex number $z$ 
are denoted by $\mathrm{Re}(z)$, $\mathrm{Im}(z)$ 
and $\overline{z}$, respectively.

Throughout this paper, $F$ denotes the archimedean local field, 
that is, $F$ is either $\bR$ or $\bC$. 
It is convenient to define the constant $\rc_F$ by 
$\rc_\bR =1$ and $\rc_\bC =2$. We define additive characters 
$\psi_{t} \colon F \to \bC^\times$ ($t\in F$) 
and a norm $|\cdot|_{F}$ on $F$ by 
\begin{align*}
&\psi_{t}(z)=\exp \bigl(\pi \rc_F\sI (tz +\overline{tz})\bigr)
=\left\{\begin{array}{ll}
\exp (2\pi \sI tz)&\text{if $F=\bR$,}\\[1mm]
\exp (2\pi \sI (tz +\overline{tz}))&\text{if $F=\bC$,}
\end{array}\right.&
\end{align*}
and $|z|_{F}=|z|^{\rc_F}$ for $z\in F$, 
where $|\cdot |$ is the ordinary absolute value. 
When $t=\varepsilon \in \{\pm 1\}$, 
we call $\psi_{\varepsilon }$ the standard character of $F$. 
We identify the additive group $F$ with its dual group 
via the isomorphism $t\mapsto \psi_{t}$, 
and denote by $d_{F}z$ the self-dual additive Haar measure on $F$, 
that is, $d_\bR z=dz$ is the ordinary Lebesgue measure on $\bR$ 
and $d_\bC z=2dx\,dy$ ($z=x+\sI y$) is twice 
the ordinary Lebesgue measure on $\bC \simeq \bR^2$.  
For $m\in \bZ$, we define 
a meromorphic function $\Gamma_F(s;m)$ of $s$ in $\bC$ by 
\begin{align*}
&\Gamma_{F}(s;m)
=\rc_F(\pi \rc_F)^{-(s\rc_F+m)/2}
\Gamma \left(\frac{s\rc_F+m}{2}\right)
=\left\{\begin{array}{ll}
\Gamma_\bR (s+m)&\text{if }F=\bR,\\
\Gamma_\bC (s+m/2)&\text{if }F=\bC,
\end{array}\right.
\end{align*}
where $\Gamma_\bR (s)=\pi^{-s/2}\Gamma (s/2)$, 
$\Gamma_\bC (s)=2(2\pi)^{-s}\Gamma (s)$ 
and $\Gamma (s)$ is the usual Gamma function.

Throughout this paper, 
$n$ and $n'$ are positive integers. 
The space of $n\times n'$ matrices over $F$ 
is denoted by $\mathrm{M}_{n,n'}(F)$. 
When $n'=n$, we denote $\mathrm{M}_{n,n}(F)$ 
simply by $\mathrm{M}_{n}(F)$. 
We denote by $d_Fz$ the measure on $\rM_{n,n'}(F)$ defined by 
\begin{align*}
&d_Fz=\prod_{i=1}^n\prod_{j=1}^{n'}d_Fz_{i,j}&
&(z=(z_{i,j})\in \rM_{n,n'}(F)). 
\end{align*}
Let $O_{n,n'}$ be the zero matrix in $\mathrm{M}_{n,n'}(F)$. 
Let $1_n$ be the unit matrix in $\rM_{n}(F)$. 
Let $e_n=(O_{1,n-1},1)\in \rM_{1,n}(F)$. 
When $n=1$, we understand $e_1=1$.

\subsection{Groups and the invariant measures}
\label{subsec:group_algebra}

Let $G_n$ be the general linear group 
$\mathrm{GL}(n,F)$ of degree $n$ over $F$. 
We fix a maximal compact subgroup $K_n$ of $G_n$ by 
\[
K_n=\left\{
\begin{array}{ll}
\mathrm{O}(n)&\text{if }{F} =\bR,\\
\mathrm{U}(n)&\text{if }{F} =\bC,
\end{array}
\right. 
\]
where $\mathrm{O}(n)$ and $\mathrm{U}(n)$ are 
the orthogonal group and the unitary group 
of degree $n$, respectively. 
Let $N_n$ and $U_n$ be the groups of 
upper and lower triangular unipotent matrices in $G_n$, 
respectively, that is, 
\begin{align*}
&N_n=\{x=(x_{i,j})\in G_n\mid 
x_{i,j}=0\ (1\leq j<i\leq n),\ \ 
x_{k,k}=1\ (1\leq k\leq n)\},\\
&U_n=\{u=(u_{i,j})\in G_n\mid 
u_{i,j}=0\ (1\leq i<j\leq n),\ \ 
u_{k,k}=1\ (1\leq k\leq n)\}.
\end{align*}
We define subgroups $M_n$ and $A_n$ of $G_n$ by 
\begin{align*}
&M_n=\{m=\diag (m_1,m_2,\cdots ,m_n)\mid m_i\in G_1=F^\times \quad 
(1\leq i\leq n)\},\\
&A_n=\{a=\diag (a_1,a_2,\cdots ,a_n)\mid a_i\in \bR_{+}^\times  \quad 
(1\leq i\leq n)\}.
\end{align*}
Let $Z_n$ be the center of $G_n$. Then we have 
$Z_n=\{t1_n\mid t\in G_1=F^\times \}$. 
We denote by $C^\infty (G_n)$ 
the space of ($\bC$-valued) smooth functions on $G_n$. 
We regard $C^\infty (G_n)$ as a $G_n$-module via the right translation 
\begin{align*}
&(R(g)f)(h)=f(hg)&
&(g,h\in G_n,\ f\in C^\infty (G_n)). 
\end{align*}

Let $dk$, $dx$, $du$ and $da$ be the Haar measures on $K_n$, $N_n$, $U_n$ 
and $A_n$, respectively. In this paper, we normalize these Haar measures by  
\begin{align*}
&\int_{K_n}dk=1,&
&dx =\prod_{1\leq i<j\leq n}d_Fx_{i,j},&
&du =\prod_{1\leq j<i\leq n}d_Fu_{i,j},&
&da =\prod_{i=1}^n\frac{2\rc_F\,da_i}{a_i}
\end{align*}
with $x=(x_{i,j})\in N_n$, $u=(u_{i,j})\in U_n$ and 
$a=\diag (a_1,a_2,\cdots ,a_n)\in A_n$. 
When $n=1$, we understand $N_1=U_1=\{1\}$ and 
\begin{align*}
&\int_{N_1}f(x)\,dx=\int_{U_1}f(u)\,du=f(1)
\end{align*}
for a function $f$ on $\{1\}$. 
We normalize the Haar measure $dg$ on $G_n$ so that 
\begin{align*}
\int_{G_n}f(g)\,dg
&=\int_{K_n}\int_{U_n}\int_{A_n}f(auk)\, da\,du\, dk
=\int_{A_n}\int_{U_n}\int_{K_n}f(kua)\,dk\, du\, da
\end{align*}
for any integrable function $f$ on $G_n$. 
We normalize the right $G_n$-invariant measure 
$dg$ on $N_n\backslash G_n$ so that 
\begin{align}
\label{eqn:quot_GN_measure}
&\int_{G_n}f(g)\,dg
=\int_{N_n\backslash G_n}\left(\int_{N_n}f(xg)\,dx\right) dg
\end{align}
for any integrable function $f$ on $G_n$. 
We normalize the right $G_n$-invariant measure 
$dg$ on $Z_nN_n\backslash G_n$ so that 
\begin{align}
\label{eqn:quot_GN_measure}
&\int_{N_n\backslash G_n}f(g)\,dg
=\int_{Z_nN_n\backslash G_n}
\left(\int_{G_1}f(hg)\,dh\right) dg
\end{align}
for any integrable function $f$ on $N_n\backslash G_n$.

\subsection{Principal series representations of $G_n$}
\label{subsec:Cn_def_ps}

Following Jacquet \cite{Jacquet_001}, 
we will define principal series representations of $G_n$ 
as representations induced from characters of 
the lower triangular Borel subgroup $U_nM_n$ of $G_n$ 
in this paper.

Let $d=(d_1,d_2,\cdots, d_n)\in \bZ^n$ and 
$\nu =(\nu_1,\nu_2,\cdots ,\nu_n)\in \bC^n$. 
For $l\in \bZ$ and $t\in F^\times$, 
we set $\chi_{l}(t)=(t/|t|)^{l}$. 
We define characters $\chi_{d}$ and $\eta_{\nu}$ of $M_n$ by 
\begin{align*}
&\chi_{d}(m)
=\prod_{i=1}^n\chi_{d_i}(m_i)=\prod_{i=1}^n
\left(\frac{m_i}{|m_i|}\right)^{\!d_i},&
&\eta_{\nu}(m)
=\prod_{i=1}^n|m_i|_F^{\nu_i}
=\prod_{i=1}^n|m_i|^{\nu_i\rc_F}
\end{align*}
for $m=\diag (m_1,m_2,\cdots ,m_n)\in M_n$. Let 
$\rho_n=(\rho_{n,1},\rho_{n,2},\cdots ,\rho_{n,n})\in \bQ^n$ 
with 
\begin{align*}
&\rho_{n,i}=\tfrac{n+1}{2}-i&&(1\leq i\leq n). 
\end{align*}
Let $I(d,\nu )$ be the subspace of $C^\infty (G_n)$ consisting of 
all functions $f$ such that  
\begin{align}
\label{eqn:def_ps}
&f(umg)=\chi_{d}(m)\eta_{\nu -\rho_n}(m)f(g)&
&(u\in U_n,\ m\in M_n,\ g\in G_n), 
\end{align}
on which $G_n$ acts by the right translation $\Pi_{d,\nu}=R$. 
We equip $I(d,\nu )$ with the usual Fr\'{e}chet topology. 
We call $(\Pi_{d,\nu},I(d,\nu ))$ 
a (smooth) principal series representation of $G_n$. 
We denote by $I(d,\nu )_{K_n}$ the subspace of $I(d,\nu )$ 
consisting of all $K_n$-finite vectors. 
When $F=\bR$, we note that 
\begin{align}
\label{eqn:ps_isom_Rshift}
&\chi_{d+l}=\chi_{d}&
&I(d+l,\nu) =I(d,\nu ) &
&(l\in 2\bZ^n). 
\end{align} 
When $I(d,\nu )$ is irreducible, 
for any element $\sigma$ of 
the symmetric group $\mathfrak{S}_n$ of degree $n$, 
we have  
\begin{align}
\label{eqn:ps_isom_weyl}
&I(d,\nu ) \simeq I((d_{\sigma (1)},d_{\sigma (2)},\cdots ,
d_{\sigma (n)}),\, 
(\nu_{\sigma (1)},\nu_{\sigma (2)},\cdots ,
\nu_{\sigma (n)}))
\end{align} 
as representations of $G_n$ 
(\textit{cf.} \cite[Corollary 2.8]{Speh_Vogan_001}).

Let $I(d)$ be the space of smooth functions $f$ on $K_n$ satisfying 
\begin{align*}
&f(mk)=\chi_d(m)f(k)&
&(m\in M_n\cap K_n,\ k\in K_n), 
\end{align*}
and we equip this space with the usual Fr\'{e}chet topology. 
Because of the decomposition $G_n=U_nA_nK_n$ and (\ref{eqn:def_ps}), 
we can identify the space $I(d,\nu )$ with $I(d)$ via the restriction map 
$I(d,\nu )\ni f\mapsto f|_{K_n}\in I(d)$ to $K_n$. 
The inverse map $I(d)\ni f\mapsto f_\nu \in I(d,\nu )$ 
of the restriction map is given by 
\begin{align}
\label{eqn:def_st_section}
&f_\nu (uak)=\eta_{\nu -\rho_n}(a)f(k)&
&(u\in U_n,\ a\in A_n,\ k\in K_n).
\end{align}
We regard $I(d)$ as a $G_n$-module via this identification, and 
we denote the action of $G_n$ on $I(d)$ corresponding to 
$\Pi_{d,\nu}$ by $\Pi_{\nu}$, that is,  
\begin{align*}
&(\Pi_{\nu}(g)f)(k)=f_\nu (kg)&
&(g\in G_n,\ k\in K_n,\ f\in I(d)).
\end{align*}
Here we note that $\Pi_{\nu}|_{K_n}$ is the right translation 
and does not depend on $\nu$. 
We denote by $I(d)_{K_n}$ the subspace of $I(d)$ 
consisting of all $K_n$-finite vectors. 
For $f\in I(d)$, we call a map 
$\bC^n\ni \nu \mapsto f_\nu \in C^\infty (G_n)$ 
defined by (\ref{eqn:def_st_section}) 
the standard section corresponding to $f$.

\begin{rem}
For the study of automorphic forms such as the Eisenstein series, 
it is convenient to realize 
principal series representations of $G_n$ 
as representations $(\Pi_{B_n,d,\nu},I_{B_n}(d,\nu ))$ induced from 
characters of the upper triangular Borel subgroup $B_n=N_nM_n$, 
that is, $I_{B_n}(d,\nu )$ is the subspace of $C^\infty (G_n)$ 
consisting of all functions $f$ such that  
\begin{align*}
&f(xmg)=\chi_{d}(m)\eta_{\nu +\rho_n}(m)f(g)&
&(x\in N_n,\ m\in M_n,\ g\in G_n), 
\end{align*}
and the action $\Pi_{B_n,d,\nu}$ of $G_n$ is the right translation $R$. 
The results in this paper can be translated into this realization 
via the $G_n$-isomorphism 
\begin{align*}
I_{B_n}(d,\nu )\ni f\mapsto f^{w_n}\in 
I((d_n,d_{n-1},\cdots ,d_1),\, (\nu_n,\nu_{n-1},\cdots ,\nu_1))
\end{align*}
with $f^{w_n}(g)=f(w_ng)$\ \ $(g\in G_n)$. 
Here $w_n$ is the anti-diagonal matrix of size $n$ 
with $1$ at all anti-diagonal entries. 
\end{rem}

\subsection{Whittaker functions} 

Let $\varepsilon \in \{\pm 1\}$, 
and let $\psi_{\varepsilon}$ be the standard character of $F$ defined 
in \S \ref{subsec:notation}. 
Let $\psi_{\varepsilon ,n}$ be a character of $N_n$ defined by 
\begin{align*}
&\psi_{\varepsilon,n} (x)=\psi_\varepsilon (x_{1,2}+x_{2,3}+\cdots +x_{n-1,n})&
&(x=(x_{i,j})\in N_n ).
\end{align*}
When $n=1$, we understand that $\psi_{\varepsilon,1}$ 
is the trivial character of $N_1=\{1\}$.

Let $d\in \bZ^n$ and 
$\nu =(\nu_1,\nu_2,\cdots ,\nu_n)\in \bC^n$. 
A $\psi_\varepsilon $-form on $I(d,\nu )$ is a continuous $\bC$-linear form 
$\cT \colon I(d,\nu ) \to \bC$ satisfying 
\begin{align*}
&\cT (\Pi_{d,\nu}(x)f)=\psi_{\varepsilon ,n} (x)\cT (f)&
&(x\in N_n,\ f\in I(d,\nu )).
\end{align*}
Kostant \cite{Kostant_001} shows that the space of $\psi_\varepsilon $-forms 
on $I(d,\nu )$ is one dimensional. 
Let us recall the construction of nonzero $\psi_\varepsilon $-forms on 
principal series representations of $G_n$, 
which are called the Jacquet integrals. 
If $\nu$ satisfies 
\begin{align}
\label{eqn:cdn_conv_Jac}
&\mathrm{Re}(\nu_{i+1}-\nu_i)>0&
&(1\leq i\leq n-1), 
\end{align}
we define the Jacquet integral 
$\cJ_{\varepsilon}\colon I(d,\nu )\to \bC$ by 
the convergent integral 
\begin{align*}
&\cJ_{\varepsilon }(f)=
\int_{N_n}f(x)\psi_{-\varepsilon,n}(x)\,dx&
&(f\in I(d,\nu )).
\end{align*}
When $n=1$, we understand $\cJ_{\varepsilon }(f)=f(1)$ 
($f\in I(d,\nu )$). 
For $\nu \in \bC$ satisfying (\ref{eqn:cdn_conv_Jac}), we set 
$\cJ_{\varepsilon }^{(d,\nu )}(f)
=\cJ_{\varepsilon }(f_\nu )$ ($f\in I(d)$), 
where $f_{\nu}$ is the standard section corresponding to $f$. 
By \cite[Theorem 15.4.1]{Wallach_003}, we know that 
$\cJ_{\varepsilon }^{(d,\nu )}(f)$ has the holomorphic continuation 
to whole $\nu \in \bC^n$ for every $f\in I(d)$, 
and $\bC^n\times I(d)\ni (\nu ,f)\mapsto 
\cJ_{\varepsilon }^{(d,\nu )}(f)\in \bC$ is continuous. 
Furthermore, this extends $\cJ_{\varepsilon }^{(d,\nu )}$ 
to all $\nu \in \bC^n$ as a nonzero continuous $\bC$-linear form 
on $I(d)$ satisfying 
\begin{align*}
&\cJ_{\varepsilon }^{(d,\nu )} (\Pi_{\nu}(x)f)
=\psi_{\varepsilon ,n} (x)\cJ_{\varepsilon }^{(d,\nu )}(f)&
&(x\in N_n,\ f\in I(d)).
\end{align*} 
We extends the Jacquet integral 
$\cJ_{\varepsilon}\colon I(d,\nu )\to \bC$ 
to whole $\nu \in \bC^n$ by 
\begin{align*}
\cJ_{\varepsilon }(f)
=\cJ_{\varepsilon }^{(d,\nu )}(f|_{K_n})&
&(f\in I(d,\nu ))
\end{align*}
which is a nonzero $\psi_\varepsilon $-form on $I(d,\nu )$. 
We set 
\begin{align}
\label{eqn:def_JacWhittaker}
&\rW_{\varepsilon}(f)(g)=\cJ_{\varepsilon}(\Pi_{d,\nu}(g)f)&
&(f\in I(d,\nu ),\ g\in G_n). 
\end{align}
For $f\in I(d,\nu )$, 
$\rW_{\varepsilon}(f)$ is called 
a Whittaker function for $(\Pi_{d,\nu},\psi_\varepsilon )$, and satisfies 
\begin{align}
\label{eqn:Wh_rel_inv}
&\rW_{\varepsilon}(f)(xg)
=\psi_{\varepsilon ,n}(x)\rW_{\varepsilon}(f)(g)&
&(x\in N_n,\ g\in G_n).
\end{align}
We note that 
$\rW_{\varepsilon}(f_\nu )(g)
=\cJ_{\varepsilon}^{(d,\nu )}(\Pi_{\nu}(g)f)$ 
is an entire function of $\nu$ for $g\in G_n$ and 
the standard section $f_{\nu}$ corresponding to $f\in I(d)_{K_n}$. 
Let 
\[
\cW (\Pi_{d,\nu},\psi_\varepsilon )=\{\rW_{\varepsilon}(f)
\mid f\in I(d,\nu )\}. 
\]
When $\Pi_{d,\nu}$ is irreducible, 
this is a Whittaker model of $\Pi_{d,\nu}$.

\subsection{The Gelfand--Tsetlin type basis}
\label{subsec:GT_basis}

In this subsection, 
we introduce a Gelfand--Tsetlin type basis of 
an irreducible holomorphic finite dimensional representation of 
$\mathrm{GL}(n,\bC )$. 
Let $\g \gl (n,\bC )=\mathrm{M}_n(\bC )$ be 
the associated Lie algebra of $\mathrm{GL}(n,\bC )$. 
For $1\leq i<j\leq n$, we denote by $E_{i,j}$ 
the matrix unit in $\g \gl (n,\bC )$ 
with $1$ at the $(i,j)$-th entry and $0$ at other entries. 
We define the set $\Lambda_n$ of dominant weights by 
\[
\Lambda_n=\{\lambda =(\lambda_1,\lambda_2,\cdots ,\lambda_n)\in \bZ^n 
\mid \lambda_1\geq \lambda_2\geq \cdots \geq \lambda_n\}. 
\]
Let $(\tau_\lambda ,V_\lambda )$ be an irreducible holomorphic 
finite dimensional representation 
of $\mathrm{GL}(n,\bC)$ with highest weight 
$\lambda =(\lambda_1,\lambda_2,\cdots ,\lambda_n)\in \Lambda_n$, 
and we fix a $\rU (n)$-invariant hermitian inner product 
$\langle \cdot ,\cdot \rangle$ on $V_\lambda$. 
Then we have  
\[
\dim V_{\lambda}=\prod_{1\leq i<j\leq n}
\frac{\lambda_i-\lambda_j+j-i}{j-i}
\]
by Weyl's dimension formula \cite[Theorem 4.48]{Knapp_002}.

Let us recall the orthonormal basis on $V_\lambda$, 
which is constructed by Gelfand and Tsetlin \cite{Gelfand_Tsetlin_001} 
(see Zhelobenko \cite{Zhelobenko_001} for a detailed proof). 
We call  
\begin{align*}
&M=({m}_{i,j})_{1\leq i\leq j\leq n}=\left(\begin{array}{c}
{m}_{1,n}\quad {m}_{2,n}\quad \quad  \cdots \quad \quad {m}_{n,n}\\
{m}_{1,n-1}\ \ \cdots \ \ {m}_{n-1,n-1}\\
\cdots \ \ \cdots \ \ \cdots \\
{m}_{1,2}\ \ {m}_{2,2}\\
{m}_{1,1}
\end{array}\right)&
&(m_{i,j}\in \bZ )
\end{align*}
a integral triangular array of size $n$, and call $m_{i,j}$ 
the $(i,j)$-th entry of $M$. 
We denote by $\rG (\lambda)$ the set of 
integral triangular arrays $M=({m}_{i,j})_{1\leq i\leq j\leq n}$ of size $n$
such that 
\begin{align*}
&{m}_{i,n}=\lambda_i\quad  (1\leq i\leq n),&
&{m}_{j,k}\geq {m}_{j,k-1}\geq {m}_{j+1,k}\quad  (1\leq j<k\leq n).
\end{align*}
For $M=({m}_{i,j})_{1\leq i\leq j\leq n}\in \rG (\lambda )$, 
we define $\gamma^M=(\gamma^M_1,\gamma^M_2,\cdots ,\gamma^M_n)$ by 
\begin{align}
\label{eqn:def_wt_M}
&\gamma_{j}^M=\sum_{i=1}^j{m}_{i,j}-\sum_{i=1}^{j-1}{m}_{i,j-1}&
&(1\leq j\leq n).
\end{align}
We call $\gamma^M$ the weight of $M$. 
Gelfand and Tsetlin construct 
an orthonormal basis $\{\zeta_M\}_{M\in \rG (\lambda )}$ 
of $V_\lambda$ with the following formulas of 
$\g \gl (n,\bC)$-actions:
\begin{align}
\label{eqn:GT_act_wt}
&\tau_\lambda (E_{k,k})\zeta_{M}
=\gamma_{k}^M \zeta_{M},\\
\label{eqn:GT_act+}
&\tau_\lambda (E_{j,j+1})\zeta_{M}=
\underset{M+\Delta_{i,j}\in \rG (\lambda)}{\sum_{1\leq i\leq j}}
\tilde{\ra}_{i,j}^+(M)\zeta_{M+\Delta_{i,j}},\\
\label{eqn:GT_act-}
&\tau_\lambda (E_{j+1,j})\zeta_{M}=
\underset{M-\Delta_{i,j}\in \rG (\lambda)}{\sum_{1\leq i\leq j}}
\tilde{\ra}_{i,j}^-(M)\zeta_{M-\Delta_{i,j}}
\end{align}
for $1\leq k\leq n$, $1\leq j\leq n-1$ and 
$M=(m_{i,j})_{1\leq i<j\leq n}\in \rG (\lambda )$, where 
$\Delta_{i,j}$ is the integral triangular array of size $n$ 
with $1$ at the $(i,j)$-th entry and $0$ at other entries, and 
\begin{align*}
\tilde{\ra}_{i,j}^+(M)=&\left|
\frac{\bigl(\prod_{h=1}^{j+1}(m_{h,j+1}-m_{i,j}-h+i)\bigr) 
\prod_{h=1}^{j-1}(m_{h,j-1}-m_{i,j}-h+i-1)}
{\prod_{1\leq h\leq j,\ h\neq i}
(m_{h,j}-m_{i,j}-h+i)(m_{h,j}-m_{i,j}-h+i-1)}
\right|^{\frac{1}{2}}\!,\\
\tilde{\ra}_{i,j}^-(M)=&\left|
\frac{\bigl(\prod_{h=1}^{j+1}(m_{h,j+1}-m_{i,j}-h+i+1)\bigr) 
\prod_{h=1}^{j-1}(m_{h,j-1}-m_{i,j}-h+i)}
{\prod_{1\leq h\leq j,\ h\neq i}
(m_{h,j}-m_{i,j}-h+i)(m_{h,j}-m_{i,j}-h+i+1)}
\right|^{\frac{1}{2}}\!. 
\end{align*}

We denote by $H(\lambda )$ a unique element of  
$\rG (\lambda )$ whose weight is  $\lambda $, 
that is, 
\begin{align}
\label{eqn:def_Hlambda}
&H(\lambda )=(h_{i,j})_{1\leq i\leq j\leq n}\in 
\rG (\lambda )\qquad  \text{with}\quad 
h_{i,j}=\lambda_i. 
\end{align}
Then $\zeta_{H(\lambda )}$ is 
a highest vector in $V_\lambda$, that is, 
\begin{align*}
&\tau_\lambda (E_{i,i})\zeta_{H(\lambda )}
=\lambda_i\zeta_{H(\lambda )} 
\quad (1\leq i\leq n),\ \ &
&\tau_\lambda (E_{j,k})\zeta_{H(\lambda )}=0 
\quad (1\leq j<k\leq n).
\end{align*}
There is a $\bQ$-rational structure of $V_\lambda $ associated to 
the highest weight vector $\zeta_{H(\lambda )}$. 
It comes from the natural $\bQ$-rational structure of 
a tensor power of the standard representation of $\mathrm{GL}(n,\bC )$. 
We fix an embedding of $V_\lambda$ into a tensor power 
of the standard representation of $\mathrm{GL}(n,\bC )$ so that 
the image of $\zeta_{H(\lambda )}$ is $\bQ$-rational, 
and give a $\bQ$-rational structure of $V_\lambda$ via this embedding. 

Let us construct a Gelfand--Tsetlin type $\bQ$-rational basis 
of $V_\lambda$. 
We set 
\begin{align}
\label{eqn:def_Qbasis}
&\xi_M=\sqrt{\rr (M)}\zeta_M&
&(M=(m_{i,j})_{1\leq i<j\leq n}\in \rG (\lambda ))
\end{align}
with the rational constant 
\begin{equation}
\label{eqn:def_rM}
\rr (M)=\prod_{1\leq i\leq j<k\leq n}
\frac{(m_{i,k}-m_{j,k-1}-i+j)!(m_{i,k-1}-m_{j+1,k}-i+j)!}
{(m_{i,k-1}-m_{j,k-1}-i+j)!(m_{i,k}-m_{j+1,k}-i+j)!}.
\end{equation}
Then $\{\xi_M\}_{M\in \rG (\lambda)}$ is an orthogonal basis 
of $V_\lambda$ such that $\langle \xi_M,\xi_M\rangle =\rr (M)$ 
($M\in \rG (\lambda)$). 
For a integral triangular array 
$M=({m}_{i,j})_{1\leq i\leq j\leq n}$, 
we define the dual triangular array  
$M^\vee=({m}_{i,j}^\vee )_{1\leq i\leq j\leq n}$ of $M$ by 
${m}_{i,j}^\vee =-{m}_{j+1-i,j}$. 
The formulas corresponding to (\ref{eqn:GT_act_wt}), (\ref{eqn:GT_act+}) and 
(\ref{eqn:GT_act-}) are given respectively by 
\begin{align}
\label{eqn:xi_act_wt}
&\tau_\lambda (E_{k,k})\xi_{M}=\gamma_{k}^M \xi_{M},\\
\label{eqn:xi_act_+}
&\tau_\lambda (E_{j,j+1})\xi_M
=\underset{M+\Delta_{i,j}\in \rG (\lambda)}{\sum_{1\leq i\leq j}}
\ra_{i,j}(M)\xi_{M+\Delta_{i,j}},\\
\label{eqn:xi_act_-}
&\tau_\lambda (E_{j+1,j})\xi_M=
\underset{M+\Delta_{i,j}^\vee \in \rG (\lambda)}{\sum_{1\leq i\leq j}}
\ra_{i,j}(M^\vee)\xi_{M+\Delta_{i,j}^\vee }
\end{align}
for $1\leq k\leq n$, $1\leq j\leq n-1$ and 
$M=(m_{i,j})_{1\leq i<j\leq n}\in \rG (\lambda )$, 
where $\ra_{i,j}(M)$ is a rational number given by 
\begin{align*}
\ra_{i,j}(M)=
\frac{\prod_{h=1}^{i}(m_{h,j+1}-m_{i,j}-h+i)}
{\prod_{h=1}^{i-1}(m_{h,j}-m_{i,j}-h+i)}
\left(\prod_{h=2}^{i}\frac{m_{h-1,j-1}-m_{i,j}-h+i}
{m_{h-1,j}-m_{i,j}-h+i}\right).
\end{align*}
By these formulas and the equality $\xi_{H(\lambda )}=\zeta_{H(\lambda )}$, 
we know that 
$\{\xi_M\}_{M\in \rG (\lambda )}$ is a $\bQ$-rational basis 
of $V_\lambda$.

Until the end of this subsection, 
we assume $n>1$. Let  
\begin{align*}
&\Xi^+ (\lambda )=\{
\mu =(\mu_1,\mu_2,\cdots , \mu_{n-1})\in \Lambda_{n-1}\mid 
\lambda_i\geq \mu_i\geq \lambda_{i+1}\ \ (1\leq i\leq n-1)\}. 
\end{align*}
We regard $\mathrm{GL}_{n-1}(\bC )$ 
as a subgroup of $\mathrm{GL}_n(\bC )$ 
via the embedding 
\begin{align}
\label{eqn:embed_Gn-1_Gn}
&\iota_n\colon \mathrm{GL}(n-1,\bC )\ni g\mapsto 
\begin{pmatrix}
g&O_{n-1,1}\\ 
O_{1,n-1}&1
\end{pmatrix}\in \mathrm{GL}(n,\bC ).
\end{align}
We set $\widehat{M}=(m_{i,j})_{1\leq i\leq j\leq n-1}$ 
for $M=(m_{i,j})_{1\leq i\leq j\leq n}\in \rG (\lambda )$. 
By the construction of $\{\xi_{M}\}_{M\in \rG (\lambda )}$, 
we know that $V_{\lambda}$ has the irreducible decomposition 
\begin{align}
\label{eqn:decomp_n_n-1}
&V_{\lambda}=\bigoplus_{\mu \in \Xi^+ (\lambda )}V_{\lambda,\mu},&
&V_{\lambda,\mu}=\bigoplus_{M\in \rG (\lambda ;\mu )}\bC \xi_M 
\simeq V_\mu ,  
\end{align}
as a $\mathrm{GL}(n-1,\bC )$-module, where  
\[
\rG (\lambda ;\mu )=\{M\in \rG (\lambda )\mid 
\widehat{M}\in \rG (\mu )\}.
\]
Let $\mu \in \Xi^+ (\lambda )$. 
We define a $\bC$-linear map 
$\rR^{\lambda}_\mu \colon V_{\lambda }\to V_{\mu}$ 
by 
\begin{align}
\label{eqn:def_projR}
&\rR^{\lambda}_\mu (\xi_M)=\left\{\begin{array}{ll}
\xi_{\widehat{M}}&\text{if}\ M\in \rG (\lambda ;\mu ),\\
0&\text{otherwise}
\end{array}\right.&
&(M\in \rG (\lambda )).
\end{align}
By the formulas (\ref{eqn:xi_act_wt}), 
(\ref{eqn:xi_act_+}) and (\ref{eqn:xi_act_-}), 
we know that $\rR^{\lambda}_\mu$ is 
a surjective $\mathrm{GL}(n-1,\bC )$-homomorphism. 
For $M\in \rG (\mu)$, we denote by $M[\lambda ]$ 
the element of $\rG (\lambda ;\mu )$ characterized by 
$\widehat{M[\lambda ]}=M$, that is, 
\begin{align}
\label{eqn:def_Mlambda}
&M[\lambda ]=\begin{pmatrix}\lambda \\ M\end{pmatrix}\in \rG (\lambda ;\mu ). 
\end{align}
Then we have 
$\displaystyle H(\mu )[\lambda ]=
\begin{pmatrix}\lambda \\ H(\mu )\end{pmatrix}$, 
and $\xi_{H(\mu )[\lambda ]}$ is the highest weight vector in 
the $\mathrm{GL} (n-1,\bC)$-module $V_{\lambda,\mu}$. 
For later use, we prepare the following lemma. 

\begin{lem}
\label{lem:unitary_isom_lambda_mu}
Retain the notation. A $\bC$-linear map 
\begin{align*}
&V_\mu \ni \zeta_M\mapsto \zeta_{M[\lambda ]}\in V_{\lambda, \mu}&
&(M\in \rG (\mu )) 
\end{align*}
is a $\mathrm{GL} (n-1,\bC)$-isomorphism which preserves 
the fixed inner products $\langle \cdot ,\cdot \rangle $. 
\end{lem}
\begin{proof}
The assertion follows from 
(\ref{eqn:GT_act_wt}), (\ref{eqn:GT_act+}) and 
(\ref{eqn:GT_act-}). 
\end{proof}

\subsection{Complex conjugate representations}
\label{subsec:com_conj_rep}

For a finite dimensional representation $(\tau,V_\tau )$ of 
$\mathrm{GL}(n,\bC )$, 
we define the complex conjugate representation 
$(\overline{\tau},\overline{V_\tau})$ 
of $\tau$ as follows:
\begin{itemize}
\item Let $\overline{V_\tau}$ be a set with a fixed bijective map 
$V_\tau \ni v\mapsto \overline{v}\in \overline{V_\tau}$. 
We regard $\overline{V_\tau}$ as a $\bC$-vector space via 
the following addition and 
scalar multiplication:  
\begin{align*}
&\overline{v_1}+\overline{v_2}=\overline{v_1+v_2}\quad 
(v_1,v_2\in V_\tau ),&
&c\overline{v}=\overline{\overline{c}v}\quad 
(c\in \bC,\ v\in V_\tau ),
\end{align*}
where $\overline{c}$ is the complex conjugate of $c$. 

\item The action $\overline{\tau }$ is defined by 
$\overline{\tau }(g)\overline{v}
=\overline{\tau (g)v}$\quad 
($g\in \mathrm{GL}(n,\bC),\ v\in V_{\tau}$).
\end{itemize}
By definition, the following assertions hold 
for finite dimensional representations $(\tau,V_\tau )$ and 
$(\tau',V_{\tau'})$ of $\mathrm{GL}(n,\bC )$: 
\begin{itemize}
\item The complex conjugate representation 
$(\overline{\overline{\tau}},\overline{\overline{V_\tau}})$ 
of $\overline{\tau}$ is naturally identified with 
$(\tau,V_\tau )$ via the correspondence 
$\overline{\overline{v}}
\leftrightarrow v$ \ ($v\in V_\tau $).

\item If $\langle \cdot ,\cdot \rangle$ is 
a $\rU (n)$-invariant hermitian inner product on $V_\tau$, 
then 
\begin{align*}
V_\tau \otimes_\bC \overline{V_\tau}
\ni v_1\otimes \overline{v_2}\mapsto 
\langle v_1,v_2\rangle \in \bC
\end{align*}
is a non-degenerate $\bC$-bilinear 
$\mathrm{U}(n)$-invariant pairing.

\item The complex conjugate representation 
$(\overline{\tau \otimes \tau'},
\overline{V_\tau \otimes_\bC V_{\tau'}})$ of $\tau \otimes \tau '$ 
is naturally identified with 
$(\overline{\tau}\otimes \overline{\tau'},
\overline{V_\tau} \otimes_\bC \overline{V_{\tau'}})$ 
via the correspondence 
\begin{align*}
&\overline{v_1\otimes v_2}\leftrightarrow 
\overline{v_1}\otimes \overline{v_2}&
&(v_1\in V_\tau,\ v_2\in V_{\tau'}).
\end{align*}

\item For any subgroup $S$ of $\mathrm{GL}(n,\bC )$, 
there is a bijective $\bC$-antilinear map 
\[
\Hom_{S}(V_\tau ,V_{\tau'})\ni \Psi\mapsto 
\overline{\Psi}\in \Hom_{S}
(\overline{V_\tau} ,\overline{V_{\tau'}})
\]
defined by $\overline{\Psi}(\overline{v})
=\overline{\Psi (v)}\in \overline{V_{\tau'}}$ ($v\in V_\tau$).

\end{itemize}

Let $\lambda =(\lambda_1,\lambda_2,\cdots ,\lambda_n) \in \Lambda_n$. 
We consider the complex conjugate representation 
$(\overline{\tau_\lambda},\overline{V_\lambda})$ 
of $\tau_\lambda $. 
We denote by $\gu (n)$ the associated Lie algebra of $\rU (n)$. 
The complexification $\gu (n)_\bC =\gu (n)\otimes_{\bR}\bC$ of $\gu (n)$ 
is isomorphic to $\g \gl (n,\bC)$ via the 
correspondence $E_{i,j}^{\gu (n)}\leftrightarrow E_{i,j}$ 
($1\leq i,j\leq n$) with 
\[
E_{i,j}^{\gu (n)}
=\frac{1}{2}\{(E_{i,j}-E_{j,i})\otimes 1
-\sI (E_{i,j}+E_{j,i})\otimes \sI\}\in \gu (n)_\bC.
\]
For $1\leq i,j\leq n$ and $v\in V_\lambda$, we have 
\begin{align}
\label{eqn:conconj_Liealg_rep}
&\tau_{\lambda }(E_{i,j}^{\gu (n)})v
=\tau_{\lambda }(E_{i,j})v,&
&\overline{\tau_{\lambda }}(E_{i,j}^{\gu (n)})\overline{v}
=-\overline{\tau_{\lambda }(E_{j,i})v}.
\end{align}
By the pairing $V_\lambda \otimes_\bC \overline{V_\lambda}
\ni v_1\otimes \overline{v_2}\mapsto 
\langle v_1,v_2\rangle \in \bC$, we can identified 
$(\overline{\tau_\lambda},\overline{V_\lambda})$ 
with the contragredient 
representation $(\tau_\lambda^\vee ,V_{\lambda}^\vee )$ 
of $\tau_\lambda $ 
as a $\mathrm{U} (n)$-module. 
Let $\lambda^\vee =(-\lambda_n,-\lambda_{n-1},\cdots,-\lambda_1)
\in \Lambda_n$. 
Since $V_{\lambda}^\vee \simeq V_{\lambda^\vee}$ 
as $\mathrm{GL}(n,\bC )$-modules, we have 
$\overline{V_\lambda}\simeq V_{\lambda^\vee}$ 
as $\mathrm{U} (n)$-modules. 
In fact, by (\ref{eqn:xi_act_wt}), (\ref{eqn:xi_act_+}), (\ref{eqn:xi_act_-}) 
and (\ref{eqn:conconj_Liealg_rep}), 
we can confirm that the $\bC$-linear map 
\begin{align*}
&\overline{V_\lambda}\ni \overline{\xi_{M}}\mapsto 
(-1)^{\sum_{1\leq i\leq j\leq n}m_{i,j}}\, \xi_{M^\vee}
\in V_{\lambda^\vee}& 
&(M=(m_{i,j})_{1\leq i\leq j\leq n}\in \rG (\lambda ))
\end{align*} 
is a $\mathrm{U} (n)$-isomorphism. 
Via this isomorphism, we derive 
the $\bQ$-rational structure of $\overline{V_\lambda}$ 
from that of $V_{\lambda^\vee}$. Then 
$\{\overline{\xi_M}\}_{M\in \rG (\lambda )}$ 
is a $\bQ$-rational basis of $\overline{V_\lambda}$. 

\begin{rem}
We note that $\{E_{i,j}-E_{j,i}\}_{1\leq i<j\leq n}$ forms a basis of 
the associated Lie algebra $\go (n)$ of $\rO (n)$. 
By (\ref{eqn:xi_act_+}), (\ref{eqn:xi_act_-}) 
and (\ref{eqn:conconj_Liealg_rep}), 
we know that 
\begin{align*}
&\overline{V_\lambda}\ni \overline{\xi_{M}}\mapsto 
\xi_{M}
\in V_{\lambda}& 
&(M\in \rG (\lambda ))
\end{align*} 
defines a $\bQ$-rational $\rO (n)$-isomorphism. 
\end{rem}

\subsection{The minimal $K_n$- and $K_{n'}$-types}
\label{subsec:minKtype}

We define a subset $\Lambda_{n,F}$ of $\Lambda_n$ by 
$\Lambda_{n,\bR }=\Lambda_n\cap \{0,1\}^n$ and 
$\Lambda_{n,\bC }=\Lambda_n$. 
In \S \ref{subsec:lem_Krep}, 
we study the $\mathrm{O}(n)$-module structure of $V_\lambda$ 
for $\lambda \in \Lambda_{n,\bR}$, 
and prove the following lemma. 

\begin{lem}
\label{lem:rep_Kn_irred}
Let $\lambda \in \Lambda_{n,F}$. Then 
$V_\lambda$ is an irreducible $K_n$-module. 
Moreover, for any $\lambda' \in \Lambda_{n,F}$ such that 
$\lambda'\neq \lambda $, we have 
$V_\lambda \not\simeq V_{\lambda'}$ as $K_n$-modules. 
\end{lem}

Let $(\Pi_{d,\nu},I(d,\nu ))$ be a principal series representations of $G_n$ 
with 
\begin{equation*}
\begin{aligned}
&d=(d_1,d_2,\cdots ,d_n)\in \bZ^n, \qquad &
&\nu =(\nu_1,\nu_2,\cdots ,\nu_n)\in \bC^n
\end{aligned}
\end{equation*} 
such that $d\in \Lambda_{n,F}$. 
By the formula (\ref{eqn:xi_act_wt}) and 
the Frobenius reciprocity law \cite[Theorem 1.14]{Knapp_002}, we know that 
$\tau_d|_{K_n}$ is the minimal $K_n$-type of $\Pi_{d,\nu}$, and 
$\Hom_{K_n}(V_d,I(d,\nu ))$ is 1 dimensional. 
Let $\rf_{d,\nu}\colon V_d\to I(d,\nu )$ be 
the $K_n$-homomorphism normalized by 
$\rf_{d,\nu}(\xi_{H(d)})(1_n)=1$, that is, 
\begin{align}
\label{eqn:def_minKtype1}
&\rf_{d,\nu }(v)(uak)=\eta_{\nu -\rho_n}(a)
\langle \tau_d(k)v,\xi_{H(d)}\rangle 
\end{align}
for $u\in U_n$, $a\in A_n$, $k\in K_n$ and $v\in V_d$. 
Here $H(\lambda )$ ($\lambda \in \Lambda_n$) are 
defined by (\ref{eqn:def_Hlambda}). 
For $v\in V_d$, 
we note that $\rf_{d,\nu}(v)$ is the standard section 
corresponding to $\rf_d(v)\in I(d)$ defined by 
$\rf_d(v)(k)=\langle \tau_d(k)v,\xi_{H(d)}\rangle $ \ 
($k\in K_n$).

Let $(\Pi_{d',\nu'},I(d',\nu' ))$ be a principal series representations 
of $G_{n'}$ 
with 
\begin{equation*}
\begin{aligned}
&d'=(d_1',d_2',\cdots ,d_{n'}')\in \bZ^{n'}, \qquad &
&\nu' =(\nu_1',\nu_2',\cdots ,\nu_{n'}')\in \bC^{n'}
\end{aligned}
\end{equation*} 
such that $-d'\in \Lambda_{n',F}$. 
By the formula (\ref{eqn:xi_act_wt}) and the Frobenius reciprocity law 
\cite[Theorem 1.14]{Knapp_002}, 
we know that $\overline{\tau_{-d'}}|_{K_{n'}}$ is the minimal $K_{n'}$-type of 
$\Pi_{d',\nu'}$ and 
$\Hom_{K_{n'}}(\overline{V_{-d'}},I(d',\nu'))$ is 1 dimensional. 
Let $\bar{\rf}_{d',\nu'}\colon \overline{V_{-d'}}\to I(d',\nu')$ be 
the $K_{n'}$-homomorphism normalized by 
$\bar{\rf}_{d',\nu'}(\overline{\xi_{H(-d')}})(1_{n'})=1$, that is, 
\begin{align}
\label{eqn:def_minKtype2}
&\bar{\rf}_{d',\nu'}(\overline{v})(uak)
=\eta_{\nu' -\rho_{n'}}(a)
\overline{\langle \tau_{-d'}(k)v,\xi_{H(-d')}\rangle}
\end{align}
for $u\in U_{n'}$, $a\in A_{n'}$, $k\in K_{n'}$ and $v\in V_{-d'}$. 
For $v\in V_{-d'}$, we note that $\bar{\rf}_{d',\nu'}(\overline{v})$ 
is the standard section corresponding to 
$\bar{\rf}_{d'}(\overline{v})\in I(d')$ defined by 
$\bar{\rf}_{d'}(\overline{v})(k)=
\overline{\langle \tau_{-d'}(k)v,\xi_{H(-d')}\rangle}$ \ 
($k\in K_{n'})$. 

We define the archimedean $L$-factor 
for $\Pi_{d,\nu}\times \Pi_{d',\nu'}$ by 
\begin{align*}
L(s,\Pi_{d,\nu}\times \Pi_{d',\nu'})
=\prod_{i=1}^n\prod_{j=1}^{n'}
\Gamma_F(s+\nu_i+\nu_j;\,|d_i+d_j'|),
\end{align*}
where the functions $\Gamma_F(s;m)$ ($m\in \bZ$) 
are defined in \S \ref{subsec:notation}. 
Moreover, we set 
\begin{align*}
&\boldsymbol{\Gamma}_F(\nu;d)
=\prod_{1\leq i<j\leq n}
\Gamma_{F}(\nu_j-\nu_i+1;\,|d_i-d_j|), \\
&\boldsymbol{\Gamma}_F(\nu';d')
=\prod_{1\leq i<j\leq n'}
\Gamma_{F}(\nu_j'-\nu_i'+1;\,|d_i'-d_j'|).
\end{align*} 
Here $1/\boldsymbol{\Gamma}_F(\nu;d)$ and 
$1/\boldsymbol{\Gamma}_F(\nu';d')$ are both nonzero 
if $\Pi_{d,\nu}$ and $\Pi_{d',\nu'}$ are irreducible. 
This fact follows from Corollary \ref{cor:god_explicit_2} 
in \S \ref{subsec:god_explicit}. 

\subsection{Archimedean Rankin--Selberg integrals for $G_n\times G_{n-1}$}

In this subsection, we assume $n>1$. 
Let $(\Pi_{d,\nu},I(d,\nu ))$ and 
$(\Pi_{d',\nu'},I(d',\nu'))$ be principal series representations 
of $G_n$ and $G_{n-1}$, respectively, with parameters 
\begin{equation*}
\begin{aligned}
&d=(d_1,d_2,\cdots ,d_n)\in \bZ^n, \qquad &
&\nu =(\nu_1,\nu_2,\cdots ,\nu_n)\in \bC^n, \\
&d'=(d_1',d_2',\cdots ,d_{n-1}')\in \bZ^{n-1},& 
&\nu'=(\nu_1',\nu_2',\cdots ,\nu_{n-1}')\in \bC^{n-1}.
\end{aligned}
\end{equation*}
We assume $d\in \Lambda_{n,F}$ and $-d'\in \Lambda_{n-1,F}$. 
If $\Pi_{d,\nu}$ and $\Pi_{d',\nu'}$ are irreducible, 
these are not serious assumptions because of 
(\ref{eqn:ps_isom_Rshift}) and (\ref{eqn:ps_isom_weyl}). 
We take $\rf_{d,\nu }$, $\bar{\rf}_{d',\nu'}$, $\boldsymbol{\Gamma}_F(\nu;d)$, 
$\boldsymbol{\Gamma}_F(\nu';d')$ and $L(s,\Pi_{d,\nu}\times \Pi_{d',\nu'})$ 
as in \S \ref{subsec:minKtype} with $n'=n-1$.

Let $\varepsilon \in \{\pm 1\}$, 
$W\in \cW (\Pi_{d,\nu},\psi_\varepsilon )$ and 
$W'\in \cW (\Pi_{d',\nu'},\psi_{-\varepsilon})$. 
Let $s\in \bC$ such that 
$\mathrm{Re}(s)$ is sufficiently large. 
We define the archimedean Rankin--Selberg integral 
$Z(s,W,W')$ for $\Pi_{d,\nu} \times \Pi_{d',\nu'}$ by 
\begin{align}
\label{eqn:def_zeta_int1}
Z(s,W,W')=\int_{N_{n-1}\backslash G_{n-1}}
W(\iota_n(g))W'(g)|\det g|_F^{s-1/2}\,dg, 
\end{align}
where $\iota_n$ is defined by (\ref{eqn:embed_Gn-1_Gn}). 
Here we note 
\begin{align}
\label{eqn:zeta+_Kact}
&Z(s,R(\iota_n(k))W,R(k)W')
=Z(s,W,W')&
&(k\in K_{n-1}). 
\end{align}
By (\ref{eqn:zeta+_Kact}), we know that 
\begin{align}
\label{eqn:zeta_hom_1}
&v_1\otimes \overline{v_2}
\mapsto Z(s,\rW_{\varepsilon}(\rf_{d,\nu}(v_1)),
\rW_{-\varepsilon}(\bar{\rf}_{d',\nu'}(\overline{v_2})))
\end{align}
defines an element of 
$\Hom_{K_{n-1}}(V_{d}\otimes_\bC \overline{V_{-d'}},\bC_{\mathrm{triv}})$. 
Here $\rW_{\varepsilon}$ is defined by (\ref{eqn:def_JacWhittaker}), 
and $\bC_{\mathrm{triv}}=\bC$ is the trivial $K_{n-1}$-module. 
The following theorem is the first main result of this paper, 
which gives the explicit expression of 
the $K_{n-1}$-homomorphism (\ref{eqn:zeta_hom_1}). 

\begin{thm}
\label{thm:main1}
Retain the notation. 
For $v_1\in V_{d}$ and $v_2\in V_{-d'}$, we have 
\begin{align*}
&Z(s,\rW_{\varepsilon}(\rf_{d,\nu}(v_1)),
\rW_{-\varepsilon}(\bar{\rf}_{d',\nu'}(\overline{v_2})))\\
&=\frac{(-\varepsilon \sI )^{\sum_{i=1}^{n-1}(n-i)(d_i+d_i')}}
{(\dim V_{-d'})\boldsymbol{\Gamma}_F(\nu;d)
\boldsymbol{\Gamma}_F(\nu';d')}
L(s,\Pi_{d,\nu}\times \Pi_{d',\nu'})
\langle \rR_{-d'}^{d}(v_1),\,v_2\rangle 
\end{align*}
if $-d'\in \Xi^+(d)$, and 
$Z(s,\rW_{\varepsilon}(\rf_{d,\nu}(v_1)),
\rW_{-\varepsilon}(\bar{\rf}_{d',\nu'}(\overline{v_2})))=0$ 
otherwise. 
Here $\rR_{-d'}^{d}$ 
is given explicitly by (\ref{eqn:def_projR}). 
In particular, we have 
\begin{equation}
\label{eqn:thm_main1}
\begin{aligned}
&Z(s,\rW_{\varepsilon}(\rf_{d,\nu}(\xi_{H(-d')[d]})),
\rW_{-\varepsilon}(\bar{\rf}_{d',\nu'}(\overline{\xi_{H(-d')}})))\\
&=\frac{(-\varepsilon \sI )^{\sum_{i=1}^{n-1}(n-i)(d_i+d_i')}}
{(\dim V_{-d'})\boldsymbol{\Gamma}_F(\nu;d)
\boldsymbol{\Gamma}_F(\nu';d')}L(s,\Pi_{d,\nu}\times \Pi_{d',\nu'})
\end{aligned}
\end{equation}
if $-d'\in \Xi^+(d)$. 
Here $H(-d')$ and $H(-d')[d]$ are defined by (\ref{eqn:def_Hlambda}) 
and (\ref{eqn:def_Mlambda}). 
\end{thm}
\begin{rem}
Retain the notation. 
The space 
$\Hom_{K_{n-1}}(V_{d}\otimes_\bC \overline{V_{-d'}},\bC_{\mathrm{triv}})$ 
is $1$ dimensional if $-d'\in \Xi^+(d)$, and is equal to $\{0\}$ 
otherwise. 
This fact follows from Lemma \ref{lem:schur_pairing+} 
in \S \ref{subsec:lem_Krep}. 
\end{rem}

\begin{rem}
\label{rem:cohomological_1}
We set $F=\bC$. 
By \cite[Proposition 2.14 and Theorem 2.21]{Raghuram_001}, we note that 
the compatible pairs of cohomological representations 
of $G_n$ and $G_{n-1}$ in Sun \cite[\S 6]{Sun_001} can be regarded as 
pairs of some irreducible principal series 
representations $\Pi_{d,\nu}$ and $\Pi_{d',\nu'}$
with $d\in \Lambda_{n,\bC}$ and $-d'\in \Xi^+(d)$. 
Hence, Theorem \ref{thm:main1} gives another proof 
of the non-vanishing hypothesis for 
$G_n\times G_{n-1}$, 
which were originally proved by Sun \cite[Theorem C]{Sun_001}. 
\end{rem}

\begin{cor}
Retain the notation, and assume $-d'\in \Xi^+(d)$.
Then 
\[
\sum_{M\in \rG (-d')}\rr (M)^{-1}\,
\xi_{M[d]}\otimes \overline{\xi_M}
\]
is a unique $\bQ$-rational $K_{n-1}$-invariant vector in 
$V_d\otimes_\bC \overline{V_{-d'}}$ up to scalar multiple, 
and its image under 
the $K_{n-1}$-homomorphism (\ref{eqn:zeta_hom_1}) is 
\begin{align*}
&\frac{(-\varepsilon \sI )^{\sum_{i=1}^{n-1}(n-i)(d_i+d_i')}}
{\boldsymbol{\Gamma}_F(\nu;d)
\boldsymbol{\Gamma}_F(\nu';d')}L(s,\Pi_{d,\nu}\times \Pi_{d',\nu'}).
\end{align*}
Here $\rr (M)$ and $M[d]$ are defined by (\ref{eqn:def_rM}) and 
(\ref{eqn:def_Mlambda}), respectively. 
\end{cor}

This corollary follows from Theorem \ref{thm:main1} with   
Lemma \ref{lem:schur_pairing+}, 
and is similar to Corollary \ref{cor:main2}  
which gives the explicit description of 
the archimedean part of Greni\'{e}'s theorem 
\cite[Theorem 2]{Grenie_001}.

\subsection{Schwartz functions}
\label{subsec:def_schwartz}

Let $\cS (\rM_{n,n'}(F))$ be the space of Schwartz functions on 
$\rM_{n,n'}(F)$. We define $\me_{(n,n')}\in \cS (\rM_{n,n'}(F))$ by 
\begin{align}
&\me_{(n,n')} (z)=\exp \bigl(-\pi \rc_F\mathrm{Tr}({}^t\overline{z}z)\bigr)
=\left\{\begin{array}{ll}
\exp (-\pi \mathrm{Tr}({}^t\!zz))&\text{if }F=\bR,\\[1.5mm]
\exp (-2\pi \mathrm{Tr}({}^t\overline{z}z))&\text{if }F=\bC
\end{array}\right.
\end{align}
for $z\in \rM_{n,n'}(F)$. 
We denote $\me_{(n,n)}$ simply by $\me_{(n)}$. 
Let $\cS_0(\rM_{n,n'}(F))$ be the subspace of 
$\cS (\rM_{n,n'}(F))$ consisting of all 
functions $\phi$ of the form 
\begin{align*}
&\phi (z)=p(z,\overline{z})\me_{(n,n')}(z)&
&(z\in \rM_{n,n'}(F)),
\end{align*}
where $p$ is a polynomial function. 
We call elements of $\cS_0(\rM_{n,n'}(F))$ 
standard Schwartz functions on $\rM_{n,n'}(F)$.

Let $C(\rM_{n,n'}(F))$ be the space of continuous 
functions on $\rM_{n,n'}(F)$. 
We define actions of $G_n$ and $G_{n'}$ on $C(\rM_{n,n'}(F))$ by 
\begin{align*}
&(L(g)f)(z)=f(g^{-1}z),&&(R(h)f)(z)=f(zh)&
\end{align*}
for $g\in G_n$, $h\in G_{n'}$, $f\in C(\rM_{n,n'}(F))$ and 
$z\in \rM_{n,n'}(F)$. 
Since $\me_{(n,n')}$ is $K_n\times K_{n'}$-invariant, 
we note that $\cS_0 (\rM_{n,n'}(F))$ is closed under 
the action $L\boxtimes R$ of $K_n\times K_{n'}$, 
and all elements of $\cS_0 (\rM_{n,n'}(F))$ are $K_n\times K_{n'}$-finite.

Let $l\in \bN_0$, and we consider  
the representation $(\tau_{(l,\mathbf{0}_{n-1})},V_{(l,\mathbf{0}_{n-1})})$. 
Here we put $\mathbf{0}_{n-1}=(0,0,\cdots ,0)\in \Lambda_{n-1}$ if $n>1$, 
and erase $\mathbf{0}_{n-1}$ if $n=1$. 
We set 
\begin{align}
\label{eqn:def_ell_gamma}
&\ell (\gamma )=\gamma_1+\gamma_2+\cdots +\gamma_n&
&(\gamma =(\gamma_1,\gamma_2,\cdots ,\gamma_n)\in \bZ^n).
\end{align}
For $\gamma =(\gamma_1,\gamma_2,\cdots ,\gamma_n)\in \bN_0^n$ 
such that $\ell (\gamma )=l$, 
we denote by $Q(\gamma )$ a unique element of 
$\rG ((l,\mathbf{0}_{n-1}))$ whose weight 
is $\gamma $, that is, 
\begin{align}
\label{eqn:def_Qgamma}
Q(\gamma )=(q_{i,j})_{1\leq i\leq j\leq n}\qquad 
\text{with}\quad q_{i,j}=\left\{\begin{array}{ll}
\sum_{k=1}^j\gamma_k&\text{if $i=1$},\\
0&\text{otherwise}. 
\end{array}\right.
\end{align}
Then we have 
$\rG ((l,\mathbf{0}_{n-1}))=\{Q (\gamma )\mid 
\gamma \in \bN_0^n,\ 
\ell (\gamma )=l\}$. 
We define $\bC$-linear maps 
$\varphi_{1,n}^{(l)}\colon 
V_{(l,\mathbf{0}_{n-1})}\to \cS_0 (\rM_{1,n}(F))$ and 
$\overline{\varphi}_{1,n}^{(l)}\colon \overline{V_{(l,\mathbf{0}_{n-1})}}
\to \cS_0 (\rM_{1,n}(F))$ by 
\begin{align}
\label{eqn:def_varphi}
&\varphi_{1,n}^{(l)}(\xi_{Q(\gamma )}) (z)=
z_1^{\gamma_1}z_2^{\gamma_2}\cdots z_n^{\gamma_n}\me_{(1,n)}(z),\\
\label{eqn:def_varphi_bar}
&\overline{\varphi}_{1,n}^{(l)}(\overline{\xi_{Q(\gamma )}})(z)=
\overline{z_1^{\gamma_1}z_2^{\gamma_2}\cdots z_n^{\gamma_n}}\,\me_{(1,n)}(z)
\end{align}
for $z=(z_1,z_2,\cdots ,z_n)\in \rM_{1,n}(F)$ 
and $\gamma =(\gamma_1,\gamma_2,\cdots ,\gamma_n)\in \bN_0^n$ 
such that $\ell (\gamma )=l$. 
In \S \ref{subsec:standard_schwartz}, we show that 
$\varphi_{1,n}^{(l)}$ and $\overline{\varphi}_{1,n}^{(l)}$ are 
$K_n$-homomorphisms, where 
$\cS_0 (\rM_{1,n}(F))$ is regarded as a $K_n$-module via the action $R$.

\subsection{Injector}
\label{subsec:injector}

Let $\lambda =(\lambda_1,\lambda_2,\cdots ,\lambda_n)\in \Lambda_n$ 
and $l\in \bN_0$. 
In this subsection, we specify each irreducible components of 
the tensor product $V_{\lambda }\otimes_{\bC} V_{(l,\mathbf{0}_{n-1})}$. 
Let 
\[
\Xi^\circ (\lambda )=\{\lambda'=(\lambda_1',\lambda_2',\cdots ,\lambda_n')
\in \Lambda_n
\mid \lambda_1'\geq \lambda_1\geq \lambda_2'\geq \lambda_2\geq \cdots 
\geq \lambda_n'\geq \lambda_n \},
\]
and $\Xi^\circ (\lambda ;l)=\{\lambda'\in \Xi^\circ (\lambda )\mid 
\ell (\lambda'-\lambda )=l\}$. 
Then Pieri's rule \cite[Corollary 9.2.4]{Goodman_Wallach_001} 
asserts that 
$V_{\lambda }\otimes_{\bC} V_{(l,\mathbf{0}_{n-1})}$ 
has the irreducible decomposition 
\begin{align}
\label{eqn:pieri}
&V_{\lambda }\otimes_{\bC} V_{(l,\mathbf{0}_{n-1})}
\simeq \bigoplus_{\lambda' \in \Xi^\circ (\lambda ;l)}
V_{\lambda'} 
\end{align}
as $\mathrm{GL}(n,\bC )$-modules. 

For $\lambda'=(\lambda_1',\lambda_2',\cdots ,\lambda_n')
\in \Xi^\circ (\lambda )$, we set 
\begin{align}
\label{eqn:def_S0}
&\rS^\circ (\lambda' ,\lambda )=
\frac{\prod_{1\leq i\leq j\leq n}(\lambda_{i}'-\lambda_{j}-i+j)!}
{\prod_{1\leq i\leq j<n}(\lambda_{i}-\lambda_{j+1}'-i+j)!}.
\end{align}
When $n>1$, for $\mu =(\mu_1,\mu_2,\cdots ,\mu_{n-1})\in \Xi^+ (\lambda )$, 
we set 
\begin{align}
\label{eqn:def_S+}
&\rS^+(\lambda ,\mu )=\prod_{1\leq i\leq j<n}
\frac{(\lambda_{i}-\mu_{j}-i+j)!}{(\mu_{i}-\lambda_{j+1}-i+j)!}.
\end{align}
Let $\lambda'=(\lambda_1',\lambda_2',\cdots ,\lambda_n')
\in \Xi^\circ (\lambda ;l)$. 
Based on the result of Jucys \cite{Jucis_001}, 
we will construct a $\bQ$-rational $\mathrm{GL}(n,\bC)$-homomorphism 
$\rI^{\lambda,l}_{\lambda'}\colon 
V_{\lambda'}\to V_{\lambda }\otimes_{\bC}
V_{(l,\mathbf{0}_{n-1})}$ in \S \ref{subsec:CG_coeff}. 
The explicit expression of $\rI^{\lambda,l}_{\lambda'}$ is given by 
\begin{align}
\label{eqn:injector_explicit}
&\rI^{\lambda,l}_{\lambda'} (\xi_{M'})
=\sum_{M\in \rG (\lambda)}\,\sum_{P\in \rG ((l,\mathbf{0}_{n-1}))}
\rc^{M,P}_{M'}\xi_{M}\otimes \xi_{P}&
&(M'\in \rG (\lambda')), 
\end{align}
where $\rc^{M,P}_{M'}$ 
($M'\in \rG (\lambda')$, 
$M\in \rG (\lambda )$, 
$P\in \rG ((l,\mathbf{0}_{n-1}))$) are 
rational numbers determined by the following conditions, 
recursively:
\begin{itemize}
\item When $n=1$, 
we have $\rc^{\lambda_1,l}_{\lambda_1+l}=1$ for $\lambda_1,l\in \bZ$. 
\item When $n>1$, 
for $\mu'\in \Xi^+ (\lambda ')$, $M'\in \rG (\lambda';\mu')$, 
$\mu \in \Xi^+ (\lambda )$, $M\in \rG (\lambda ;\mu)$, 
$0\leq q\leq l$ and 
$P\in \rG ((l,\mathbf{0}_{n-1});(q,\mathbf{0}_{n-2}))$, 
we have 
\begin{align*}
\rc^{M,P}_{M'}=\,&\rc^{\widehat{M},\widehat{P}}_{\widehat{M'}}\,
\rS^\circ (\lambda',\lambda')\rS^\circ (\mu ,\mu )
\frac{l!}{q!}
\frac{\prod_{1\leq i\leq j<n}(\lambda_{i}-\lambda_{j+1}-i+j)!}
{\prod_{1\leq i\leq j\leq n}(\lambda_{i}'-\lambda_{j}-i+j)!}
\\
&\times \left(\prod_{1\leq i\leq j<n}
\frac{(\mu_{i}'-\mu_{j}-i+j)!(\mu_{i}'-\lambda_{j+1}'-i+j)!}
{(\mu_{i}'-\mu_{j}'-i+j)!(\mu_{i}-\lambda_{j+1}-i+j)!}
\right) \\
&\times 
\underset{\mu'\in \Xi^\circ (\alpha ),\  \alpha \in \Xi^\circ (\mu )}
{\sum_{\alpha \in \Xi^+ (\lambda )\cap \Xi^+ (\lambda')}}
\frac{(-1)^{\ell (\alpha -\mu )}\rS^\circ (\alpha ,\alpha )}
{\rS^\circ (\mu',\alpha )\rS^\circ (\alpha ,\mu )}
\frac{\rS^+(\lambda',\alpha )}{\rS^+(\lambda,\alpha )}
\end{align*}
if $\mu'\in \Xi^\circ (\mu ;q)$, and $\rc^{M,P}_{M'}=0$ otherwise.  
\end{itemize}
We define a $\rU (n)$-invariant hermitian inner product on 
$V_{\lambda}\otimes_\bC V_{(l,\mathbf{0}_{n-1})}$ by 
\begin{align*}
&\langle v_1\otimes v_1',\ v_2\otimes v_2'\rangle 
=\langle v_1,v_2\rangle \, \langle v_1',v_2'\rangle &
&(v_1,v_2\in V_{\lambda },\ v_1',v_2'\in V_{(l,\mathbf{0}_{n-1})}). 
\end{align*}
Then we have 
\begin{align}
\label{eqn:cdn_inj_3}
&\bigl\langle\, \mathrm{I}^{\lambda,l}_{\lambda'}(\xi_{H(\lambda')})
,\ \xi_{H(\lambda)}\otimes \xi_{Q(\lambda'-\lambda )}
\, \bigr\rangle 
=\rC^\circ (\lambda';\lambda ),\\
\label{eqn:cdn_inj_2}
&\bigl\langle\, \mathrm{I}^{\lambda,l}_{\lambda'}(v)
,\ \mathrm{I}^{\lambda,l}_{\lambda'}(v')
\, \bigr\rangle  
=\rb (\lambda'-\lambda )\rC^\circ (\lambda';\lambda )
\langle v,v'\rangle \hspace{10mm}(v,v'\in V_{\lambda'}),
\end{align}
where 
\begin{align}
\label{eqn:const_CG_coeff}
&\rC^\circ (\lambda';\lambda )=\prod_{1\leq i<j\leq n}
\frac{(\lambda_{i}'-\lambda_{j}'-i+j)!
(\lambda_{i}-\lambda_{j}-i+j-1)!}
{(\lambda_{i}'-\lambda_{j}-i+j)!
(\lambda_{i}-\lambda_{j}'-i+j-1)!},\\
\label{eqn:def_c_gamma}
&\rb (\gamma )
=\frac{(\gamma_1+\gamma_2+\cdots +\gamma_n)!}
{\gamma_1!\gamma_2!\cdots \gamma_n!}\hspace{10mm} 
(\gamma =(\gamma_1,\gamma_2,\cdots ,\gamma_n)\in \bN_0^n).
\end{align}

\subsection{Archimedean Rankin--Selberg integrals for $G_n\times G_n$}
\label{subsec:main2}

Let $(\Pi_{d,\nu},I(d,\nu ))$ and 
$(\Pi_{d',\nu'},I(d',\nu'))$ be principal series representations 
of $G_n$ with parameters 
\begin{equation*}
\begin{aligned}
&d=(d_1,d_2,\cdots ,d_n)\in \bZ^n, \qquad &
&\nu =(\nu_1,\nu_2,\cdots ,\nu_n)\in \bC^n, \\
&d'=(d_1',d_2',\cdots ,d_{n}')\in \bZ^{n},& 
&\nu'=(\nu_1',\nu_2',\cdots ,\nu_{n}')\in \bC^{n}.
\end{aligned}
\end{equation*}
We assume $d\in \Lambda_{n,F}$ and $-d'\in \Lambda_{n,F}$. 
If $\Pi_{d,\nu}$ and $\Pi_{d',\nu'}$ are irreducible, 
these are not serious assumptions because of 
(\ref{eqn:ps_isom_Rshift}) and (\ref{eqn:ps_isom_weyl}). 
We take $\rf_{d,\nu }$, $\bar{\rf}_{d',\nu'}$, $\boldsymbol{\Gamma}_F(\nu;d)$, 
$\boldsymbol{\Gamma}_F(\nu';d')$ and $L(s,\Pi_{d,\nu}\times \Pi_{d',\nu'})$ 
as in \S \ref{subsec:minKtype} with $n'=n$. 

Let $\varepsilon \in \{\pm 1\}$, 
$W\in \cW (\Pi_{d,\nu},\psi_\varepsilon )$, 
$W'\in \cW (\Pi_{d',\nu'},\psi_{-\varepsilon})$ 
and $\phi \in \cS (\rM_{1,n}(F))$. 
Let $s\in \bC$ such that 
$\mathrm{Re}(s)$ is sufficiently large. 
We define the archimedean Rankin--Selberg integral 
$Z(s,W,W',\phi)$ for $\Pi_{d,\nu} \times \Pi_{d',\nu'}$ by 
\begin{align}
\label{eqn:def_zeta_int2}
Z(s,W,W',\phi)=\int_{N_{n}\backslash G_{n}}
W(g)W'(g)\phi (e_ng)|\det g|_F^{s}\,dg, 
\end{align}
where we put $e_n=(O_{1,n-1},1)\in \rM_{1,n}(F)$ 
as in \S \ref{subsec:notation}. 
Here we note 
\begin{align}
\label{eqn:zeta0_symmetry}
&Z(s,W,W',\phi)=Z(s,W',W,\phi),\\
\label{eqn:zeta0_Kact}
&Z(s,R(k)W,R(k)W',R(k)\phi)=Z(s,W,W',\phi)&
&(k\in K_n). 
\end{align}
Let $l$ be a integer determined by 
\[
\left\{\begin{array}{ll}
l\in \{0,1\}\text{ and } l\equiv -\ell (d+d')\bmod 2  
&\text{if $F=\bR$ and $\ell (d+d')\leq 0$},\\[1mm] 
l\in \{0,-1\}\text{ and } l\equiv -\ell (d+d')\bmod 2 
&\text{if $F=\bR$ and $\ell (d+d')\geq 0$},\\[1mm] 
l=-\ell (d+d')&\text{if $F=\bC$},
\end{array}\right.
\]
where $\ell (\gamma)$ ($\gamma \in \bZ^n$) are defined by 
(\ref{eqn:def_ell_gamma}). 
By (\ref{eqn:zeta0_Kact}), we know that 
\begin{align}
\label{eqn:zeta_hom_21}
&v_1\otimes \overline{v_2}\otimes v_3
\mapsto 
Z(s,\rW_{\varepsilon}(\rf_{d,\nu}(v_1)),
\rW_{-\varepsilon}(\bar{\rf}_{d',\nu'}(\overline{v_2})),
\varphi_{1,n}^{(l)}(v_3))
\end{align}
defines an element of 
$\Hom_{K_{n}}(V_{d}\otimes_\bC \overline{V_{-d'}}
\otimes_\bC V_{(l,\mathbf{0}_{n-1})},\bC_{\mathrm{triv}})$ 
if $l\geq 0$, and 
\begin{align}
\label{eqn:zeta_hom_22}
&v_1\otimes \overline{v_2}\otimes \overline{v_3}
\mapsto 
Z(s,\rW_{\varepsilon}(\rf_{d,\nu}(v_1)),
\rW_{-\varepsilon}(\bar{\rf}_{d',\nu'}(\overline{v_2})),
\overline{\varphi}_{1,n}^{(-l)}(\overline{v_3}))
\end{align}
defines an element of 
$\Hom_{K_{n}}(V_{d}\otimes_\bC \overline{V_{-d'}}
\otimes_\bC \overline{V_{(-l,\mathbf{0}_{n-1})}},\bC_{\mathrm{triv}})$ 
if $l\leq 0$. 
Here $\rW_{\varepsilon}$, $\varphi_{1,n}^{(l)}$, 
$\overline{\varphi}_{1,n}^{(-l)}$ are defined by 
(\ref{eqn:def_JacWhittaker}), (\ref{eqn:def_varphi}), 
(\ref{eqn:def_varphi_bar}), respectively, and 
$\bC_{\mathrm{triv}}=\bC$ is the trivial $K_n$-module. 
The following theorem is the second main result of this paper,  
which gives the explicit expressions of 
the $K_n$-homomorphisms (\ref{eqn:zeta_hom_21}) 
and (\ref{eqn:zeta_hom_22}).


\begin{thm}
\label{thm:main2}
Retain the notation. 
\\[1mm]
\noindent 
(1) Assume $l\geq 0$. For $v_1\in V_{d}$, $v_2\in V_{-d'}$ and 
$v_3\in V_{(l,\mathbf{0}_{n-1})}$, we have 
\begin{equation*}
\begin{aligned}
&Z\bigl(s,
\rW_{\varepsilon}(\rf_{d,\nu}(v_1)),
\rW_{-\varepsilon}(\bar{\rf}_{d',\nu'}(\overline{v_2})),
\varphi_{1,n}^{(l)}(v_3)\bigr)\\
&=\frac{(-\varepsilon \sI )^{\sum_{i=1}^{n-1}(n-i)(d_i+d_i')}}
{(\dim V_{-d'})\boldsymbol{\Gamma}_F(\nu;d)
\boldsymbol{\Gamma}_F(\nu';d')}
L(s,\Pi_{d,\nu}\times \Pi_{d',\nu'})\,
\langle v_1\otimes v_3,\,\mathrm{I}^{d,l}_{-d'}(v_2)\rangle .
\end{aligned}
\end{equation*}
if $-d'\in \Xi^\circ (d)$, 
and $Z\bigl(s,
\rW_{\varepsilon}(\rf_{d,\nu}(v_1)),
\rW_{-\varepsilon}(\bar{\rf}_{d',\nu'}(\overline{v_2})),
\varphi_{1,n}^{(l)}(v_3)\bigr)=0$ otherwise. 
Here $\mathrm{I}^{d,l}_{-d'}$ 
is given explicitly by (\ref{eqn:injector_explicit}). 
In particular, if $-d'\in \Xi^\circ (d)$, we have 
\begin{equation}
\label{eqn:main2_22} 
\begin{aligned}
&Z\bigl(s,
\rW_{\varepsilon}(\rf_{d,\nu}(\xi_{H(d)})),
\rW_{-\varepsilon}(\bar{\rf}_{d',\nu'}(\overline{\xi_{H(-d')}})),
\varphi_{1,n}^{(l)}(\xi_{Q(-d-d')})\bigr)\\
&=\frac{(-\varepsilon \sI )^{\sum_{i=1}^{n-1}(n-i)(d_i+d_i')}
\rC^\circ (-d';d)}
{(\dim V_{-d'})\boldsymbol{\Gamma}_F(\nu;d)
\boldsymbol{\Gamma}_F(\nu';d')}
L(s,\Pi_{d,\nu}\times \Pi_{d',\nu'}).
\end{aligned}
\end{equation}

\noindent (2) Assume $l\leq 0$. 
For $v_1\in V_{d}$, $v_2\in V_{-d'}$ and 
$v_3\in V_{(-l,\mathbf{0}_{n-1})}$, we have 
\begin{align*}
&Z(s,\rW_{\varepsilon}(\rf_{d,\nu}(v_1)),
\rW_{-\varepsilon}(\bar{\rf}_{d',\nu'}(\overline{v_2})),
\overline{\varphi}_{1,n}^{(-l)}(\overline{v_3}))\\
&=\frac{(-\varepsilon \sI )^{\sum_{i=1}^{n-1}(n-i)(d_i+d_i')}}
{(\dim V_{d})\boldsymbol{\Gamma}_F(\nu;d)
\boldsymbol{\Gamma}_F(\nu';d')}
L(s,\Pi_{d,\nu}\times \Pi_{d',\nu'})\,
\langle \mathrm{I}^{-d',-l}_{d}(v_1),\,v_2\otimes v_3\rangle 
\end{align*}
if $d\in \Xi^\circ (-d')$, and 
$Z(s,\rW_{\varepsilon}(\rf_{d,\nu}(v_1)),
\rW_{-\varepsilon}(\bar{\rf}_{d',\nu'}(\overline{v_2})),
\overline{\varphi}_{1,n}^{(-l)}(\overline{v_3}))=0$ otherwise. 
Here $\mathrm{I}^{-d',-l}_{d}$ 
is given explicitly by (\ref{eqn:injector_explicit}). 
In particular, if $d\in \Xi^\circ (-d')$, we have 
\begin{equation}
\label{eqn:main2_12} 
\begin{aligned}
&Z(s,\rW_{\varepsilon}(\rf_{d,\nu}(\xi_{H(d)})),
\rW_{-\varepsilon}(\bar{\rf}_{d',\nu'}(\overline{\xi_{H(-d')}})),
\overline{\varphi}_{1,n}^{(-l)}(\overline{\xi_{Q(d+d')}}))\\
&=\frac{(-\varepsilon \sI )^{\sum_{i=1}^{n-1}(n-i)(d_i+d_i')}
\rC^\circ (d;-d')}{(\dim V_{d})\boldsymbol{\Gamma}_F(\nu;d)
\boldsymbol{\Gamma}_F(\nu';d')}
L(s,\Pi_{d,\nu}\times \Pi_{d',\nu'}).
\end{aligned}
\end{equation}
\end{thm}
\begin{rem}
Retain the notation. Under the assumption $l\geq 0$, 
the space 
\[
\Hom_{K_{n}}(V_{d}\otimes_\bC \overline{V_{-d'}}
\otimes_\bC V_{(l,\mathbf{0}_{n-1})},\bC_{\mathrm{triv}})
\]
is $1$ dimensional if $-d'\in \Xi^\circ (d)$, 
and is equal to $\{0\}$ otherwise. 
Under the assumption $l\leq 0$, the space 
\[
\Hom_{K_{n}}(V_{d}\otimes_\bC \overline{V_{-d'}}
\otimes_\bC \overline{V_{(-l,\mathbf{0}_{n-1})}},\bC_{\mathrm{triv}})
\]
is $1$ dimensional if $d\in \Xi^\circ (-d')$, 
and is equal to $\{0\}$ otherwise. 
These facts follow from Lemma \ref{lem:schur_pairing0} 
in \S \ref{subsec:lem_Krep}. 
\end{rem}

Let $P_{n}$ be a maximal parabolic subgroup of $G_n$ defined by 
\[
P_{n}=\{p=(p_{i,j})\in G_n\mid 
p_{n,j}=0\ (1\leq i\leq n-1)\}, 
\]
which contains 
the upper triangular Borel subgroup $B_n=N_nM_n$. 
We put $\chi_{l}(t)=(t/|t|)^{l}$ \ ($t\in F^\times$) 
as in \S \ref{subsec:Cn_def_ps}, and 
set $\nu''=-\sum_{i=1}^n(\nu_i+\nu_i')$. 
We define a subspace $I_{P_{n}}(l,\nu'',s)$ 
of $C^\infty (G_n)$ consisting of 
all functions $f$ such that  
\begin{align*}
&f(pg)=\chi_l(p_{n,n})
|p_{n,n}|_F^{\nu''-ns}|\det p|_F^{s}f(g)&
&(p=(p_{i,j})\in P_{n},\ g\in G_n), 
\end{align*}
on which $G_n$ acts by the right translation 
$\Pi_{P_{n},l,\nu'',s}=R$. 
The representation 
$(\Pi_{P_{n},l,\nu'',s},I_{P_{n}}(l,\nu'',s))$ is called 
a degenerate principal series representation of $G_n$.  
Similar to the proof of \cite[Proposition 7]{Grenie_001}, 
we can specify the minimal $K_n$-type of $\Pi_{P_{n},l,\nu'',s}$, 
which occur in $\Pi_{P_{n},l,\nu'',s}|_{K_n}$ 
with multiplicity $1$.  
If $l\geq 0$, we know that $\tau_{(l,\mathbf{0}_{n-1})}|_{K_n}$ is 
the minimal $K_n$-type of $\Pi_{P_{n},l,\nu'',s}$, 
and there is a $K_n$-homomorphism $\rf_{P_{n},l,\nu'',s}\colon 
V_{(l,\mathbf{0}_{n-1})}\to I_{P_{n}}(l,\nu'',s)$ 
characterized by 
\begin{align*}
&\rf_{P_{n},l,\nu'',s}(\xi_{Q(\gamma )})(g)
=\frac{|\det g|_F^{s}\prod_{i=1}^n g_{n,i}^{\gamma_i}}
{(\sum_{i=1}^n|g_{n,i}|^2)^{(ns\rc_F-\nu'' \rc_F+l)/2}} 
\end{align*}
for $g=(g_{i,j})\in G_n$ and 
$\gamma =(\gamma_1,\gamma_2,\cdots ,\gamma_n)\in \bN_0^n$ 
such that $\ell (\gamma )=l$. 
If $l\leq 0$, we know that 
$\overline{\tau_{(-l,\mathbf{0}_{n-1})}}|_{K_n}$ is 
the minimal $K_n$-type of $\Pi_{P_{n},l,\nu'',s}$, 
and there is a $K_n$-homomorphism  $\bar{\rf}_{P_{n},l,\nu'',s}\colon 
\overline{V_{(-l,\mathbf{0}_{n-1})}}\to I_{P_{n}}(l,\nu'',s)$ 
characterized by 
\begin{align*}
&\bar{\rf}_{P_{n},l,\nu'',s}(\overline{\xi_{Q(\gamma )}})(g)
=\frac{|\det g|_F^{s}\prod_{i=1}^n \overline{g_{n,i}}^{\,\gamma_i}}
{(\sum_{i=1}^n|g_{n,i}|^2)^{(ns\rc_F-\nu'' \rc_F-l)/2}} 
\end{align*}
for $g=(g_{i,j})\in G_n$ and 
$\gamma =(\gamma_1,\gamma_2,\cdots ,\gamma_n)\in \bN_0^n$ 
such that $\ell (\gamma )=-l$.

For $f\in I_{P_{n}}(l,\nu'',s)$, we define an integral 
\begin{align}
\label{eqn:def_zeta_int3}
&Z_{P_n}(W,W',f)
=\int_{Z_nN_{n}\backslash G_{n}}
W(g)W'(g)f(g)
\,dg. 
\end{align}
This integral is equivalent to (\ref{eqn:def_zeta_int2}) 
via the correspondence 
\begin{align}
\label{eqn:rel_zeta_int_23}
&Z(s,W,W',\phi)=Z_{P_n}(W,W',\rg_{P_{n},l,\nu'',s}(\phi ))
\end{align}
with $\rg_{P_{n},l,\nu'',s}(\phi )\in 
I_{P_{n}}(l,\nu'' ,s)$ defined by 
\begin{align*}
&\rg_{P_{n},l,\nu'',s}(\phi )(g)
=|\det g|_F^{s}
\int_{G_1}
\chi_{-l}(h)\phi (he_ng)|h|_F^{ns-\nu''}\,dh&
&(g\in G_n).
\end{align*}
For $g\in G_n$ and 
$\gamma =(\gamma_1,\gamma_2,\cdots ,\gamma_n)\in \bN_0^n$ 
such that $\ell (\gamma )=|l|$, we have 
\begin{align*}
&\rg_{P_{n},l,\nu,s}(\varphi_{1,n}^{(l)}(\xi_{Q(\gamma )}))(g)
=\Gamma_F(ns-\nu ;l)
\rf_{P_{n},l,\nu,s}(\xi_{Q(\gamma )})(g)&
&\text{if $l\geq 0$},\\
&\rg_{P_{n},l,\nu,s}(\overline{\varphi}_{1,n}^{(-l)}
(\overline{\xi_{Q(\gamma )}}))(g)
=\Gamma_F(ns-\nu ;-l)
\bar{\rf}_{P_{n},l,\nu,s}(\overline{\xi_{Q(\gamma )}})(g)&
&\text{if $l\leq 0$}
\end{align*}
using 
\begin{equation}
\label{eqn:Gamma_F_integral}
\begin{aligned}
\int_0^{\infty} \exp (-\pi \rc_Frt^2)t^{s\rc_F+m}
\,\frac{2\rc_F\,dt}{t}
=\frac{\Gamma_F(s;m)}{r^{(s\rc_F+m)/2}}&\\
(r\in \bR_{+}^\times ,\ m\in \bZ,\ 
\mathrm{Re}(s\rc_F+m)>0)&.
\end{aligned}
\end{equation}
Hence, Theorem \ref{thm:main2} gives the 
explicit descriptions of 
the integral (\ref{eqn:def_zeta_int3}) at the minimal 
$K_n\times K_{n}\times K_{n}$-type of 
$\Pi_{d,\nu}\boxtimes \Pi_{d',\nu'}\boxtimes \Pi_{P_{n},l,\nu'',s}$. 
We note that 
\begin{align}
\label{eqn:zeta_hom_2}
&v_1\otimes \overline{v_2}\otimes v_3
\mapsto 
Z_{P_n}(s,\rW_{\varepsilon}(\rf_{d,\nu}(v_1)),
\rW_{-\varepsilon}(\bar{\rf}_{d',\nu'}(\overline{v_2})),
\rf_{P_{n},l,\nu'',s}(v_3))
\end{align}
defines an element of 
$\Hom_{K_{n}}(V_{d}\otimes_\bC \overline{V_{-d'}}
\otimes_\bC V_{(l,\mathbf{0}_{n-1})},\bC_{\mathrm{triv}})$ 
if $l\geq 0$, and 
\begin{align}
\label{eqn:zeta_hom_3}
&v_1\otimes \overline{v_2}\otimes \overline{v_3}
\mapsto 
Z_{P_n}(s,\rW_{\varepsilon}(\rf_{d,\nu}(v_1)),
\rW_{-\varepsilon}(\bar{\rf}_{d',\nu'}(\overline{v_2})),
\bar{\rf}_{P_{n},l,\nu'',s}(\overline{v_3}))
\end{align}
defines an element of 
$\Hom_{K_{n}}(V_{d}\otimes_\bC \overline{V_{-d'}}
\otimes_\bC \overline{V_{(-l,\mathbf{0}_{n-1})}},\bC_{\mathrm{triv}})$ 
if $l\leq 0$. 
By Theorem \ref{thm:main2} with Lemma \ref{lem:schur_pairing0}, 
we obtain the following corollary.

\begin{cor}
\label{cor:main2}
Retain the notation. 
\vspace{1mm}

\noindent 
(1) Assume $-d'\in \Xi^\circ (d)$. Then 
\begin{align*}
&\sum_{M\in \rG (d)}\,\sum_{M'\in \rG (-d')}\,
\sum_{P\in \rG ((l,\mathbf{0}_{n-1}))}
\frac{\rc^{M,P}_{M'}}{\rr (M')}\,
\xi_M\otimes \overline{\xi_{M'}}\otimes \xi_{P}, 
\end{align*}
is a unique $\bQ$-rational $K_n$-invariant vector in 
$V_d\otimes_\bC \overline{V_{-d'}}
\otimes_\bC V_{(l,\mathbf{0}_{n-1})}$ up to scalar multiple, 
and its image under 
the $K_{n}$-homomorphism (\ref{eqn:zeta_hom_2}) is 
\begin{align*}
&\frac{(-\varepsilon \sI )^{\sum_{i=1}^{n-1}(n-i)(d_i+d_i')}
\rb (-d-d')\rC^\circ (-d';d)}
{\boldsymbol{\Gamma}_F(\nu;d)\boldsymbol{\Gamma}_F(\nu';d')}
\frac{L(s,\Pi_{d,\nu}\times \Pi_{d',\nu'})}{\Gamma_F(ns-\nu'';l)}.
\end{align*}
Here $\rb (-d-d')$ and $\rC^\circ (-d';d)$ are the nonzero rational 
constants, which are given by 
(\ref{eqn:def_c_gamma}) and 
(\ref{eqn:const_CG_coeff}), respectively. \vspace{1mm}

\noindent 
(2) Assume $d\in \Xi^\circ (-d')$. 
Then 
\begin{align*}
&\sum_{M\in \rG (d)}\,\sum_{M'\in \rG (-d')}\,
\sum_{P\in \rG ((-l,\mathbf{0}_{n-1}))}
\frac{\rc^{M',P}_{M}}{\rr (M)}\,
\xi_M\otimes \overline{\xi_{M'}}\otimes \overline{\xi_{P}}, 
\end{align*}
is a unique $\bQ$-rational $K_n$-invariant vector in 
$V_d\otimes_\bC \overline{V_{-d'}}
\otimes_\bC \overline{V_{(-l,\mathbf{0}_{n-1})}}$ up to scalar multiple, 
and its image under 
the $K_{n}$-homomorphism (\ref{eqn:zeta_hom_3}) is 
\begin{align*}
\frac{(-\varepsilon \sI )^{\sum_{i=1}^{n-1}(n-i)(d_i+d_i')}
\rb (d+d')\rC^\circ (d;-d')}
{\boldsymbol{\Gamma}_F(\nu;d)\boldsymbol{\Gamma}_F(\nu';d')}
\frac{L(s,\Pi_{d,\nu}\times \Pi_{d',\nu'})}{\Gamma_F(ns-\nu'';-l)}.
\end{align*}
Here $\rb (d+d')$ and $\rC^\circ (d;-d')$ are the nonzero rational 
constants, which are  given by 
(\ref{eqn:def_c_gamma}) and 
(\ref{eqn:const_CG_coeff}), respectively. 
\end{cor}

\begin{rem}
\label{rem:cohomological_2}
We set $F=\bC$. By \cite[Proposition 3.3]{Dong_Xue_001}, we note that 
the compatible pairs of cohomological representations 
of $G_n$ in Greni\'{e} \cite{Grenie_001} 
can be regarded as 
pairs of some irreducible principal series 
representations $\Pi_{d,\nu}$ and $\Pi_{d',\nu'}$ 
with $d,-d'\in \Lambda_{n,F}$ such that 
either $-d'\in \Xi^\circ (d)$ or $d\in \Xi^\circ (-d')$ holds. 
Hence, Theorem \ref{thm:main2} gives a proof of 
Greni\'{e}'s conjecture \cite[Conjecture 1]{Grenie_001} 
at all critical points 
(Dong and Xue \cite{Dong_Xue_001} proved 
this conjecture only at the central 
critical point by another method). 
Corollary \ref{cor:main2} gives 
the explicit descriptions of 
the archimedean part of Greni\'{e}'s theorem 
\cite[Theorem 2]{Grenie_001}. 
\end{rem}

\begin{rem}
Although we use the orthonormal basis $\{\zeta_M\}_{M\in \rG (\lambda )}$ 
rather than $\{\xi_M\}_{M\in \rG (\lambda )}$ in the proofs, 
we state the main theorems in terms of 
the $\bQ$-rational basis $\{\xi_M\}_{M\in \rG (\lambda )}$ 
because of the applications in the above remark. 
\end{rem}

\section{Recurrence relations}
\label{sec:relation}

\subsection{The Godement section ($G_{n-1}\to G_n$)}
\label{subsec:god_Jacquet}

Let us recall the Godement section, 
which is defined by Jacquet in \cite[\S 7.1]{Jacquet_001}. 
Assume $n>1$. 
Let $d=(d_1,d_2,\cdots ,d_n)\in \bZ^n$ and 
$\nu =(\nu_1,\nu_2,\cdots ,\nu_n)\in \bC^n$. 
We set $\widehat{d} =(d_1,d_2,\cdots ,d_{n-1})\in \bZ^{n-1}$ and 
$\widehat{\nu} =(\nu_1,\nu_2,\cdots ,\nu_{n-1})\in \bC^{n-1}$. 
Let $f\in I(\widehat{d})_{K_{n-1}}$, 
and we denote by $f_{\widehat{\nu}}$ the standard section 
corresponding to $f$. 
Let $\phi \in \cS_0 (\rM_{n-1,n}(F))$. 
When $\mathrm{Re}(\nu_n-\nu_i)>-1$ ($1\leq i\leq n-1$), 
we define the Godement section $\rg^+_{d_n,\nu_n}(f_{\widehat{\nu}},\phi )$ 
by the convergent integral 
\begin{equation*}
\begin{aligned}
&\rg^+_{d_n,\nu_n}(f_{\widehat{\nu}},\phi )(g)=
\chi_{d_n}(\det g)|\det g|_F^{\nu_n+(n-1)/2}\\
&\hspace{8mm}
\times \int_{G_{n-1}}
\phi((h,O_{n-1,1})g)
f_{\widehat{\nu}}(h^{-1})\chi_{d_n}(\det h)
|\det h|_F^{\nu_n+n/2}\,dh
\end{aligned}
\end{equation*}
for $g\in G_n$. Here we set $\chi_{l}(t)=(t/|t|)^{l}$ 
($l\in \bZ$, $t\in F^\times$) 
as in \S \ref{subsec:Cn_def_ps}. 
Jacquet shows that $\rg^+_{d_n,\nu_n}(f_{\widehat{\nu}},\phi )(g)$ 
extends to a meromorphic function of $\nu_n$ in $\bC$, 
which is a holomorphic multiple of 
\[
\prod_{1\leq i\leq n-1}\Gamma_F(\nu_n-\nu_i+1;\,|d_n-d_i|).
\]
Moreover, 
$\rg^+_{d_n,\nu_n}(f_{\widehat{\nu}},\phi )$ is an element of 
$I(d,\nu )_{K_n}$ if it is defined. 
For later use, we prepare the following lemma. 
\begin{lem}
Retain the notation. Then we have 
\begin{align}
\label{eqn:god+_Kact1}
&\Pi_{d,\nu}(k)\rg^+_{d_n,\nu_n}(f_{\widehat{\nu}},\phi )
=(\det k)^{d_n}\rg^+_{d_n,\nu_n}(f_{\widehat{\nu}},R(k)\phi )&
&(k\in K_n),\\
\label{eqn:god+_Kact2}
&(\det k')^{-d_n}
\rg^+_{d_n,\nu_n}
\bigl(\Pi_{\widehat{d},\widehat{\nu}}(k')f_{\widehat{\nu}},\,
L(k')\phi \bigr)
=\rg^+_{d_n,\nu_n}(f_{\widehat{\nu}},\phi )&
&(k'\in K_{n-1}).
\end{align}
\end{lem}
\begin{proof}
When $\mathrm{Re}(\nu_n-\nu_i)>-1$ ($1\leq i\leq n-1$), 
the equalities (\ref{eqn:god+_Kact1}) and (\ref{eqn:god+_Kact2}) 
follow immediately from the definition. 
Hence, by the uniqueness of the analytic continuations, 
we obtain the assertion. 
\end{proof}

Let $\varepsilon \in \{\pm 1\}$. 
In \cite[\S 7.2]{Jacquet_001}, 
Jacquet gives convenient integral representations of 
Whittaker functions. If $\nu$ satisfies (\ref{eqn:cdn_conv_Jac}), 
then for $g\in G_n$, we have 
\begin{equation}
\label{eqn:W_god+}
\begin{aligned}
&\rW_{\varepsilon}
\bigl(\rg^+_{d_n,\nu_n}(f_{\widehat{\nu}},\phi )\bigr)(g)=
\chi_{d_n}(\det g)|\det g|_F^{\nu_n+(n-1)/2}\\
&\hspace{10mm}\times 
\int_{G_{n-1}}\left(\int_{\rM_{n-1,1}(F)}
\phi\left(\left(h,hz\right)g\right)
\psi_{-\varepsilon }(e_{n-1}z)\,dz\right)\\
&\hspace{10mm}\times 
\rW_{\varepsilon}(f_{\widehat{\nu}})(h^{-1})
\chi_{d_n}(\det h)|\det h|_F^{\nu_n+n/2}
\,dh,
\end{aligned}
\end{equation}
where $e_{n-1}=(O_{1,n-2},1)\in \rM_{1,n-1}(F)$. 
The right hand side of (\ref{eqn:W_god+}) 
converges absolutely for all $\nu \in \bC^n$, 
and defines an entire function of $\nu$. 
Thus the equality holds for all $\nu $. 
In the appendix of this paper, 
we show that the integral representation (\ref{eqn:W_god+}) 
can be regarded as a generalization of 
the recursive formula \cite[Theorem 14]{Ishii_Stade_001} 
of spherical Whittaker functions.

\subsection{The Godement section ($G_{n}\to G_n$)}

In this subsection, we define a new kind of the Godement section. 
Let $d\in \bZ^n$ and $\nu \in \bC^n$. 
Let $f\in I(d)_{K_n}$ and $\phi \in \cS_0 (\rM_{n}(F))$. 
We denote by $f_{\nu}$ the standard section corresponding to $f$. 
For $s\in \bC$, $l\in \bZ$ and $g\in G_n$, we set 
\begin{align}
\label{eqn:def_god}
&\rg^\circ_{l,s}(f_\nu ,\phi )(g)=
\int_{G_{n}}
f_\nu (gh)\phi (h)
\chi_{l}(\det h)|\det h|_F^{s+(n-1)/2}\,dh.
\end{align}

\begin{prop}
\label{prop:god0_Wgod0}
Let $d\in \bZ^n$, $l\in \bZ$ and $\varepsilon \in \{\pm 1\}$. 
Let $\Omega$ be an open, relatively compact subset of $\bC^n$. 
Then there is a constant $c_0$ such that, for any  
$f\in I(d)_{K_n}$ and $\phi \in \cS_0(\rM_{n}(F))$, 
the following assertions (i) and (ii) hold: \vspace{1.5mm}

\noindent (i) The integral (\ref{eqn:def_god}) converges absolutely and 
uniformly on any compact subset of 
$\{(s,\nu,g)\in \bC \times \Omega \times G_n \mid 
\mathrm{Re}(s)>c_0\}$. \vspace{1.5mm}

\noindent (ii) 
Let $\nu \in \Omega$ and $s\in \bC$ such that $\mathrm{Re}(s)>c_0$. 
Then $\rg^\circ_{l,s}(f_\nu ,\phi )$ is an element 
of $I(d,\nu )_{K_n}$ satisfying 
\begin{align}
\label{eqn:god0_Kact1}
&\Pi_{d,\nu}(k)\rg^\circ_{l,s}(f_\nu ,\phi )
=(\det k)^{-l}\rg^\circ_{l,s}(f_\nu ,L(k)\phi )&
&(k\in K_n),\\
\label{eqn:god0_Kact2}
&(\det k')^{l}\,
\rg^\circ_{l,s}\bigl(\Pi_{d,\nu}(k')f_\nu ,\,
R(k')\phi \bigr)
=\rg^\circ_{l,s}(f_\nu ,\phi )&
&(k'\in K_{n}).
\end{align}
Moreover, for $g\in G_n$, we have 
\begin{equation}
\label{eqn:W_god0}
\begin{aligned}
&\rW_{\varepsilon}(\rg^\circ_{l,s}(f_\nu ,\phi ))(g)
=\int_{G_{n}}\rW_{\varepsilon}(f_\nu)(gh)\phi (h)
\chi_{l}(\det h)|\det h|_F^{s+(n-1)/2}\,dh. 
\end{aligned}
\end{equation}
Here $f_{\nu}$ is the standard section corresponding to $f$. 
\end{prop}
\begin{proof}
For $g\in G_n$, we set 
$\| g\|=\mathrm{Tr}(g{}^t\overline{g})
+\mathrm{Tr}((g^{-1}){}^{t}\!(\overline{g^{-1}}))$ 
and denote by 
\begin{align*}
&g=\ru (g)\ra (g)\rk (g)&
&(\ru (g)\in U_n,\ \ra (g)\in A_n,\ \rk (g)\in K_n)
\end{align*}
the decomposition of $g$ according to $G_n=U_nA_nK_n$. 
It is easy to see that  
\begin{align}
\label{eqn:pf_god_001}
&\|\ra (g)\| \leq \|g\|=\|kgk'\|,&
&\|gh\|\leq \|g\|\,\|h\|,&
&\ra (gh)=\ra (g)\ra (\rk (g)h)&
\end{align}
for $g,h\in G_n$ and $k,k'\in K_n$. 
Since $G_n\ni g\mapsto \eta_{\nu -\rho_n}(\ra (g))\in \bC$ 
is an element of $I(\mathbf{0}_n,\nu )$, we have 
\begin{align}
\label{eqn:pf_god_004}
&\int_{N_n}|\eta_{\nu -\rho_n}(\ra (x))|\,dx<\infty &
&\text{($\nu \in \bC^n$ satisfying (\ref{eqn:cdn_conv_Jac}))}.
\end{align}
by the absolute convergence of the Jacquet integral 
\cite[Theorem 15.4.1]{Wallach_003}.

We take $d$, $l$, $\varepsilon $ and $\Omega$ 
as in the statement. Replacing $\Omega$ with its superset if necessary, 
we may assume that $\Omega$ contains an element $\nu$ satisfying 
(\ref{eqn:cdn_conv_Jac}). 
By (\ref{eqn:pf_god_001}) and \cite[Proposition 3.2]{Jacquet_001}, 
there are a constant $c_1$ and a continuous semi-norm $\cQ$ on 
$I(d)$ such that, 
for any  $\nu \in \Omega$, $g\in G_n$ and $f\in I(d)$, 
the following inequalities hold:
\begin{align}
\label{eqn:pf_god_002}
&|\eta_{\nu -\rho_n}(\ra (g))|\leq \|g\|^{c_1},&
&|\rW_{\varepsilon}(f_\nu)(g)|\leq \|g\|^{c_1}\cQ (f). 
\end{align}
By \cite[Lemma 3.3 (ii)]{Jacquet_001}, 
there is a positive constant $c_0$ such that, 
for any $t>c_0$ and $\phi \in \cS (\rM_{n}(F))$, the integral 
\begin{align}
\label{eqn:pf_god_005}
&\int_{G_n}\|h\|^{c_1}\phi (h)|\det h|_F^{t+(n-1)/2}\,dh
\end{align}
converges absolutely. 

Let $f\in I(d)_{K_n}$ and $\phi \in \cS_0(\rM_{n}(F))$. 
By (\ref{eqn:pf_god_001}), (\ref{eqn:pf_god_002}) 
and the definition of $I(d,\nu )$, 
for $\nu \in \Omega$, $x\in N_n$ and $g,h\in G_n$, 
we have an estimate 
\begin{align}
\label{eqn:pf_god_003}
&|f_\nu (xgh)|\leq |\eta_{\nu -\rho_n}(\ra (x))|\,
\|g\|^{c_1}\|h\|^{c_1}
\sup_{k\in K_n}|f(k)|.
\end{align}
By the absolute convergence of (\ref{eqn:pf_god_005}) 
and the estimate (\ref{eqn:pf_god_003}) 
with $x=1_n$, we obtain the assertion (i). 

Let $\nu \in \Omega$ and $s\in \bC$ such that $\mathrm{Re}(s)>c_0$. 
By definition, we have the equalities 
(\ref{eqn:god0_Kact1}), (\ref{eqn:god0_Kact2}) and 
\begin{align}
\label{eqn:pf_god_006}
&\rg^\circ_{l,s}(f_\nu ,\phi )(umg)
=\chi_{d}(m)\eta_{\nu -\rho_n}(m)\rg^\circ_{l,s}(f_\nu ,\phi )(g)
\end{align}
for $u\in U_n$, $m\in M_n$ and $g\in G_n$. 
Since $\Pi_{d,\nu}$ is admissible and 
$\rg^\circ_{l,s}(f_\nu ,\phi )$ is a continuous $K_n$-finite function 
on $G_n$ satisfying (\ref{eqn:pf_god_006}), 
we know that $\rg^\circ_{l,s}(f_\nu ,\phi )$ is smooth and 
an element of $I(d,\nu )_{K_n}$ by \cite[Propositions 8.4 and 8.5]{Knapp_002}. 

Let $g\in G_n$. If $\nu \in \Omega$ satisfies 
(\ref{eqn:cdn_conv_Jac}), we obtain the equality 
(\ref{eqn:W_god0}) as follows:  
\begin{align*}
&\rW_{\varepsilon}(\rg^\circ_{l,s}(f_\nu ,\phi ))(g)=
\int_{N_n}\rg^\circ_{l,s}(f_\nu ,\phi )(xg)
\psi_{-\varepsilon,n} (x)\,dx\\
&=\int_{N_n}\left(\int_{G_{n}}
f(xgh)\phi (h)\chi_{l}(\det h)|\det h|_F^{s+(n-1)/2}\,dh
\right)\psi_{-\varepsilon,n}(x)\,dx\\
&=\int_{G_{n}}\left(\int_{N_n}
f(xgh)\psi_{-\varepsilon,n}(x)\,dx
\right)\phi (h)
\chi_{l}(\det h)|\det h|_F^{s+(n-1)/2}\,dh\\
&=\int_{G_{n}}\rW_{\varepsilon}(f_\nu)(gh)\phi (h)
\chi_{l}(\det h)|\det h|_F^{s+(n-1)/2}\,dh.
\end{align*}
Here the third equality is justified by 
Fubini's theorem, since the double integral 
converges absolutely by (\ref{eqn:pf_god_004}), 
(\ref{eqn:pf_god_003}) and 
the absolute convergence of (\ref{eqn:pf_god_005}). 

In order to complete the proof, it suffices to show that 
the both sides of (\ref{eqn:W_god0}) are holomorphic 
functions of $(s,\nu )$ on a domain 
\begin{equation}
\label{eqn:pf_god_007}
\{(s,\nu )\in \bC \times \Omega \mid \mathrm{Re}(s)>c_0\}. 
\end{equation}
By (\ref{eqn:pf_god_001}), (\ref{eqn:pf_god_002}) 
and the absolute convergence of (\ref{eqn:pf_god_005}), 
the integral in the right hand side of (\ref{eqn:W_god0}) 
converges absolutely and uniformly on any compact subset of 
the domain (\ref{eqn:pf_god_007}), 
and defines a holomorphic function on the domain 
(\ref{eqn:pf_god_007}). 

Let $\cS_{\phi,l}$ be a subspace of $\cS_0(\rM_n(F))$ 
spanned by $L(k)\phi $ \ ($k\in K_n$), 
and we regard $\cS_{\phi,l}$ as a $K_n$-module via the action 
$\det^{-l}\otimes L$. 
Let $I_{\phi,l}$ be a subspace of $I(d)_{K_n}$ spanned by 
$\{T(\phi')\mid \phi'\in \cS_{\phi,l}$, 
$T\in \Hom_{K_n}(\cS_{\phi,l},\, I(d)_{K_n})\}$. 
Then we have $\rg^\circ_{l,s}(f_\nu ,\phi )|_{K_n}\in I_{\phi,l}$ 
by (\ref{eqn:god0_Kact1}). 
Since $\phi$ is $K_n$-finite and $\Pi_{d,\nu}$ is admissible, 
the space $I_{\phi,l}$ is finite dimensional. 
Let $\{f_{\phi,i}\}_{i=1}^m$ be an orthonormal basis of 
$I_{\phi,l}$ with respect to the $L^2$-inner product 
\begin{align*}
&\langle f_1,f_2\rangle_{L^2}=\int_{K_n}f_1(k)\overline{f_2(k)}\,dk&
&(f_1,f_2\in I(d)). 
\end{align*}
Since 
$\rg^\circ_{l,s}(f_\nu ,\phi )|_{K_n}=\sum_{i=1}^m
\langle \rg^\circ_{l,s}(f_\nu ,\phi )|_{K_n},
f_{\phi,i}\rangle_{L^2}\,f_{\phi,i}$, 
we have 
\begin{align*}
&\rW_{\varepsilon}(\rg^\circ_{l,s}(f_\nu ,\phi ))(g)
=\sum_{i=1}^m\langle \rg^\circ_{l,s}(f_\nu ,\phi )|_{K_n},\,
f_{\phi,i}\rangle_{L^2}\,
\rW_{\varepsilon}(f_{\phi,i,\nu})(g), 
\end{align*}
where $f_{\phi,i,\nu}$ is the standard section 
corresponding to $f_{\phi,i}$. 
By this expression and the statement (i), we know that 
the right hand side of (\ref{eqn:W_god0}) 
is holomorphic on the domain (\ref{eqn:pf_god_007}). 
\end{proof}

\begin{rem}
The equality (\ref{eqn:W_god0}) with $l=0$ can be regarded 
as the local theta correspondence for a principal series 
representation $\Pi_{d,\nu}$ in \cite[\S 2]{Watanabe_001}. 
\end{rem}

\subsection{Recurrence relations with the Godement sections}
\label{subsec:rec_rel_godement}

Let $\varepsilon \in \{\pm 1\}$. For $\phi \in \cS (\rM_{n,1}(F))$, 
we define $\cF_\varepsilon (\phi )\in \cS (\rM_{1,n}(F))$ by 
\begin{align}
\label{eqn:def_fourier}
&\cF_\varepsilon (\phi )(t)=\int_{\rM_{n,1}(F)}\phi (z)
\psi_{-\varepsilon }(tz)\,d_Fz&
&(t\in \rM_{1,n}(F)). 
\end{align}
Let 
\begin{equation*}
\begin{aligned}
&d=(d_1,d_2,\cdots ,d_n)\in \bZ^n, \qquad &
&\nu =(\nu_1,\nu_2,\cdots ,\nu_n)\in \bC^n, \\
&d'=(d_1',d_2',\cdots ,d_{n'}')\in \bZ^{n'},& 
&\nu'=(\nu_1',\nu_2',\cdots ,\nu_{n'}')\in \bC^{n'}.
\end{aligned}
\end{equation*}
If $n>1$, we set $\widehat{d}=(d_1,d_2,\cdots ,d_{n-1})$ and 
$\widehat{\nu}=(\nu_1,\nu_2,\cdots ,\nu_{n-1})$. 
If $n'>1$, we set $\widehat{d'}=(d_1',d_2',\cdots ,d_{n'-1}')$ and 
$\widehat{\nu'}=(\nu_1',\nu_2',\cdots ,\nu_{n'-1}')$.

\begin{prop}[{$G_n\times G_{n}\to G_{n}\times G_{n-1}$}]
\label{prop:zeta_god_recurrence1}
Retain the notation, and assume $n'=n>1$. 
Let $f\in I(d)_{K_n}$ and $f'\in I(\widehat{d'})_{K_{n-1}}$. 
We denote by $f_\nu$ and $f'_{\widehat{\nu'}}$ 
the standard sections corresponding to $f$ and $f'$, respectively. 
Let $\phi_1\in \cS_0 (\rM_{n-1,n}(F))$ and 
$\phi_2\in \cS_0 (\rM_{1,n}(F))$. 
For $s\in \bC$ such that $\mathrm{Re}(s)$ is sufficiently large, we have 
\begin{align*}
&Z\bigl(s,\,\rW_{\varepsilon} (f_\nu ),\,
\rW_{-\varepsilon}\bigl(\rg^+_{d_n',\nu_n'}(f_{\widehat{\nu'}}',
\phi_1)\bigr),\,
\phi_2\bigr)\\
&=Z\bigl(s,\,\rW_{\varepsilon}
\bigl(\rg^\circ_{d_n',s+\nu_n'}(f_\nu ,\phi_0 )\bigr),\,
\rW_{-\varepsilon}(f_{\widehat{\nu'}}')\bigr),
\end{align*}
where $\phi_0\in \cS_0 (\rM_{n}(F))$ is defined by 
\begin{align*}
&\phi_0(z)=\phi_1((1_{n-1},O_{n-1,1})z)
\phi_2(e_nz)&&(z\in \rM_{n}(F)).
\end{align*}
\end{prop}
\begin{proof}
Using (\ref{eqn:W_god+}), Jacquet shows the following equality 
\cite[(8.1)]{Jacquet_001}: 
\begin{align*}
&Z\bigl(s,\,\rW_{\varepsilon} (f_\nu ),\,
\rW_{-\varepsilon}\bigl(\rg^+_{d_n',\nu_n'}(f_{\widehat{\nu'}}',
\phi_1)\bigr),\,
\phi_2\bigr)\\
&=\int_{N_{n-1}\backslash G_{n-1}}\biggl(\int_{G_n}
\rW_{\varepsilon}(f_\nu )(\iota_n(h)g)\phi_0(g)
\chi_{d_n'}(\det g)|\det g|_F^{s+\nu_n'+(n-1)/2}\,dg\biggr)\\
&\quad \times \rW_{-\varepsilon}(f_{\widehat{\nu'}}')(h)|\det h|_F^{s-1/2}dh.
\end{align*}
Hence, we obtain the assertion by 
Proposition \ref{prop:god0_Wgod0}. 
\end{proof}

\begin{prop}[{$G_n\times G_{n-1}\to G_{n-1}\times G_{n-1}$}]
\label{prop:zeta_god_recurrence2}
Retain the notation, and assume $n'=n-1$. 
Let $f\in I(\widehat{d})_{K_{n-1}}$ and $f'\in I(d')_{K_{n-1}}$. 
We denote by $f_{\widehat{\nu}}$ and $f_{\nu'}'$ 
the standard sections corresponding to $f$ and $f'$, respectively. 
Let $\phi_1\in \cS_0 (\rM_{n-1}(F))$ and 
$\phi_2\in \cS_0 (\rM_{n-1,1}(F))$. 
For $s\in \bC$ such that $\mathrm{Re}(s)$ is sufficiently large, we have 
\begin{align*}
&Z\bigl(s,\,
\rW_{\varepsilon}
\bigl(\rg^+_{d_n,\nu_n}(f_{\widehat{\nu}},\phi_0)\bigr),\,
\rW_{-\varepsilon}(f_{\nu'}')\bigr)\\
&=Z\bigl(s,\,
\rW_{\varepsilon}(f_{\widehat{\nu}}),\,
\rW_{-\varepsilon}(\rg^\circ_{d_n,s+\nu_n}(f_{\nu'}',\phi_1)),\,
\cF_\varepsilon (\phi_2)\bigr),
\end{align*}
where $\phi_0\in \cS_0 (\rM_{n-1,n}(F))$ is defined by 
\begin{align*}
&\phi_0(z)=\phi_1(z\,{}^t\!(1_{n-1},O_{n-1,1}))\phi_2(z\,{}^t\!e_{n})&
&(z\in \rM_{n-1,n}(F)). 
\end{align*}
\end{prop}
\begin{proof}
Using (\ref{eqn:W_god+}), Jacquet shows the 
following equality \cite[(8.3)]{Jacquet_001}:
\begin{align*}
&Z\bigl(s,\,
\rW_{\varepsilon}
\bigl(\rg^+_{d_n,\nu_n}(f_{\widehat{\nu}},\phi_0)\bigr),\,
\rW_{-\varepsilon}(f_{\nu'}')\bigr)\\
&=\int_{N_{n-1}\backslash G_{n-1}}
\biggl(\int_{G_{n-1}}
\rW_{-\varepsilon}(f_{\nu'}')(gh)\phi_1(h)
\chi_{d_n}(\det h)|\det h|_F^{s+\nu_n+(n-2)/2}dh
\biggr)\\
&\quad \times \rW_{\varepsilon}(f_{\widehat{\nu}})(g)
\cF_\varepsilon (\phi_2)(e_{n-1}g)|\det g|_F^{s}\,dg.
\end{align*}
Hence, we obtain the assertion by 
Proposition \ref{prop:god0_Wgod0}. 
\end{proof}

\section{Finite dimensional representations}
\label{sec:fin_dim}

\subsection{The Clebsch--Gordan coefficients}
\label{subsec:CG_coeff}

Let $\lambda =(\lambda_1,\lambda_2,\cdots ,\lambda_n)\in \Lambda_n$, 
$l\in \bN_0$ and 
$\lambda'=(\lambda_1',\lambda_2',\cdots ,\lambda_n')
\in \Xi^\circ (\lambda ;l)$. 
By Pieri's rule (\ref{eqn:pieri}), we can take 
a $\mathrm{GL}(n,\bC )$-homomorphism 
$\tilde{\rI}^{\lambda,l}_{\lambda'}\colon 
V_{\lambda'}\to V_{\lambda }\otimes_{\bC}
V_{(l,\mathbf{0}_{n-1})}$ 
satisfying 
\begin{align}
\label{eqn:Inj_normalize}
&\bigl\langle\, \tilde{\rI}^{\lambda,l}_{\lambda'}(v)
,\ \tilde{\rI}^{\lambda,l}_{\lambda'}(v')
\, \bigr\rangle 
=\langle v,v'\rangle \quad (v,v'\in V_{\lambda'}). 
\end{align} 
Such $\tilde{\rI}^{\lambda,l}_{\lambda'}$ is unique 
up to multiplication by scalars in $\rU (1)$. 
We set 
\begin{align*}
&\tilde{\rI}^{\lambda,l}_{\lambda'} (\zeta_{M'})
=\sum_{M\in \rG (\lambda)}\,\sum_{P\in \rG ((l,\mathbf{0}_{n-1}))}
\rC^{M,P}_{M'}\zeta_{M}\otimes \zeta_{P}&
&(M'\in \rG (\lambda')). 
\end{align*}
Then we call $\rC^{M,P}_{M'}$ 
($M\in \rG (\lambda)$, $P\in \rG ((l,\mathbf{0}_{n-1}))$, 
$M'\in \rG (\lambda')$) the Clebsch--Gordan coefficients. 
When $n=1$, we may normalize $\rC^{\lambda_1,l}_{\lambda_1+l}=1$, since 
\begin{align*}
&\Xi^\circ (\lambda_1;l)=\{\lambda_1+l\},&
&\rG (\lambda_1)=\{\lambda_1\},&
&\rG (l)=\{l\},&
&\rG (\lambda_1+l)=\{\lambda_1+l\}. 
\end{align*}
We consider the case $n>1$. 
Let $\mu \in \Xi^+ (\lambda )$ and $0\leq q\leq l$. 
By Lemma \ref{lem:unitary_isom_lambda_mu}, 
there are some constants 
\begin{align}
\label{eqn:isoscalar_factor}
&\left(\begin{array}{cc|c}
\lambda ,&l&\lambda'\\
\mu ,&q&\mu'\\
\end{array}\right)&
&(\mu'\in \Xi^+ (\lambda ')\cap \Xi^\circ (\mu ;q))
\end{align}
such that, for any $M'\in \rG (\lambda';\mu')$, 
$M\in \rG (\lambda ;\mu)$, 
$P\in \rG ((l,\mathbf{0}_{n-1});(q,\mathbf{0}_{n-2}))$ and 
$\mu'\in \Xi^+ (\lambda ')$, the following equality hold:
\begin{align}
\label{eqn:coef_decomp}
\rC^{M,P}_{M'}
=\left\{
\begin{array}{ll}
\displaystyle 
\left(\begin{array}{cc|c}
\lambda ,&l&\lambda'\\
\mu ,&q&\mu'\\
\end{array}\right)
\rC^{\widehat{M},\widehat{P}}_{\widehat{M'}}&
\text{if}\ \mu'\in \Xi^\circ (\mu ;q),\\[4mm]
0&\text{otherwise}.
\end{array}\right.
\end{align}
The constants (\ref{eqn:isoscalar_factor}) 
are called the isoscalar factors. 
In \cite{Jucis_001} (see also \cite{Alisauskas_Jucys_001} 
and \cite[Chapter 18]{Vilenkin_Klimyk_001}), 
Jucys gives the following expressions of them 
under some normalization of $\tilde{\rI}^{\lambda,l}_{\lambda'}$: 
\begin{equation}
\label{eqn:isoscalar_explicit}
\begin{aligned}
\left(\begin{array}{cc|c}
\lambda ,&l&\lambda'\\
\mu ,&q&\mu'\\
\end{array}\right)=\,&\sqrt{
\frac{(l-q)!\,\rS^\circ (\lambda',\lambda')\rS^+(\lambda ,\mu )
\rS^\circ (\mu',\mu )\rS^\circ (\mu ,\mu )}
{\rS^\circ (\lambda',\lambda )\rS^+(\lambda',\mu')}}\\
&\times 
\underset{\mu'\in \Xi^\circ (\alpha ),\  \alpha \in \Xi^\circ (\mu )}
{\sum_{\alpha \in \Xi^+ (\lambda )\cap \Xi^+ (\lambda')}}
\frac{(-1)^{\ell (\alpha -\mu )}\rS^\circ (\alpha ,\alpha )}
{\rS^\circ (\mu',\alpha )\rS^\circ (\alpha ,\mu )}
\frac{\rS^+(\lambda',\alpha )}{\rS^+(\lambda,\alpha )}
\end{aligned}
\end{equation}
for $\mu'\in \Xi^+ (\lambda ')\cap \Xi^\circ (\mu ;q)$, 
where the symbols $\rS^\circ (\lambda' ,\lambda )$ and $\rS^+(\lambda ,\mu )$ 
are defined by (\ref{eqn:def_S0}) and (\ref{eqn:def_S+}), respectively. 
Hereafter, 
we assume that $\tilde{\rI}^{\lambda,l}_{\lambda'}$ is normalized 
so that (\ref{eqn:isoscalar_explicit}) holds. 
Then we have 
\begin{align}
\label{eqn:CG_explicit}
&\rC^{M,Q((\mathbf{0}_{n-1},l))}_{M'}
=\sqrt{\frac{l!\,\rS^\circ (\lambda',\lambda')\rS^+(\lambda',\mu )}
{\rS^\circ (\lambda',\lambda )\rS^+(\lambda,\mu )}},&
&\rC^{H(\lambda ),Q(\lambda'-\lambda )}_{H(\lambda')}
=\sqrt{\rC^\circ (\lambda';\lambda )}
\end{align}
for $M\in \rG (\lambda )$ and $M'\in \rG (\lambda' )$ such that 
$\widehat{M}=\widehat{M'}\in \rG (\mu )$. 
Here $H(\lambda )$, 
$Q(\gamma )$ ($\gamma \in \bN_0^n$) and $\rC^\circ (\lambda';\lambda )$ 
are defined by (\ref{eqn:def_Hlambda}), (\ref{eqn:def_Qgamma}) and 
(\ref{eqn:const_CG_coeff}), respectively.  
All the Clebsch--Gordan coefficients $\rC^{M,P}_{M'}$ 
are real numbers, and we have 
\begin{align}
\label{eqn:decomp_pure_tensor}
&\begin{aligned}
\zeta_{M}\otimes \zeta_{P}
&=\sum_{\lambda' \in \Xi^\circ (\lambda ;l)}
\,\sum_{M'\in \rG (\lambda')}
\bigl\langle\, \zeta_M\otimes \zeta_P,\ 
\tilde{\rI}^{\lambda,l}_{\lambda'}(\zeta_{M'})\, \bigr\rangle\, 
\tilde{\rI}^{\lambda,l}_{\lambda'}(\zeta_{M'})\\
&=\sum_{\lambda' \in \Xi^\circ (\lambda ;l)}
\,\sum_{M'\in \rG (\lambda')}{\rC^{M,P}_{M'}}\, 
\tilde{\rI}^{\lambda,l}_{\lambda'}(\zeta_{M'})
\end{aligned}
\end{align}
for $M\in \rG (\lambda )$ and $P\in \rG ((l,\mathbf{0}_{n-1}))$. 
We set 
\begin{equation*}
\rI^{\lambda,l}_{\lambda'}=
\sqrt{\rb (\lambda'-\lambda )\rC^\circ (\lambda';\lambda )}
\,\tilde{\rI}^{\lambda,l}_{\lambda'}, 
\end{equation*}
where $\rb (\lambda'-\lambda )$ is defined by (\ref{eqn:def_c_gamma}). 
Then the explicit expression (\ref{eqn:injector_explicit}) 
of $\rI^{\lambda,l}_{\lambda'}$ 
follows from (\ref{eqn:def_rM}) and (\ref{eqn:isoscalar_explicit}), 
since 
\[
\rc^{M,P}_{M'}
=\sqrt{\frac{\rb (\lambda'-\lambda )\rC^\circ (\lambda';\lambda )\rr (M')}
{\rr (M) \rr (P)}}\rC^{M,P}_{M'}
\]
for $M\in \rG (\lambda)$, $P\in \rG ((l,\mathbf{0}_{n-1}))$ and 
$M'\in \rG (\lambda')$. 
Here $\rr (M)$ is defined by (\ref{eqn:def_rM}). 
The equalities 
(\ref{eqn:cdn_inj_2}) and (\ref{eqn:cdn_inj_3}) 
follow from (\ref{eqn:Inj_normalize}) and 
(\ref{eqn:CG_explicit}), respectively.

\subsection{Some lemmas for representations of $K_n$}
\label{subsec:lem_Krep}

Let $\bC_{\mathrm{triv}}=\bC$ be the trivial $\mathrm{GL}(n,\bC)$-module. 
The purpose of this subsection is to give proofs of 
Lemma \ref{lem:rep_Kn_irred} and the following three lemmas.

\begin{lem}
\label{lem:schur_pairing}
Let $\lambda \in \Lambda_{n,F}$. 
\vspace{1mm}

\noindent 
(1) The space 
$\Hom_{K_n}(V_\lambda \otimes_\bC \overline{V_{\lambda}},
\bC_{\mathrm{triv}})$ is a $1$ dimensional space spanned by 
the $\bC$-linear map 
\[
V_\lambda \otimes_\bC \overline{V_\lambda}
\ni v_1\otimes \overline{v_2}\mapsto 
\langle v_1,v_2\rangle \in \bC.
\]

\noindent 
(2) Let $\lambda'\in \Lambda_{n,F}\cap \Xi^\circ (\lambda )$, and 
set $l=\ell (\lambda' -\lambda )$. For $\lambda''\in \Xi^\circ (\lambda ;l)$ 
such that $\lambda''\neq \lambda'$, we have 
$\Hom_{K_n}(V_{\lambda'} \otimes_\bC \overline{V_{\lambda''}},
\bC_{\mathrm{triv}})=\{0\}$. 
\end{lem}

\begin{lem}
\label{lem:schur_pairing+}
Assume $n>1$, and we regard $K_{n-1}$ as a subgroup of $K_n$ via 
(\ref{eqn:embed_Gn-1_Gn}). 
Let $\lambda \in \Lambda_{n,F}$ and $\mu \in \Lambda_{n-1,F}$.  
Then $\Hom_{K_{n-1}}(V_\lambda \otimes_\bC \overline{V_{\mu}},
\bC_{\mathrm{triv}})$ 
is a $1$ dimensional space spanned by the $\bC$-linear map 
\[
V_\lambda \otimes_\bC \overline{V_{\mu}}\ni v_1\otimes \overline{v_2}
\mapsto \langle \rR_{\mu}^{\lambda}(v_1),\,v_2\rangle 
\in \bC_{\mathrm{triv}}
\]
if $\mu \in \Xi^+(\lambda )$, 
and is equal to $\{0\}$ otherwise. 
Here $\rR^{\lambda}_\mu$ is defined by (\ref{eqn:def_projR}). 
Moreover, if $\mu \in \Xi^+(\lambda )$, then 
\begin{equation*}
\sum_{M\in \rG (\mu )}\rr (M)^{-1}\,
\xi_{M[\lambda ]}\otimes \overline{\xi_M}
\end{equation*}
is a unique $\bQ$-rational $K_{n-1}$-invariant vector 
in $V_\lambda \otimes_\bC \overline{V_{\mu}}$ 
up to scalar multiple. 
\end{lem}

\begin{lem}
\label{lem:schur_pairing0}
Let $\lambda ,\lambda'\in \Lambda_{n,F}$ such that 
$\ell (\lambda'-\lambda )\geq 0$. Let $l\in \bN_0$. 
If $F=\bR$, we assume $l\in \{0,1\}$. 
Then $\Hom_{K_{n}}(V_{\lambda'} \otimes_\bC \overline{V_{\lambda}}
\otimes_\bC \overline{V_{(l,\mathbf{0}_{n-1})}},
\bC_{\mathrm{triv}})$ is a $1$ dimensional space spanned by 
the $\bC$-linear map 
\[
V_{\lambda'} \otimes_\bC \overline{V_{\lambda}}
\otimes_\bC \overline{V_{(l,\mathbf{0}_{n-1})}}
\ni v_1\otimes \overline{v_2}\otimes \overline{v_3}
\mapsto \langle \mathrm{I}^{\lambda ,l}_{\lambda'}(v_1),\,
v_2\otimes v_3\rangle 
\in \bC_{\mathrm{triv}} 
\]
if $\lambda'\in \Xi^\circ (\lambda ;l)$, 
and is equal to $\{0\}$ otherwise. 
Moreover, if $\lambda'\in \Xi^\circ (\lambda ;l)$, then 
\begin{align}
\label{eqn:lem_schur_pairing0}
&\sum_{M'\in \rG (\lambda')}\,\sum_{M\in \rG (\lambda)}\,
\sum_{P\in \rG ((l,\mathbf{0}_{n-1}))}
\frac{\rc^{M,P}_{M'}}{\rr (M')}\,
\xi_{M'}\otimes \overline{\xi_{M}}\otimes \overline{\xi_{P}}, 
\end{align}
is a unique $\bQ$-rational $K_n$-invariant vector in 
$V_{\lambda'} \otimes_\bC \overline{V_{\lambda}}
\otimes_\bC \overline{V_{(l,\mathbf{0}_{n-1})}}$ up to scalar multiple. 
\end{lem}

Since proofs of Lemmas \ref{lem:rep_Kn_irred}, 
\ref{lem:schur_pairing}, \ref{lem:schur_pairing+} 
and \ref{lem:schur_pairing0} 
are easy for $F=\bC$, 
the main concern is the case of $F=\bR$. 
We have $\Lambda_{n,\bR}=\{(\mathbf{1}_j,\mathbf{0}_{n-j})\mid 
0\leq j\leq n\}$ with 
$\mathbf{1}_{j}=(1,1,\cdots ,1)\in \bZ^{j}$ 
and $\mathbf{0}_{n-j}=(0,0,\cdots ,0)\in \bZ^{n-j}$. 
Here we erase the symbol $\mathbf{1}_j$ if $j=0$, 
and erase the symbol $\mathbf{0}_{n-j}$ if $j=n$. 
Let $0\leq l\leq n$, and 
we regard the $l$-th exterior power
${\textstyle \bigwedge^l} (\rM_{n,1}(\bC))$ of $\rM_{n,1}(\bC)$ 
as a $\mathrm{GL}(n,\bC )$-module via the action derived from 
the matrix multiplication. Then we have 
$V_{(\mathbf{1}_l,\mathbf{0}_{n-l})}
\simeq {\textstyle \bigwedge^l} (\rM_{n,1}(\bC))$ 
as $\mathrm{GL}(n,\bC)$-modules via the correspondence 
\begin{align*}
&\zeta_M\leftrightarrow \gge_{i_1}\wedge \gge_{i_2}\wedge \cdots 
\wedge \gge_{i_l}&
&(M\in \rG ((\mathbf{1}_l,\mathbf{0}_{n-l})))
\end{align*}
with $1\leq i_1<i_2<\cdots <i_l\leq n$ such that 
$\gamma_{i}^M=1$ ($i\in \{i_1,i_2,\cdots ,i_l\}$). 
Here 
$\gge_{j}$ is the matrix unit 
in $\rM_{n,1}(\bC )$ with $1$ at the $(j,1)$-th entry $(1\leq j\leq n)$ 
and $0$ at other entries. 
We identify $V_{(\mathbf{1}_l,\mathbf{0}_{n-l})}$ with 
${\textstyle \bigwedge^l} (\rM_{n,1}(\bC))$ via this isomorphism.

We have $\rO (n)=\mathrm{SO}(n)\sqcup \mathrm{SO}(n)k_0$ with 
$k_0=\diag (1,1,\cdots ,1,-1)\in \rO (n)$ and 
$\mathrm{SO}(n)=\{k\in \rO (n)\mid \det k=1\}$. 
The complexification $\gs \go (n)_\bC$ of 
the associated Lie algebra $\gs \go (n)$ of $\mathrm{SO}(n)$ is given by 
$\gs \go (n)_{\bC}=\bigoplus_{1\leq i<j\leq n}\bC E_{i,j}^{\gs \go (n)}$ 
with $E_{i,j}^{\gs \go (n)}=E_{i,j}-E_{j,i}$. 
Here we understand $k_0=-1$ and $\gs \go (1)_{\bC}=\{0\}$ if $n=1$. 
Let us recall some facts 
in the highest weight theory \cite[Theorem 4.28]{Knapp_002} 
for $\mathrm{SO}(n)$. 
Let $m$ be the largest integer such that $2m\leq n$. 
When $n\geq 2$, for an irreducible representation $(\tau ,V_\tau )$ 
of $\mathrm{SO}(n)$, 
there is a nonzero vector $v_0$ in $V_\tau$ such that, 
for $1\leq i \leq m$ and $2i+1\leq j\leq n$,  
\begin{align*}
&\tau \bigl(E_{2i-1,2i}^{\gs \go (n)}\bigr)v_0
=\sI \lambda_iv_0,&
&\tau \bigl(E_{2i-1,j}^{\gs \go (n)}+\sI E_{2i,j}^{\gs \go (n)}\bigr)v_0=0
\end{align*} 
with some $\lambda =(\lambda_1,\lambda_2,\cdots ,\lambda_m)\in \bZ^m$.  
Such vector $v_0$ is unique up to nonzero scalar multiple, 
and we call $v_0$ an $\mathrm{SO}(n)$-highest weight vector 
of weight $\lambda$.  
The weight $\lambda_\tau =\lambda $ is 
called the highest weight of $\tau$, 
and $\tau \mapsto \lambda_\tau$ gives the bijection from 
the set of equivalence classes of irreducible representations 
of $\mathrm{SO}(n)$ to the set of 
$\lambda =(\lambda_1,\lambda_2,\cdots ,\lambda_{m-1},\lambda_m)
\in \Lambda_m$ satisfying 
\[
\left\{\begin{array}{ll}
(\lambda_1,\lambda_2,\cdots ,\lambda_{m-1},-\lambda_m)\in \Lambda_m
&\text{if $n$ is even},\\
\lambda_m\geq 0&\text{if $n$ is odd}.
\end{array}\right.
\]
For $\lambda =(\lambda_1 ,\lambda_2,\cdots ,\lambda_{m-1},\lambda_m)
\in \Lambda_m$ such that $\lambda_m\geq 0$, 
we take a representation 
$(\tau_{\mathfrak{so}(n),\lambda},V_{\mathfrak{so}(n),\lambda})$ 
of $\mathrm{SO}(n)$ so that 
$\tau_{\mathfrak{so}(n),\lambda}$ is 
a direct sum of two 
irreducible representations  
with highest weights $\lambda $ and 
$(\lambda_1,\lambda_2,\cdots ,\lambda_{m-1},-\lambda_m)$
if $n=2m$ and $\lambda_m>0$, 
and $\tau_{\mathfrak{so}(n),\lambda}$ is 
an irreducible representation 
with highest weight $\lambda $ otherwise. 
By Weyl's dimension formula \cite[Theorem 4.48]{Knapp_002}, 
we have 
\begin{align}
\label{eqn:dim_On_rep}
\dim V_{\mathfrak{so}(n),
(i+1,\mathbf{1}_{h-1},\mathbf{0}_{m-h})}
=\frac{(2i+n)}{(i+l)(i+n-l)}
\frac{(i+n-1)!}{i!(n-1-l)!(l-1)!}
\end{align}
for $1\leq h\leq m$ and $i\in \bN_0$.

\begin{lem}
\label{lem:On_rep_1}
Retain the notation. Let $0\leq l\leq n$. As an $\rO (n)$-module, 
$V_{(\mathbf{1}_l,\mathbf{0}_{n-l})}$ is irreducible and 
$V_{(\mathbf{1}_l,\mathbf{0}_{n-l})} \not\simeq 
V_{(\mathbf{1}_{l'},\mathbf{0}_{n-l'})}$ 
for any $0\leq l'\leq n$ such that $l'\neq l$. 
We set $h=\min \{l,n-l\}$. 
When $n\geq 2$, we have 
\begin{align}
\label{eqn:lem_rest_On_001}
V_{(\mathbf{1}_l,\mathbf{0}_{n-l})}
\simeq V_{\gs \go (n),(\mathbf{1}_{h},\mathbf{0}_{m-h})}
\quad \text{as $\mathrm{SO}(n)$-modules}.
\end{align}
\end{lem}
\begin{proof}
For $1\leq i_1<i_2<\cdots <i_l\leq n$ and 
$\varepsilon_1,\varepsilon_2,\cdots ,\varepsilon_n\in \{\pm 1\}$, 
we have 
\begin{align*}
&\tau_{(\mathbf{1}_l,\mathbf{0}_{n-l})}
(\diag (\varepsilon_1,\varepsilon_2,\cdots ,\varepsilon_n))
\,\gge_{i_1}\wedge \gge_{i_2}\wedge \cdots \wedge \gge_{i_l}\\
&=\varepsilon_{i_1}\varepsilon_{i_2}\cdots \varepsilon_{i_l}
\,\gge_{i_1}\wedge \gge_{i_2}\wedge \cdots \wedge \gge_{i_l}.
\end{align*}
By this equality, we know that 
$\Hom_{\rO (n)}(V_{(\mathbf{1}_l,\mathbf{0}_{n-l})}, 
V_{(\mathbf{1}_{l'},\mathbf{0}_{n-l'})})=\{0\}$ 
for any $0\leq l'\leq n$ such that $l'\neq l$. 
Hence, our task is to show  
(\ref{eqn:lem_rest_On_001}) and the irreducibility of 
$V_{(\mathbf{1}_l,\mathbf{0}_{n-l})}$ 
as an $\rO (n)$-module.

In the case of $n\geq 2$ and $n\neq 2l$, the isomorphism 
(\ref{eqn:lem_rest_On_001}) follows from 
\cite[Examples in Chapter IV, \S 7]{Knapp_002}, 
and we note that 
$V_{(\mathbf{1}_l,\mathbf{0}_{n-l})}$ 
is an irreducible $\rO (n)$-module. 
In the case of $n=1$, 
the irreducibility of an $\rO (1)$-module 
$V_{(\mathbf{1}_l,\mathbf{0}_{1-l})}$ is trivial. 
Let us consider the case of $n=2l$. 
By direct computation, for $\varepsilon \in \{\pm 1\}$, 
we can confirm that 
\begin{align*}
v_\varepsilon =\,&(\gge_1+\sI \gge_2)\wedge (\gge_3+\sI \gge_4)\wedge \cdots \\
&\cdots \wedge (\gge_{n-3}+\sI \gge_{n-2})
\wedge (\gge_{n-1}+\varepsilon \sI \gge_n)
\end{align*}
is an $\mathrm{SO}(n)$-highest weight vector of weight 
$(\mathbf{1}_{l-1},\varepsilon )$ 
in $V_{(\mathbf{1}_l,\mathbf{0}_{l})}$ and satisfies 
$\tau_{(\mathbf{1}_l,\mathbf{0}_{n-l})}(k_0)v_\varepsilon =v_{-\varepsilon}$. 
Since $\dim V_{(\mathbf{1}_l,\mathbf{0}_{l})}
=\dim V_{\mathfrak{so}(n),\mathbf{1}_l}$ by (\ref{eqn:dim_On_rep}), 
we know that (\ref{eqn:lem_rest_On_001}) holds 
and $V_{(\mathbf{1}_l,\mathbf{0}_{l})}$ is 
an irreducible $\rO (n)$-module. 
\end{proof}

\begin{proof}[Proof of Lemma \ref{lem:rep_Kn_irred}]
The assertion for $F=\bR$ follows immediately from Lemma \ref{lem:On_rep_1}. 
The assertion for $F=\bC$ follows immediately from 
the highest weight theory \cite[Theorem 4.28]{Knapp_002} 
for $\rU (n)$. 
\end{proof}

\begin{lem}
Assume $n\geq 2$. 
Let $1\leq l\leq n-1$. We set $h=\min \{l,n-l\}$. 
Let $\rI_{l}\colon V_{(\mathbf{1}_{l-1},\mathbf{0}_{n-l+1})}\to 
V_{(\mathbf{1}_l,\mathbf{0}_{n-l})}\otimes_{\bC}
V_{(1,\mathbf{0}_{n-1})}$ be a $\bC$-linear map defined by 
\begin{align*}
&\rI_{l}(\gge_{i_1}\wedge \gge_{i_2}\wedge \cdots \wedge \gge_{i_{l-1}})
=\sum_{j=1}^n(\gge_{i_1}\wedge \gge_{i_2}\wedge \cdots \wedge \gge_{i_{l-1}}
\wedge \gge_j)\otimes \gge_j 
\end{align*}
for $i_1,i_2, \cdots ,i_{l-1}\in \{1,2,\cdots ,n\}$. 
Here we understand $\rI_{1}(1)
=\sum_{j=1}^n\gge_j\otimes \gge_j$ if $l=1$. 
Then $\rI_{l}$ is an $\rO (n)$-homomorphism. 
Moreover, there is an $\mathrm{SO}(n)$-submodule $V'$ 
of $V_{(\mathbf{1}_l,\mathbf{0}_{n-l})}
\otimes_\bC V_{(1,\mathbf{0}_{n-1})}$ such that 
$V'\simeq V_{\gs \go (n),(2,\mathbf{1}_{h-1},\mathbf{0}_{m-h})}$ and 
\begin{equation}
\label{eqn:decomp_On_rep}
\begin{aligned}
&V_{(\mathbf{1}_l,\mathbf{0}_{n-l})}
\otimes_\bC V_{(1,\mathbf{0}_{n-1})}\\
&=\rI^{(\mathbf{1}_{l},\mathbf{0}_{n-l}),1}_{
(\mathbf{1}_{l+1},\mathbf{0}_{n-l-1})}
(V_{(\mathbf{1}_{l+1},\mathbf{0}_{n-l-1})})
\oplus \rI_{l}(V_{(\mathbf{1}_{l-1},\mathbf{0}_{n-l+1})})
\oplus V',
\end{aligned}
\end{equation}
where $\rI^{(\mathbf{1}_{l},\mathbf{0}_{n-l}),1}_{
(\mathbf{1}_{l+1},\mathbf{0}_{n-l-1})}$ is 
the $\rU (n)$-homomorphism given by (\ref{eqn:injector_explicit}). 
\end{lem}
\begin{proof}
For $v\in V_{(\mathbf{1}_{l-1},\mathbf{0}_{n-l+1})}$ and $1\leq i<j\leq n$, 
we have 
\begin{align*}
&\rI_{l}
(\tau_{(\mathbf{1}_{l-1},\mathbf{0}_{n-l+1})}(E^{\gs \go (n)}_{i,j})v)
=(\tau_{(\mathbf{1}_{l},\mathbf{0}_{n-l})}\otimes \tau_{(1,\mathbf{0}_{n-1})})
(E^{\gs \go (n)}_{i,j})\rI_{l}(v),\\
&\rI_{l}
(\tau_{(\mathbf{1}_{l-1},\mathbf{0}_{n-l+1})}(k_0)v)
=(\tau_{(\mathbf{1}_{l},\mathbf{0}_{n-l})}\otimes \tau_{(1,\mathbf{0}_{n-1})})
(k_0)\rI_{l}(v)
\end{align*}
by direct computation. 
Hence, $\rI_{l}$ is an $\rO (n)$-homomorphism. 

For an $\mathrm{SO}(n)$-highest weight vector $v$ of weight $\lambda$ 
in $V_{(\mathbf{1}_l,\mathbf{0}_{n-l})}$, we note that 
$v\otimes (\gge_1+\sI \gge_2)$ is an $\mathrm{SO}(n)$-highest weight vector 
of weight $\lambda +(1,\mathbf{0}_{m-1})$ 
in $V_{(\mathbf{1}_l,\mathbf{0}_{n-l})}\otimes_\bC 
V_{(1,\mathbf{0}_{n-1})}$. 
Hence, by Lemma \ref{lem:On_rep_1}, 
there is an $\mathrm{SO}(n)$-submodule $V'$ 
of $V_{(\mathbf{1}_l,\mathbf{0}_{n-l})}
\otimes_\bC V_{(1,\mathbf{0}_{n-1})}$ such that 
$V'\simeq V_{\gs \go (n),(2,\mathbf{1}_{h-1},\mathbf{0}_{m-h})}$ and 
\[
\Bigl(
\rI^{(\mathbf{1}_{l},\mathbf{0}_{n-l}),1}_{
(\mathbf{1}_{l+1},\mathbf{0}_{n-l-1})}
(V_{(\mathbf{1}_{l+1},\mathbf{0}_{n-l-1})})
\oplus \rI_{l}(V_{(\mathbf{1}_{l-1},\mathbf{0}_{n-l+1})})\Bigr)\cap 
V'=\{0\}
\]
By (\ref{eqn:dim_On_rep}), we know that 
$\dim V_{(\mathbf{1}_l,\mathbf{0}_{n-l})}\otimes_\bC V_{(1,\mathbf{0}_{n-1})}$ 
is equal to 
\[
\dim V_{(\mathbf{1}_{l+1},\mathbf{0}_{n-l-1})}
+\dim V_{(\mathbf{1}_{l-1},\mathbf{0}_{n-l+1})}
+\dim V_{\gs \go (n),(2,\mathbf{1}_{h-1},\mathbf{0}_{m-h})}.
\]
This implies that (\ref{eqn:decomp_On_rep}) holds. 
\end{proof}

\begin{lem}
\label{lem:isom_pairing}
Let $(\tau,V_\tau )$ and 
$(\tau',V_{\tau'})$ be finite dimensional representations 
of $\mathrm{GL}(n,\bC )$. 
Let $\langle \cdot ,\cdot \rangle$ be 
a $\rU (n)$-invariant hermitian inner product on $V_\tau $. 
Let $\{v_i\}_{i=1}^d$ be an orthonormal basis of $V_{\tau}$. 
\vspace{1mm}

\noindent (1) A $\bC$-linear map 
$\Psi_1\colon \Hom_{\bC}(V_{\tau'},V_{\tau})\to 
\Hom_{\bC}(V_{\tau'} \otimes_\bC \overline{V_{\tau}},\bC_{\mathrm{triv}})$ 
defined by 
\begin{align*}
\Psi_1(f)(v'\otimes \overline{v})=\langle f (v'),v\rangle&
&(f\in \Hom_{\bC }(V_{\tau'},V_{\tau}),\ v'\in V_{\tau'},\ v\in V_{\tau})
\end{align*}
is bijective, and its inverse map is given by  
\begin{align*}
\Psi_1^{-1}(f)(v')=\sum_{i=1}^d
f(v'\otimes \overline{v_i})\,v_i&
&(f\in 
\Hom_{\bC}(V_{\tau'} \otimes_\bC \overline{V_{\tau}},\bC_{\mathrm{triv}}),\ 
v'\in V_{\tau'}).
\end{align*}
Moreover, we have 
$\Psi_1(\Hom_{K_n}(V_{\tau'},V_{\tau}))=
\Hom_{K_n}(V_{\tau'} \otimes_\bC \overline{V_{\tau}},\bC_{\mathrm{triv}})$. 
\vspace{1mm}

\noindent 
(2) A $\bC$-linear map 
$\Psi_2\colon V_{\tau'} \otimes_\bC \overline{V_{\tau}}\to 
\Hom_{\bC}(V_{\tau},V_{\tau'})$ defined by 
\begin{align*}
&\Psi_2(v'\otimes \overline{v_1})(v_2)
=\langle v_2,v_1\rangle v'&
&(v_1,v_2\in V_{\tau},\ v'\in V_{\tau'})
\end{align*}
is bijective, and its inverse map is given by 
\begin{align*}
\Psi_2^{-1}(f)=\sum_{i=1}^d
f(v_i)\otimes \overline{v_i}&
&(f\in \Hom_{\bC}(V_{\tau},V_{\tau'})).
\end{align*} 
Moreover, we have $\Phi_2((V_{\tau'} \otimes_\bC \overline{V_{\tau}})^{K_n})
=\Hom_{K_n}(V_{\tau},V_{\tau'})$, where 
$(V_{\tau'} \otimes_\bC \overline{V_{\tau}})^{K_n}$ is 
the subspace of $V_{\tau'} \otimes_\bC \overline{V_{\tau}}$ 
consisting of all $K_n$-invariant vectors. 
\end{lem}
\begin{proof}
The former part of the statement (1) follows from definition. 
The latter part of the statement (1) follows from 
\begin{align*}
\Psi_1(f)((\tau' \otimes \overline{\tau})(k)v'\otimes \overline{v})
&=\langle f (\tau'(k)v'),\tau (k)v\rangle 
=\langle \tau (k^{-1})f (\tau'(k)v'),v\rangle \\
&=\Psi_1(\tau (k^{-1})\circ f\circ \tau'(k))(v'\otimes \overline{v})
\end{align*}
for $f\in \Hom_{\bC }(V_{\tau'},V_{\tau})$, $v'\in V_{\tau'}$, 
$v\in V_{\tau}$ and $k\in K_n$. 

The former part of the statement (2) follows from definition. 
The latter part of the statement (2) follows from 
\begin{align*}
\Psi_2((\tau' \otimes \overline{\tau})(k)v'\otimes \overline{v_1})(v_2)
&=\langle v_2,\tau (k)v_1\rangle \tau' (k)v'
=\langle \tau (k^{-1})v_2 ,v_1 \rangle \tau' (k)v'\\
&=\tau' (k)
\Psi_2(v'\otimes \overline{v_1})(\tau (k^{-1})v_2 )
\end{align*}
for $v_1 ,v_2 \in V_{\tau }$, $v'\in V_{\tau'}$ and $k\in K_n$. 
\end{proof}

\begin{proof}[Proof of Lemma \ref{lem:schur_pairing}]
Let $\lambda \in \Lambda_{n,F}$. 
By Lemma \ref{lem:rep_Kn_irred}, we note that the space 
$\Hom_{K_n}(V_\lambda ,V_{\lambda})$ is a $1$ dimensional space spanned 
by the identity map. Hence, we obtain 
the statement (1) by Lemma \ref{lem:isom_pairing} (1). 

Let $\lambda'\in \Lambda_{n,F}\cap \Xi^\circ (\lambda )$, 
and we set $l=\ell (\lambda' -\lambda )$. 
By the decompositions (\ref{eqn:pieri}), (\ref{eqn:decomp_On_rep}) 
and Lemma \ref{lem:On_rep_1}, we have 
$\Hom_{K_n}(V_{\lambda'} ,V_{\lambda''})=\{0\}$ for 
$\lambda''\in \Xi^\circ (\lambda ;l)$ such that $\lambda''\neq \lambda'$. 
Hence, we obtain the statement (2) by Lemma \ref{lem:isom_pairing} (1). 
\end{proof}

\begin{proof}[Proof of Lemma \ref{lem:schur_pairing+}]
By (\ref{eqn:decomp_n_n-1}) and Lemma \ref{lem:rep_Kn_irred}, 
we know that $\Hom_{K_{n-1}}(V_\lambda ,V_{\mu})$ is equal 
to $\bC \,\rR_{\mu}^{\lambda}$ if $\mu \in \Xi^+(\lambda )$, 
and is equal to $\{0\}$ otherwise. 
By Lemma \ref{lem:isom_pairing}, 
we obtain the former part of the assertion, 
and know that, if $\mu \in \Xi^+(\lambda )$,  
\begin{equation*}
\sum_{M \in \rG (\mu )}
\rR_{\mu}^{\lambda}(\zeta_{M[\lambda ]})\otimes \overline{\zeta_{M[\lambda ]}}
\end{equation*}
is a unique $K_{n-1}$-invariant vector in 
$V_{\mu }\otimes_\bC 
\overline{V_{\lambda}}$ up to scalar multiple. 
Hence, by (\ref{eqn:def_Qbasis}) and the properties of complex conjugate 
representations in \S \ref{subsec:com_conj_rep}, 
we obtain the latter part of the assertion. 
\end{proof}

\begin{proof}[Proof of Lemma \ref{lem:schur_pairing0}]
By the decompositions (\ref{eqn:pieri}), (\ref{eqn:decomp_On_rep}) 
and Lemma \ref{lem:On_rep_1}, we know that 
$\Hom_{K_{n}}(V_{\lambda'},V_{\lambda}\otimes_\bC V_{(l,\mathbf{0}_{n-1})})$ 
is equal to $\bC\,\mathrm{I}^{\lambda ,l}_{\lambda'}$ 
if $\lambda'\in \Xi^\circ (\lambda ;l)$, 
and is equal to $\{0\}$ otherwise. 
By Lemma \ref{lem:isom_pairing}, 
we obtain the former part of the assertion, 
and know that, if $\lambda'\in \Xi^\circ (\lambda ;l)$,  
\begin{equation*}
\sum_{M'\in \rG (\lambda')}
\rI^{\lambda,l}_{\lambda'}(\zeta_{M'})\otimes \overline{\zeta_{M'}}
\end{equation*}
is a unique $K_{n}$-invariant vector in 
$V_{\lambda}
\otimes_\bC V_{(l,\mathbf{0}_{n-1})} \otimes_\bC 
\overline{V_{\lambda'}}$ up to scalar multiple. 
Hence, by (\ref{eqn:def_Qbasis}) and the properties of complex conjugate 
representations in \S \ref{subsec:com_conj_rep}, 
we obtain the latter part of the assertion. 
\end{proof}

\subsection{Polynomial functions}

We set 
\[
\Lambda_n^{\mathrm{poly}}=\{\lambda =
(\lambda_1,\lambda_2,\cdots ,\lambda_n)\in \Lambda_n\mid 
\lambda_n\geq 0\}. 
\]
We denote by $\cP (\rM_{n,n'}(\bC))$ the subspace of 
$C(\rM_{n,n'}(\bC))$ consisting of all polynomial functions. 
Let $l\in \bN_0$. We denote by 
$\cP_l (\rM_{n,n'}(\bC))$ the subspace of 
$\cP (\rM_{n,n'}(\bC))$ consisting of all 
degree $l$ homogeneous polynomial functions. 
We regard $\cP_l (\rM_{n,n'}(\bC))$ as a 
$\mathrm{GL}(n,\bC )\times \mathrm{GL}(n',\bC )$-module 
via the action $L\boxtimes R$ which is defined in \S 
\ref{subsec:def_schwartz}. 
Let $q=\min \{n,n'\}$. 
Then the $\mathrm{GL}(n)$-$\mathrm{GL}(n')$ duality 
\cite[Theorem 5.6.7]{Goodman_Wallach_001} asserts that 
\begin{align*}
&\cP_l(\rM_{n,n'}(\bC))\simeq 
\bigoplus_{\lambda \in \Lambda_q^{\mathrm{poly}},\ \ell (\lambda )=l}
V_{(\lambda ,\mathbf{0}_{n-q})}^\vee \boxtimes_{\bC}
V_{(\lambda ,\mathbf{0}_{n'-q})}
\end{align*}
as $\mathrm{GL}(n,\bC )\times \mathrm{GL}(n',\bC )$-modules. 
Since $V_{(\lambda ,\mathbf{0}_{n-q})}^\vee \simeq 
\overline{V_{(\lambda ,\mathbf{0}_{n-q})}}$ as $\rU (n)$-modules, 
we also have 
\begin{align}
\label{eqn:K_GL_duality}
&\cP_l(\rM_{n,n'}(\bC))\simeq 
\bigoplus_{\lambda \in \Lambda_q^{\mathrm{poly}},\ \ell (\lambda )=l}
\overline{V_{(\lambda ,\mathbf{0}_{n-q})}} \boxtimes_{\bC}
V_{(\lambda ,\mathbf{0}_{n'-q})}
\end{align}
as $\rU (n)\times \mathrm{GL}(n',\bC )$-modules. 

The purpose of this subsection is to 
construct polynomial functions, explicitly.
We define 
$\rU (n)\times \mathrm{GL}(n,\bC )$-homomorphisms 
$\rP^\circ_{\lambda} \colon 
\overline{V_{\lambda}}\boxtimes_\bC V_{\lambda}
\to \cP (\rM_{n}(\bC ))$ ($\lambda \in \Lambda_n^{\mathrm{poly}}$) 
by the following lemma.

\begin{lem}
\label{lem:K_GLn_poly1}
Let $\lambda \in \Lambda_n^{\mathrm{poly}}$. Then 
there is a $\rU (n)\times \mathrm{GL}(n,\bC )$-homomorphism 
$\rP^\circ_{\lambda} \colon 
\overline{V_{\lambda}}\boxtimes_\bC V_{\lambda}
\to \cP (\rM_{n}(\bC ))$ characterized by 
\begin{align}
\label{eqn:Kn_GLn_poly_001}
&\rP^\circ_{\lambda}(\overline{v_1}\boxtimes v_2)(g)
=\langle \tau_{\lambda}(g)v_2,v_1 \rangle&
&(v_1,v_2\in V_{\lambda},\ g\in \mathrm{GL}(n,\bC )).
\end{align}
\end{lem}
\begin{proof}
Because of (\ref{eqn:K_GL_duality}), 
there is a nonzero $\rU (n)\times \mathrm{GL}(n,\bC )$-homomorphism 
$\rP \colon 
\overline{V_{\lambda}}\boxtimes_\bC V_{\lambda}
\to \cP (\rM_{n}(\bC ))$. 
Since $\mathrm{GL}(n,\bC )$ is dense in 
$\rM_{n}(\bC)$ and 
\begin{align}
\label{eqn:pf_lem_K_GL_001}
&\rP (\overline{v_1}\boxtimes v_2)(g)
=\rP (\overline{v_1}\boxtimes \tau_{\lambda}(g)v_2)(1_{n})&
&(v_1,v_2\in V_{\lambda},\ g\in \mathrm{GL}(n,\bC )), 
\end{align}
we note that 
\[
\overline{V_{\lambda}}\boxtimes_{\bC}V_{\lambda}
\ni v_1\boxtimes v_2\mapsto 
\rP (v_1\boxtimes v_2)(1_n)\in \bC 
\]
is a nonzero $\bC$-bilinear pairing. Because of 
\begin{align*}
&\rP (\overline{\tau_\lambda}(k)\overline{v_1}
\boxtimes \tau_{\lambda}(k)v_2)(1_{n})
=\rP (\overline{v_1}\boxtimes v_2)(1_{n})&
&(k\in \rU (n))
\end{align*}
and Lemma \ref{lem:schur_pairing} (1) for $F=\bC$, 
there is a nonzero constant $c$ such that  
\begin{align}
\label{eqn:pf_lem_K_GL_002}
&\rP (\overline{v_1}\boxtimes v_2)(1_{n})
=c\, \langle v_2,v_1\rangle&
&(v_1,v_2\in V_{\lambda}).
\end{align}
By (\ref{eqn:pf_lem_K_GL_001}) and (\ref{eqn:pf_lem_K_GL_002}), 
we know that $\rP^\circ_\lambda =c^{-1}\rP$ satisfies 
(\ref{eqn:Kn_GLn_poly_001}). 
Since $\mathrm{GL}(n,\bC )$ is dense in $\rM_{n}(\bC)$, 
we note that (\ref{eqn:Kn_GLn_poly_001}) characterizes 
$\rP^\circ_\lambda $. 
\end{proof}

When $n>1$, we define $\rU (n-1)\times \mathrm{GL}(n,\bC )$-homomorphisms 
$\rP^+_{\mu} \colon 
\overline{V_{\mu}}\boxtimes_\bC V_{(\mu,0)}
\to \cP (\rM_{n-1,n}(\bC ))$ ($\mu \in \Lambda_{n-1}^{\mathrm{poly}}$) 
by the following lemma.

\begin{lem}
\label{lem:K_GLn_poly2}
Assume $n>1$ and let  
$\mu \in \Lambda_{n-1}^{\mathrm{poly}}$. 
There is a $\rU (n-1)\times \mathrm{GL}(n,\bC )$-homomorphism 
$\rP^+_{\mu} \colon 
\overline{V_{\mu}}\boxtimes_\bC V_{(\mu,0)}
\to \cP (\rM_{n-1,n}(\bC ))$ characterized by 
\begin{align}
\label{eqn:Kn_GLn_poly_002}
\rP^+_{\mu } (\overline{\zeta_{M}}\boxtimes v)
((1_{n-1},O_{n-1,1})z)
=\rP^\circ_{(\mu ,0)} (\overline{\zeta_{M[(\mu,0)]}}
\boxtimes v)(z)&
\end{align}
for $M\in \rG(\mu )$, $v\in V_{(\mu ,0)}$ and 
$z\in \rM_{n}(\bC)$. Here 
$M[(\mu,0)]$ is defined by (\ref{eqn:def_Mlambda}). 
Furthermore, we have 
\begin{align}
\label{eqn:Kn_GLn_poly_003}
&\rP^+_{\mu }
(\overline{v}\boxtimes \zeta_{M[(\mu,0)]})(z)=
\rP^\circ_{\mu}(\overline{v}\boxtimes \zeta_{M})
(z\,{}^{t\,}\!(1_{n-1},O_{n-1,1}))
\end{align}
for $v\in V_\mu $, $M\in \rG(\mu )$ and 
$z\in \rM_{n-1,n}(\bC)$. 
\end{lem}
\begin{proof}
We regard $\mathrm{GL}(n-1,\bC )$ as a subgroup of 
$\mathrm{GL}(n,\bC )$ via 
the embedding $\iota_n$ defined by (\ref{eqn:embed_Gn-1_Gn}). 
By the irreducible decomposition (\ref{eqn:K_GL_duality}) 
and Lemma \ref{lem:unitary_isom_lambda_mu}, the image of a 
$\rU(n-1)\times \mathrm{GL}(n,\bC )$-homomorphism
\begin{align*}
&\overline{V_\mu}\boxtimes_\bC V_{(\mu,0)}
\ni \overline{\zeta_{M}}\boxtimes v\mapsto 
\rP^\circ_{(\mu ,0)}
\bigl(\overline{\zeta_{M[(\mu,0)]}}
\boxtimes v\bigr)\in \cP(\rM_{n}(\bC))
\end{align*}
is contained in the image of 
an injective $\rU(n-1)\times \mathrm{GL}(n,\bC )$-homomorphism
\begin{align*}
&\cP(\rM_{n-1,n}(\bC))\ni p\mapsto \bigl(z\mapsto 
p((1_{n-1},O_{n-1,1})z)\bigr)\in \cP(\rM_{n}(\bC)). 
\end{align*}
Hence, there is a 
$\rU (n-1)\times \mathrm{GL}(n,\bC )$-homomorphism 
$\rP^+_{\mu} \colon 
\overline{V_{\mu}}\boxtimes_\bC V_{(\mu,0)}
\to \cP (\rM_{n-1,n}(\bC ))$ characterized by (\ref{eqn:Kn_GLn_poly_002}). 
By the irreducible decompositions (\ref{eqn:K_GL_duality}) and 
Lemma \ref{lem:unitary_isom_lambda_mu}, two injective 
$\rU(n-1)\times \mathrm{GL}(n-1,\bC )$-homomorphisms
\begin{align*}
&\overline{V_\mu}\boxtimes_\bC V_{\mu }
\ni \overline{v}\boxtimes \zeta_{M}
\mapsto 
\rP^+_{\mu}\bigl(\overline{v}\boxtimes \zeta_{M[(\mu ,0)]}\bigr)
\in \cP(\rM_{n-1,n}(\bC)),\\[1mm]
&\overline{V_\mu}\boxtimes_\bC V_{\mu }
\ni \overline{v_1}\boxtimes v_2\mapsto \bigl(z\mapsto 
\rP^\circ_{\mu}(\overline{v_1}\boxtimes v_2)
(z\,{}^t\!(1_{n-1},O_{n-1,1}))\bigr)\in \cP(\rM_{n-1,n}(\bC))
\end{align*}
coincide up to scalar multiple. Hence, 
(\ref{eqn:Kn_GLn_poly_003}) follows from the equalities 
\begin{align*}
\rP^+_{\mu}
\bigl(\overline{\zeta_{M}}\boxtimes \zeta_{M[(\mu,0)]}
\bigr)((1_{n-1},O_{n-1,1}))
&=\langle \zeta_{M[(\mu,0)]},\zeta_{M[(\mu,0)]}\rangle =1
\end{align*}
and $\rP^\circ_{\mu}(\overline{\zeta_{M}}
\boxtimes \zeta_{M})(1_{n-1})
=\langle \zeta_{M},\zeta_{M}\rangle =1$ 
for $M\in \rG(\mu )$. 
\end{proof}

Let $l\in \bN_0$. We define $\bC$-linear maps 
$\rp^{(l)}_{1,n}\colon V_{(l,\mathbf{0}_{n-1})}\to 
\cP (\rM_{1,n}(\bC ))$ and 
$\rp^{(l)}_{n,1}\colon \overline{V_{(l,\mathbf{0}_{n-1})}}\to 
\cP (\rM_{n,1}(\bC ))$ by 
\begin{align*}
&\rp^{(l)}_{1,n}(\zeta_{Q(\gamma )})(z)
=\rp^{(l)}_{n,1}(\overline{\zeta_{Q(\gamma)}})({}^t\!z)=
\sqrt{\rb (\gamma )}\,z_1^{\gamma_1}z_2^{\gamma_2}\cdots z_n^{\gamma_n}.
\end{align*}
for $z=(z_1,z_2,\cdots ,z_n)\in \rM_{1,n}(\bC )$ and 
$\gamma=(\gamma_1,\gamma_2,\cdots ,\gamma_n)\in \bN_0^n$ such that 
$\ell (\gamma )=l$. Here $Q(\gamma)$ and $\rb (\gamma )$ are 
defined by (\ref{eqn:def_Qgamma}) and (\ref{eqn:def_c_gamma}), 
respectively.

\begin{lem}
\label{lem:poly_n1_1n}
Let $l\in \bN_0$.  
\vspace{1mm}

\noindent (1) The group $\mathrm{GL}(n,\bC )$ acts on 
$\cP (\rM_{1,n}(\bC ))$ by the action $R$. 
Then $\rp^{(l)}_{1,n}$ 
is a $\mathrm{GL}(n,\bC )$-homomorphism such that, 
for $z\in \rM_{n}(\bC )$ and $v\in V_{(l,\mathbf{0}_{n-1})}$, 
\begin{align}
\label{eqn:poly_mc_1_n}
&\rp^{(l)}_{1,n}(v)(e_nz)
=\rP^\circ_{(l,\mathbf{0}_{n-1})}(\overline{\zeta_{Q((\mathbf{0}_{n-1},l))}}
\boxtimes v)(z).
\end{align}

\noindent (2) The group $\rU (n)$ acts on $\cP (\rM_{n,1}(\bC ))$ by 
the action $L$. 
Then $\rp^{(l)}_{n,1}$ 
is a $\rU (n)$-homomorphism such that, 
for $z\in \rM_{n}(\bC )$ and $v\in V_{(l,\mathbf{0}_{n-1})}$, 
\begin{align}
\label{eqn:poly_mc_n_1}
&\rp^{(l)}_{n,1}(\overline{v})(z\,{}^t\!e_n)
=\rP^\circ_{(l,\mathbf{0}_{n-1})}(\overline{v}
\boxtimes \zeta_{Q((\mathbf{0}_{n-1},l))})(z).
\end{align}
\end{lem}
\begin{proof}
Let $\gamma=(\gamma_1,\gamma_2,\cdots ,\gamma_n)\in \bN_0^n$. 
By direct computation, we have 
\begin{align*}
&R(E_{i,i})\rp^{(l)}_{1,n}(\zeta_{Q(\gamma )})
=\gamma_i\rp^{(l)}_{1,n}(\zeta_{Q(\gamma )}),\\
&R(E_{j,j+1})\rp^{(l)}_{1,n}(\zeta_{Q(\gamma )})
=\sqrt{\gamma_{j+1}(\gamma_j+1)}\,\rp^{(l)}_{1,n}
(\zeta_{Q(\gamma +\delta_j-\delta_{j+1})}),\\
&R(E_{j+1,j})\rp^{(l)}_{1,n}(\zeta_{Q(\gamma )})
=\sqrt{\gamma_j(\gamma_{j+1}+1)}\,
\rp^{(l)}_{1,n}(\zeta_{Q(\gamma -\delta_j+\delta_{j+1})}),\\
&L(E_{i,i})\rp^{(l)}_{n,1}(\overline{\zeta_{Q(\gamma )}})
=-\gamma_i\rp^{(l)}_{n,1}(\overline{\zeta_{Q(\gamma )}}),\\
&L(E_{j,j+1})\rp^{(l)}_{n,1}(\overline{\zeta_{Q(\gamma )}})
=-\sqrt{\gamma_j(\gamma_{j+1}+1)}\,
\rp^{(l)}_{n,1}(\overline{\zeta_{Q(\gamma -\delta_j+\delta_{j+1})}}),\\
&L(E_{j+1,j})\rp^{(l)}_{n,1}(\overline{\zeta_{Q(\gamma )}})
=-\sqrt{\gamma_{j+1}(\gamma_j+1)}\,
\rp^{(l)}_{n,1}(\overline{\zeta_{Q(\gamma +\delta_j-\delta_{j+1})}})
\end{align*}
for $1\leq i\leq n$ and $1\leq j\leq n-1$. 
Here we put $\rp^{(l)}_{1,n}(\zeta_{Q(\gamma')})
=\rp^{(l)}_{n,1}(\overline{\zeta_{Q(\gamma' )}})=0$ 
if $\gamma' \not\in \bN_0^n$, and denote by 
$\delta_i$ the element of $\bZ^n$ 
with $1$ at $i$-th entry ($1\leq i\leq n$) and $0$ at other entries. 
Comparing these formulas with 
(\ref{eqn:GT_act_wt}), (\ref{eqn:GT_act+}) and (\ref{eqn:GT_act-}), 
we know that $\rp^{(l)}_{1,n}$ 
is a $\mathrm{GL}(n,\bC )$-homomorphism. 
Using (\ref{eqn:conconj_Liealg_rep}) and 
\begin{align*}
&L(E_{i,j}^{\gu (n)})=L(E_{i,j})\quad  (1\leq i,j\leq n)\qquad  
\text{on}\ \cP (\rM_{n,1}(\bC )), 
\end{align*}
we further know that $\rp^{(l)}_{n,1}$ 
is a $\rU (n)$-homomorphism.

Next, we will prove the equality (\ref{eqn:poly_mc_1_n}). 
When $n=1$, this equality follows from 
$\rG (l)=\{l\}$ and 
\begin{align*}
&\rp^{(l)}_{1,1}(\zeta_{Q(\gamma )})(g)
=\langle \tau_{l}(g)\zeta_{l},\zeta_{l}\rangle =g^l&
&(g\in \mathrm{GL}(1,\bC )).
\end{align*}
Assume $n>1$. We regard $\mathrm{GL}(n-1,\bC )$ as a subgroup of 
$\mathrm{GL}(n,\bC )$ via 
the embedding $\iota_n$ defined by (\ref{eqn:embed_Gn-1_Gn}). 
Because of the irreducible decompositions (\ref{eqn:K_GL_duality}) and 
Lemma \ref{lem:unitary_isom_lambda_mu}, 
two injective $\rU(n-1)\times \mathrm{GL}(n,\bC )$-homomorphisms
\begin{align*}
&\overline{V_{\mathbf{0}_{n-1}}}\boxtimes_\bC 
V_{(l,\mathbf{0}_{n-1})}
\ni \overline{\zeta_{Q(\mathbf{0}_{n-1})}}\boxtimes v\mapsto 
(z\mapsto \rp^{(l)}_{1,n}(v)(e_nz))
\in \cP(\rM_{n}(\bC)),\\[1mm]
&\overline{V_{\mathbf{0}_{n-1}}}\boxtimes_\bC 
V_{(l,\mathbf{0}_{n-1})}\ni \overline{\zeta_{Q(\mathbf{0}_{n-1})}}
\boxtimes v\mapsto 
\rP^\circ_{(l,\mathbf{0}_{n-1})}
(\overline{\zeta_{Q((\mathbf{0}_{n-1},l))}}\boxtimes v)\in \cP(\rM_{n}(\bC))
\end{align*}
coincide up to scalar multiple. 
Hence, (\ref{eqn:poly_mc_1_n}) follows from the equalities 
\begin{align*}
&\rp^{(l)}_{1,n}(\zeta_{Q((\mathbf{0}_{n-1},l))})(e_n)=1,&
&\rP^\circ_{(l,\mathbf{0}_{n-1})}
(\overline{\zeta_{Q((\mathbf{0}_{n-1},l))}}\boxtimes 
\zeta_{Q((\mathbf{0}_{n-1},l))})(1_n)=1.
\end{align*}
The proof of the equality 
(\ref{eqn:poly_mc_n_1}) is similar. 
\end{proof}

For later use, we prepare the following lemmas. 

\begin{lem}
\label{lem:poly_+_1n_decomp}
Assume $n>1$. Let 
$\mu \in \Lambda_{n-1}^{\mathrm{poly}}$ and 
$\gamma \in \bN_0^n$. 
We set $l=\ell (\gamma )$ and 
\begin{align*}
&p_0(z)=\rP_{\mu }^+\bigl(
\overline{\zeta_{H(\mu )}}
\boxtimes \zeta_{H((\mu ,0))}\bigr)((1_{n-1},O_{n-1,1})z)
\rp^{(l)}_{1,n}(\zeta_{Q(\gamma )})(e_nz)
\end{align*}
for $z\in \rM_n(\bC)$. 
Then we have 
\begin{align*}
&p_0=
\sum_{\lambda'\in \Xi^\circ ((\mu ,0);l)}\,
\sum_{N,N'\in \rG (\lambda')}
\rC^{H((\mu ,0)),Q((\mathbf{0}_{n-1},l))}_{N}\,
{\rC^{H((\mu ,0)),Q(\gamma )}_{N'}}\,
\rP^\circ_{\lambda'}(\overline{\zeta_{N}}\boxtimes \zeta_{N'}),
\end{align*}
where $\rC^{M,P}_{M'}$ is the Clebsch--Gordan coefficient in 
\S \ref{subsec:CG_coeff}. 
\end{lem}
\begin{proof}
We set $Q_0=Q((\mathbf{0}_{n-1},l))$ and $Q_1=Q(\gamma )$. 
Let $g\in \mathrm{GL}(n,\bC )$. 
By Lemmas \ref{lem:K_GLn_poly1}, \ref{lem:K_GLn_poly2} 
and \ref{lem:poly_n1_1n}, we have 
\begin{align*}
p_0(g)=&\rP^\circ_{\mu }\bigl(
\overline{\zeta_{H((\mu ,0))}}
\boxtimes \zeta_{H((\mu ,0))}\bigr)(g)
\rP^\circ_{(l,\mathbf{0}_{n-1})}(\overline{\zeta_{Q_0}}
\boxtimes \zeta_{Q_1})(g)\\
=&\bigl\langle \tau_{(\mu ,0)}(g)\zeta_{H((\mu ,0))},
\zeta_{H((\mu ,0))}\bigr\rangle\,
\bigl\langle \tau_{(l,\mathbf{0}_{n-1})}(g)\zeta_{Q_1},\ 
\zeta_{Q_0}\bigr\rangle \\
=&\bigl\langle 
(\tau_{(\mu ,0)}\otimes \tau_{(l,\mathbf{0}_{n-1})})
(g)\zeta_{H((\mu ,0))}\otimes \zeta_{Q_1},\ 
\zeta_{H((\mu ,0))}\otimes \zeta_{Q_0}\bigr\rangle . 
\end{align*}
By (\ref{eqn:Inj_normalize}), (\ref{eqn:decomp_pure_tensor}) 
and Lemma \ref{lem:K_GLn_poly1}, 
we have 
\begin{align*}
p_0(g)
=&
\sum_{\lambda' \in \Xi^\circ ((\mu ,0);l)}
\,\sum_{N,N'\in \rG (\lambda')}
\rC^{H((\mu ,0)),Q_0}_{N}\, 
{\rC^{H((\mu ,0)),Q_1}_{N'}}\,
\bigl\langle \tau_{\lambda'}(g)\zeta_{N'},\zeta_N\bigr\rangle \\
=&
\sum_{\lambda' \in \Xi^\circ ((\mu ,0);l)}
\,\sum_{N,N'\in \rG (\lambda')}
\rC^{H((\mu ,0)),Q_0}_{N}\, 
{\rC^{H((\mu ,0)),Q_1}_{N'}}\,
\rP^\circ_{\lambda'}(\overline{\zeta_{N}}\boxtimes \zeta_{N'})(g). 
\end{align*}
Since $\mathrm{GL}(n,\bC )$ is dense in $\rM_n(\bC)$, 
we obtain the assertion. 
\end{proof}

\begin{lem}
\label{lem:poly_0_n1_decomp}
Assume $n>1$. 
Let $\mu \in \Lambda_{n-1}^{\mathrm{poly}}$ and 
$\gamma \in \bN_0^{n-1}$. 
We set $l=\ell (\gamma )$ 
and 
\begin{align*}
&p_0(z)=\rP^\circ_{\mu }\bigl(
\overline{\zeta_{H(\mu )}}
\boxtimes \zeta_{H(\mu )}\bigr)(z\,{}^t\!(1_{n-1},O_{n-1,1}))
\rp^{(l)}_{n-1,1}(\overline{\zeta_{Q(\gamma )}})(z\,{}^t\!e_n)
\end{align*}
for $z\in \rM_{n-1,n}(\bC)$. 
Then we have 
\begin{align*}
&p_0=
\sum_{\mu' \in \Xi^\circ (\mu ;l)}
\,\underset{N'\in \rG ((\mu',0))}
{\sum_{N\in \rG ((\mu',0);\mu')}}
\rC^{H((\mu ,0)),Q((\gamma ,0))}_{N}\,
{\rC^{H((\mu ,0)),Q((\mathbf{0}_{n-1},l))}_{N'}}\, 
\rP^+_{\mu'}
\bigl(\overline{\zeta_{\widehat{N}}}\boxtimes \zeta_{N'}),
\end{align*}
where $\rC^{M,P}_{M'}$ is the Clebsch--Gordan coefficient in 
\S \ref{subsec:CG_coeff}. 
\end{lem}
\begin{proof}
We set $Q_0=Q((\mathbf{0}_{n-1},l))$ and $Q_1=Q((\gamma ,0))$. 
Let $z=(1_{n-1},O_{n-1,1})g$ with $g\in G_n$. 
Then we have 
$\rp^{(l)}_{n-1,1}(\overline{\zeta_{Q(\gamma )}})(z\,{}^t\!e_n)
=\rp^{(l)}_{n,1}(\overline{\zeta_{Q_1}})(g\,{}^t\!e_n)$ by definition. 
Hence, by Lemmas \ref{lem:K_GLn_poly1}, \ref{lem:K_GLn_poly2} 
and \ref{lem:poly_n1_1n}, we have 
\begin{align*}
p_0(z)=&\rP^\circ_{(\mu ,0)}\bigl(
\overline{\zeta_{H((\mu ,0))}}
\boxtimes \zeta_{H((\mu ,0))}\bigr)(g)
\rP^\circ_{(l,\mathbf{0}_{n-1})}
(\overline{\zeta_{Q_1}}\boxtimes \zeta_{Q_0})(g)\\
=&\bigl\langle \tau_{(\mu ,0)}(g)\zeta_{H((\mu ,0))},
\zeta_{H((\mu ,0))}\bigr\rangle\,
\bigl\langle \tau_{(l,\mathbf{0}_{n-1})}(g)\zeta_{Q_0},\ 
\zeta_{Q_1}\bigr\rangle \\
=&\bigl\langle 
(\tau_{(\mu ,0)}\otimes \tau_{(l,\mathbf{0}_{n-1})})
(g)\zeta_{H((\mu ,0))}\otimes \zeta_{Q_0},\ 
\zeta_{H((\mu ,0))}\otimes \zeta_{Q_1}\bigr\rangle . 
\end{align*}
By (\ref{eqn:Inj_normalize}) and (\ref{eqn:decomp_pure_tensor}), 
we have 
\begin{align*}
p_0(z)
=&
\sum_{\lambda' \in \Xi^\circ ((\mu ,0);l)}
\,\sum_{N,N'\in \rG (\lambda')}
\rC^{H((\mu ,0)),Q_1}_{N}\,
{\rC^{H((\mu ,0)),Q_0}_{N'}}\, 
\bigl\langle \tau_{\lambda'}(g)\zeta_{N'},\zeta_{N}\bigr\rangle . 
\end{align*}
Because of $H((\mu ,0))\in \rG ((\mu ,0);\mu )$, 
$Q_1\in \rG ((l,\mathbf{0}_{n-1});(l,\mathbf{0}_{n-2}))$ and 
(\ref{eqn:coef_decomp}), for 
$\lambda' \in \Xi^\circ ((\mu ,0);l)$ and $N\in \rG (\lambda')$, 
we have $\rC^{H((\mu ,0)),Q_1}_{N}=0$ unless 
$\lambda'=(\mu',0)$ and 
$N\in \rG ((\mu',0);\mu')$ with some 
$\mu'\in \Xi^\circ (\mu ;l)$. Hence, we have 
\begin{align*}
p_0(z)
=&
\sum_{\mu' \in \Xi^\circ (\mu ;l)}
\,\underset{N'\in \rG ((\mu',0))}
{\sum_{N\in \rG ((\mu',0);\mu')}}
\rC^{H((\mu ,0)),Q_1}_{N}\,
{\rC^{H((\mu ,0)),Q_0}_{N'}}\, 
\bigl\langle \tau_{(\mu',0)}(g)\zeta_{N'},\zeta_{N}\bigr\rangle \\
=&
\sum_{\mu' \in \Xi^\circ (\mu ;l)}
\,\underset{N'\in \rG ((\mu',0))}
{\sum_{N\in \rG ((\mu',0);\mu')}}
\rC^{H((\mu ,0)),Q_1}_{N}\,
{\rC^{H((\mu ,0)),Q_0}_{N'}}\, 
\rP^+_{\mu'}
\bigl(\overline{\zeta_{\widehat{N}}}\boxtimes \zeta_{N'})(z)
\end{align*}
by Lemmas \ref{lem:K_GLn_poly1} and \ref{lem:K_GLn_poly2}. 
Since $\{(1_{n-1},O_{n-1,1})g\mid g\in \mathrm{GL}(n,\bC )\}$
is dense in $\rM_{n-1,n}(\bC)$, 
we obtain the assertion. 
\end{proof}

\subsection{Standard Schwartz functions}
\label{subsec:standard_schwartz}

For $\lambda \in \Lambda_{n}^{\mathrm{poly}}$, we define 
two $\bC$-linear maps 
\begin{align*}
&\Phi^\circ_\lambda \colon \overline{V_{\lambda}}\boxtimes_\bC 
V_{\lambda }\ni \overline{v_1}\boxtimes v_2
\mapsto 
\rP^\circ_\lambda (\overline{v_1}\boxtimes v_2)(z)\me_{(n)}(z)
\in \cS_0(\rM_{n}(F)),\\
&\overline{\Phi^\circ_\lambda} \colon 
V_{\lambda}\boxtimes_\bC \overline{V_{\lambda }}
\ni v_1\boxtimes \overline{v_2}\mapsto \overline{\rP^\circ_\lambda (\overline{v_1}\boxtimes v_2)(z)}
\me_{(n)}(z)\in \cS_0(\rM_{n}(F))
\end{align*}
with $z\in \rM_{n}(F)$. 
By the $K_{n}\times K_{n}$-invariance of $\me_{(n)}$ and 
Lemma \ref{lem:K_GLn_poly1}, we know that 
these are $K_{n}\times K_{n}$-homomorphisms. 

When $n>1$, for $\mu \in \Lambda_{n-1}^{\mathrm{poly}}$, 
we define two $\bC$-linear maps 
\begin{align*}
&\Phi_{\mu }^+ \colon 
\overline{V_{\mu}}\boxtimes_\bC V_{(\mu ,0)}
\ni \overline{v_1}\boxtimes v_2\mapsto 
\rP^+_{\mu}(\overline{v_1}\boxtimes v_2)(z)\me_{(n-1,n)}(z)\in 
\cS_0(\rM_{n-1,n}(F)),\\
&\overline{\Phi_{\mu }^+} \colon 
V_{\mu}\boxtimes_\bC \overline{V_{(\mu ,0)}}
\ni v_1\boxtimes \overline{v_2}\mapsto 
\overline{\rP^+_{\mu}(\overline{v_1}\boxtimes v_2)(z)}\me_{(n-1,n)}(z)\in 
\cS_0(\rM_{n-1,n}(F))
\end{align*}
with $z\in \rM_{n-1,n}(F)$. By the $K_{n-1}\times K_{n}$-invariance of 
$\me_{(n-1,n)}$ and Lemma \ref{lem:K_GLn_poly2}, we know that 
these are $K_{n-1}\times K_{n}$-homomorphisms.

We regard $\cS_0 (\rM_{1,n}(F))$ and $\cS_0 (\rM_{n,1}(F))$ 
as $K_n$-modules via the actions $R$ and $L$, respectively. 
Let $l\in \bN_0$. Since $\rb (\gamma )=\rr (Q(\gamma ))^{-1}$ 
for $\gamma \in \bN_0^n$, 
the $\bC$-linear maps $\varphi_{1,n}^{(l)}$ and 
$\overline{\varphi}_{1,n}^{(l)}$ in \S \ref{subsec:def_schwartz} 
satisfy 
\begin{align*}
&\varphi_{1,n}^{(l)}(v)(z)=\rp_{1,n}^{(l)}(v)(z)\me_{(1,n)}(z),&
&\overline{\varphi}_{1,n}^{(l)}(\overline{v})(z)=
\overline{\rp_{1,n}^{(l)}(v)(z)}\me_{(1,n)}(z)
\end{align*}
for $v\in V_{(l,\mathbf{0}_{n-1})}$ and $z\in \rM_{1,n}(F)$. 
By the $K_n$-invariance of $\me_{(1,n)}$ and 
Lemma \ref{lem:poly_n1_1n}, we know that these are 
$K_n$-homomorphisms.  
We define two $\bC$-linear maps 
$\varphi_{n,1}^{(l)}\colon 
\overline{V_{(l,\mathbf{0}_{n-1})}}\to \cS_0 (\rM_{1,n}(F))$ and 
$\overline{\varphi}_{n,1}^{(l)}\colon V_{(l,\mathbf{0}_{n-1})}
\to \cS_0 (\rM_{1,n}(F))$ by 
\begin{align*}
&\varphi_{n,1}^{(l)}\colon \overline{V_{(l,\mathbf{0}_{n-1})}}
\ni \overline{v}\mapsto 
\rp_{n,1}^{(l)}(\overline{v})(z)\me_{(n,1)}(z)\in \cS_0 (\rM_{1,n}(F)),\\
&\overline{\varphi}_{n,1}^{(l)}\colon V_{(l,\mathbf{0}_{n-1})}
\ni v\mapsto 
\overline{\rp_{n,1}^{(l)}(\overline{v})(z)}\me_{(n,1)}(z)
\in \cS_0 (\rM_{1,n}(F))
\end{align*}
with $z\in \rM_{n,1}(F)$. 
By the $K_n$-invariance of $\me_{(n,1)}$ and 
Lemma \ref{lem:poly_n1_1n}, we know that these are 
$K_n$-homomorphisms.

\begin{lem}
\label{lem:fourier_explicit}
Let $\varepsilon \in \{\pm 1\}$ and $l\in \bN_0$. 
Assume $l \in \{0,1\}$ if $F=\bR$. 
Then, for $v\in V_{(l,\mathbf{0}_{n-1})}$, we have 
\begin{align*}
&\cF_\varepsilon (\varphi_{n,1}^{(l)}(\overline{v}))=
(-\varepsilon \sI )^{l}
\overline{\varphi}_{1,n}^{(l)}(\overline{v}),&
&\cF_\varepsilon (\overline{\varphi}_{n,1}^{(l)}(v))=
(-\varepsilon \sI )^{l}
\varphi_{1,n}^{(l)}(v),
\end{align*}
where the Fourier transform $\cF_\varepsilon$ is defined in 
(\ref{eqn:def_fourier}). 
\end{lem}
\begin{proof}
It suffices to show the assertion for $v=\zeta_{Q(\gamma )}$ 
with $\gamma \in \bN_0^n$ such that $\ell (\gamma )=l$. 
For $t=(t_1,t_2,\cdots ,t_n)\in \rM_{1,n}(F)$, we have 
\begin{align*}
&\cF_\varepsilon (\varphi_{n,1}^{(l)}(\overline{\zeta_{Q(\gamma )}}))(t)
=\int_{\rM_{n,1}(F)}\varphi_{n,1}^{(l)}(\overline{\zeta_{Q(\gamma )}}) (z)
\psi_{-\varepsilon }(tz)\,d_Fz\\
&\hspace{15mm}
=\sqrt{\rb (\gamma )}\prod_{i=1}^n\int_{F}z_i^{\gamma_i}
\exp \bigl(-\pi \rc_F\overline{z_i}z_i
-\pi \varepsilon \rc_F\sI (t_iz_i +\overline{t_iz_i})\bigr)\,d_Fz_i\\
&\hspace{15mm}
=\sqrt{\rb (\gamma )}
(-\varepsilon \sI \,\overline{t_i})^{\gamma_i}
\exp (-\pi \rc_F \overline{t_i}t_i)
=(-\varepsilon \sI )^l
\overline{\varphi}_{1,n}^{(l)}(\overline{\zeta_{Q(\gamma )}})(t).
\end{align*}
Here the third equality follows from the elementary formula 
\begin{align}
\label{eqn:el_fmla_int}
\int_{F}z^m
\exp \bigl(-\pi \rc_F \overline{z}z
+\pi \rc_F\sI (zt+\overline{zt})\bigr) d_{F}z
=(\sI \,\overline{t})^{m}
\exp (-\pi \rc_F\overline{t}\,t)
\end{align}
for ($t\in \bR$, $m\in \{0,1\}$) or 
($t\in \bC$, $m\in \bN_0$) according as $F=\bR$ or $F=\bC$. 
Moreover, we have 
\begin{align*}
\cF_\varepsilon (\overline{\varphi}_{n,1}^{(l)}(\zeta_{Q(\gamma )}))(t)
&=\overline{\cF_{-\varepsilon}
(\varphi_{n,1}^{(l)}(\overline{\zeta_{Q(\gamma )}}))(t)}\\
&=\overline{(\varepsilon \sI )^l
\overline{\varphi}_{1,n}^{(l)}(\overline{\zeta_{Q(\gamma )}})(t)}
=(-\varepsilon \sI )^l
\varphi_{1,n}^{(l)}(\zeta_{Q(\gamma )})(t),
\end{align*}
and completes the proof. 
\end{proof}

\section{The proofs of the Main theorems}
\label{sec:poof_main}

\subsection{Explicit calculations for the Godement sections}
\label{subsec:god_explicit}

In this subsection, we calculate the Godement sections 
in \S \ref{sec:relation}, explicitly, 
at the minimal $K_n$-types of principal series representations.

\begin{lem}
\label{lem:Pint}
Let $a=\diag (a_1,a_2,\cdots ,a_n)\in A_n$, $u\in U_n$, 
$\lambda \in \Lambda_n$ and $M\in \rG (\lambda)$. 
Then we have the following equalities 
\begin{align}
\label{eqn:Pact_inner}
&\bigl\langle \tau_{\lambda }(ua)\zeta_{M},\zeta_{M}\bigr\rangle 
=\bigl\langle \tau_{\lambda }(au)\zeta_{M},\zeta_{M}\bigr\rangle 
=\prod_{i=1}^na_i^{\gamma^M_i},\\
\label{eqn:Pact_int}
&\eta_{\rho_{n}}(a)\int_{U_{n}}
\me_{(n)}(ua)\,du
=\eta_{-\rho_{n}}(a)\int_{U_{n}}
\me_{(n)}(au)\,du
=\frac{\me_{(n)}(a)}{|\det a|_F^{(n-1)/2}},
\end{align}
where $\gamma^M=(\gamma^M_1,\gamma^M_2,\cdots ,\gamma^M_n)$ is 
the weight of $M$ defined by (\ref{eqn:def_wt_M}). 
\end{lem}
\begin{proof}
By (\ref{eqn:GT_act_wt}) and (\ref{eqn:GT_act-}), we have 
\begin{align*}
&\tau_{\lambda}(a)\zeta_{M}=\left(\prod_{i=1}^na_i^{\gamma^M_i}\right)\zeta_M,&
&\tau_{\lambda}(u)\zeta_{M}=\zeta_{M}+\sum_{N\in \rG (\lambda),\ 
\gamma^M>_{\mathrm{lex}}\gamma^N}p_{M,N}(u)\zeta_N,
\end{align*}
where $p_{M,N}$ is a polynomial function on $U_n$ 
and $>_{\mathrm{lex}}$ is the lexicographical order. 
The equality (\ref{eqn:Pact_inner}) follows from these equalities and 
the orthonormality of $\{\zeta_M\}_{M\in \rG (\lambda )}$. 
The equality (\ref{eqn:Pact_int}) follows from 
direct computation 
\begin{align*}
&\eta_{\rho_{n}}(a)\int_{U_{n}}
\me_{(n)}(ua)\,du
=\eta_{-\rho_{n}}(a)\int_{U_{n}}
\me_{(n)}(au)\,du\\
&\hspace{10mm}
=\prod_{i=1}^na_i^{-(n+1-2i)\rc_F/2}
\exp (-\pi \rc_Fa_i^2)
\prod_{j=1}^{i-1}\int_{F}
\exp (-\pi \rc_Fa_i^2\overline{u_{i,j}}u_{i,j})\,
d_Fu_{i,j}\\
&\hspace{10mm}
=\prod_{i=1}^na_i^{-(n-1)\rc_F/2}
\exp (-\pi \rc_Fa_i^2)
=\frac{\me_{(n)}(a)}{|\det a|_F^{(n-1)/2}}.
\end{align*}
Here the first equality follows from the substitution $u\to aua^{-1}$, 
and the third equality follows from the substitution 
$u_{i,j}\to a_i^{-1}u_{i,j}$ and 
the elementary formula (\ref{eqn:el_fmla_int}) with $t=m=0$. 
\end{proof}

For $\lambda =(\lambda_1,\lambda_2,\cdots, \lambda_n)\in \Lambda_n$, 
$M=(m_{i,j})_{1\leq i\leq j\leq n}\in \rG (\lambda )$ and 
$l\in \bZ$, we define $\lambda +l\in \Lambda_n$ and 
$M+l\in \rG (\lambda +l)$ by 
\begin{align*}
&\lambda +l=(\lambda_1+l,\lambda_2+l,\cdots, \lambda_n+l),& 
&M+l=(m_{i,j}+l)_{1\leq i\leq j\leq n},
\end{align*}
and denote $\lambda +(-l)$ and $M+(-l)$ simply by $\lambda -l$ and $M-l$, 
respectively. For $\lambda \in \Lambda_n$, $l\in \bZ$, 
$g\in \mathrm{GL}(n,\bC)$ and $M,N\in \rG (\lambda )$, we have  
\begin{align*}
&(\det g)^l\langle \tau_{\lambda}(g)\zeta_{M},\zeta_{N}\rangle 
=\langle \tau_{\lambda +l}(g)\zeta_{M+l},\zeta_{N+l}\rangle .
\end{align*}

\begin{lem}
\label{lem:god_explicit}
Let $d=(d_1,d_2,\cdots ,d_n)\in \bZ^n$ and 
$\nu =(\nu_1,\nu_2,\cdots ,\nu_n)\in \bC^n$. \vspace{1mm}

\noindent (1) Assume $n>1$. We take 
$\widehat{d}$ 
and $\widehat{\nu}$ as in \S \ref{subsec:god_Jacquet}. 
If $d\in \Lambda_{n,F}$, we have 
\begin{equation}
\label{eqn:god+_minKtype1}
\begin{aligned}
&\rg^+_{d_n,\nu_n}
\bigl(\,\rf_{\widehat{d},\widehat{\nu}}(\zeta_{H(\widehat{d})}),
\Phi^+_{\widehat{d}-d_n}(\overline{\zeta_{H(\widehat{d})-d_n}}
\boxtimes \zeta_{M-d_n})\,\bigr)\\
&=\frac{1}{\dim V_{\widehat{d}}}
\left(\prod_{i=1}^{n-1}\Gamma_{F}(
\nu_n-\nu_i+1;\, d_i-d_n)\right)
\rf_{d,\nu}(\zeta_M)
\end{aligned}
\end{equation}
for $M\in \rG (d)$. 
If $-d\in \Lambda_{n,F}$, we have  
\begin{equation}
\label{eqn:god+_minKtype2}
\begin{aligned}
&\rg^+_{d_n,\nu_n}
\bigl(\,\bar{\rf}_{\widehat{d},\widehat{\nu}}
(\overline{\zeta_{H(-\widehat{d})}}),
\overline{\Phi^+_{-\widehat{d}+d_n}}(\zeta_{H(-\widehat{d})+d_n}
\boxtimes \overline{\zeta_{M+d_n}})\,\bigr)\\
&=\frac{1}{\dim V_{-\widehat{d}}}
\left(\prod_{i=1}^{n-1}\Gamma_{F}(
\nu_n-\nu_i+1;\, d_n-d_i)\right)
\bar{\rf}_{d,\nu}(\overline{\zeta_M})
\end{aligned}
\end{equation}
for $M\in \rG (-d)$. 
Here $\rf_{d,\nu}$ and $\bar{\rf}_{d,\nu}$ are defined by 
(\ref{eqn:def_minKtype1}) and (\ref{eqn:def_minKtype2}), respectively. 
\\[1mm]
\noindent (2) Let $l\in \bZ$ and 
$s\in \bC$ such that $\mathrm{Re}(s)$ is sufficiently large. 
If $d\in \Lambda_{n,F}$ and $d+l\in \Lambda_n^{\mathrm{poly}}$, we have  
\begin{equation}
\label{eqn:god0_minKtype1}
\begin{aligned}
&\rg^\circ_{l,s}\bigl(\,\rf_{d,\nu}(\zeta_{H(d)}),\,
\overline{\Phi_{d+l}} (\zeta_{M+l}\boxtimes 
\overline{\zeta_{H(d)+l}})\,\bigr)\\
&=\frac{1}{\dim V_d}
\left(\prod_{i=1}^{n}
\Gamma_{F}(s+\nu_i;\,d_i+l)\right)
\rf_{d,\nu}(\zeta_{M})\hspace{10mm}(M\in \rG (d)).
\end{aligned}
\end{equation}
If $-d\in \Lambda_{n,F}$ and $-d-l\in \Lambda_n^{\mathrm{poly}}$, we have 
\begin{equation}
\label{eqn:god0_minKtype2}
\begin{aligned}
&\rg^\circ_{l,s}\bigl(\,\bar{\rf}_{d,\nu}(\overline{\zeta_{H(-d)}}),\,
\Phi^\circ_{-d-l}(\overline{\zeta_{M-l}}\boxtimes \zeta_{H(-d)-l})\,\bigr)\\
&=\frac{1}{\dim V_{-d}}\left(
\prod_{i=1}^{n}
\Gamma_{F}(s+\nu_i;\,-d_i-l)\right)
\bar{\rf}_{d,\nu}(\overline{\zeta_{M}})\hspace{10mm}(M\in \rG (-d)).
\end{aligned}
\end{equation}
\end{lem}
\begin{proof}
First, we consider the proof of the statement (1). 
Since the proofs of (\ref{eqn:god+_minKtype1}) 
and (\ref{eqn:god+_minKtype2}) are similar, 
here we will prove only (\ref{eqn:god+_minKtype1}). 
Assume $n>1$ and $d\in \Lambda_{n,F}$. We define a 
$\bC$-linear map $\rg_+ \colon V_d\to I(d,\nu )$ by 
\begin{align*}
&\rg_+ (\zeta_{M})
=\rg^+_{d_n,\nu_n}
\bigl(\,\rf_{\widehat{d},\widehat{\nu}}(\zeta_{H(\widehat{d})}),
\Phi^+_{\widehat{d}-d_n}(\overline{\zeta_{H(\widehat{d})-d_n}}
\boxtimes \zeta_{M-d_n})\,\bigr)&
&(M\in \rG (d)).
\end{align*}
Then $\rg_+$ is a $K_n$-homomorphism 
because of (\ref{eqn:god+_Kact1}). 
Since $\Hom_{K_n}(V_d, I(d,\nu ))$ is $1$ dimensional, 
there is a constant $c_+$ such that 
$\rg_+ =c_+\rf_{d,\nu}$. 
Let us calculate $c_+$. 
Since (\ref{eqn:Kn_GLn_poly_001}) and (\ref{eqn:Kn_GLn_poly_003}) imply 
\begin{align*}
&\Phi^+_{\widehat{d}-d_n}(\overline{\zeta_{H(\widehat{d})-d_n}}
\boxtimes \zeta_{H(d)-d_n})((h,O_{n-1,1}))\\
&=\bigl\langle \tau_{\widehat{d}-d_n}(h)\zeta_{H(\widehat{d})-d_n},
\zeta_{H(\widehat{d})-d_n}\bigr\rangle \me_{(n-1)}(h)
\qquad (h\in G_{n-1}),
\end{align*}
we have  
\begin{align*}
c_+=&\,c_+\rf_{d,\nu}(H(d))(1_n)=\rg_+ (\zeta_{H(d)})(1_n)\\
=&\int_{G_{n-1}}
\bigl\langle \tau_{\widehat{d}-d_n}(h)\zeta_{H(\widehat{d})-d_n},
\zeta_{H(\widehat{d})-d_n}\bigr\rangle \me_{(n-1)}(h)\\
&\times \rf_{\widehat{d},\widehat{\nu}}(\zeta_{H(\widehat{d})})(h^{-1})
\chi_{d_n}(\det h)|\det h|_F^{\nu_n+n/2}
\,dh.
\end{align*}
Decomposing $h=kua$ ($k\in K_{n-1}$, $u\in U_{n-1}$, 
$a\in A_{n-1}$) and applying Schur's orthogonality 
\cite[Corollary 1.10]{Knapp_002}
for the integration on $K_{n-1}$ 
with the equalities 
\begin{align*}
&\chi_{d_n}(\det h)
\rf_{\widehat{d},\widehat{\nu}}(\zeta_{H(\widehat{d})})(h^{-1})
=\eta_{\widehat{\nu}-\rho_{n-1}}(a^{-1})
\overline{\bigl\langle \tau_{\widehat{d}-d_n}(k)\zeta_{H(\widehat{d})-d_n},
\zeta_{H(\widehat{d})-d_n}\bigr\rangle}
\end{align*}
and $\dim V_{\widehat{d}-d_n}=\dim V_{\widehat{d}}$, we have 
\begin{align*}
c_+=&\frac{1}{\dim V_{\widehat{d}}}
\int_{A_{n-1}}\biggl(\int_{U_{n-1}}
\bigl\langle \tau_{\widehat{d}-d_n}(ua)\zeta_{H(\widehat{d})-d_n},
\zeta_{H(\widehat{d})-d_n}\bigr\rangle 
\me_{(n-1)}(ua)\,du\biggr)\\
&\times \eta_{\widehat{\nu}-\rho_{n-1}}(a^{-1})
|\det a|_F^{\nu_n+n/2}\,da.
\end{align*}
By Lemma \ref{lem:Pint} and (\ref{eqn:Gamma_F_integral}), we have 
\begin{align*}
c_+=&\frac{1}{\dim V_{\widehat{d}}}\prod_{i=1}^{n-1}
\int_{0}^\infty \exp (-\pi \rc_F a_i^2)
a_i^{(\nu_n-\nu_i+1)\rc_F+d_i-d_n}
\frac{2\rc_F\,da_i}{a_i}\\
=&\frac{1}{\dim V_{\widehat{d}}}
\prod_{i=1}^{n-1}
\Gamma_{F}\bigl(\nu_n-\nu_i+1;\, d_i-d_n\bigr).
\end{align*}
Hence, the equality (\ref{eqn:god+_minKtype1}) 
follows from $\rg_+ =c_+\rf_{d,\nu}$.

Next, we consider the proof of the statement (2). 
Since the proofs of (\ref{eqn:god0_minKtype1}) 
and (\ref{eqn:god0_minKtype2}) are similar, 
here we will prove only (\ref{eqn:god0_minKtype1}). 
Assume $d\in \Lambda_{n,F}$ and $d+l\in \Lambda_n^{\mathrm{poly}}$. 
We define a 
$\bC$-linear map $\rg_\circ \colon V_d\to I(d,\nu )$ by 
\begin{align*}
&\rg_\circ (\zeta_{M})
=\rg^\circ_{l,s}\bigl(\,\rf_{d,\nu}(\zeta_{H(d)}),\,
\overline{\Phi_{d+l}} (\zeta_{M+l}\boxtimes 
\overline{\zeta_{H(d)+l}})\,\bigr)&
&(M\in \rG (d)).
\end{align*}
Then $\rg_\circ$ is a $K_n$-homomorphism 
because of (\ref{eqn:god0_Kact1}). 
Since $\Hom_{K_n}(V_d, I(d,\nu ))$ is $1$ dimensional, 
there is a constant $c_\circ$ such that 
$\rg_\circ =c_\circ \rf_{d,\nu}$. 
Let us calculate $c_\circ $. 
By (\ref{eqn:Kn_GLn_poly_001}), we have  
\begin{align*}
c_\circ =&\,c_\circ \rf_{d,\nu}(H(d))(1_n)=\rg_\circ  (\zeta_{H(d)})(1_n)\\
=&\int_{G_{n}}
\rf_{d,\nu}(\zeta_{H(d)})(h)
\overline{\bigl\langle \tau_{d+l}(h)\zeta_{H(d)+l},
\zeta_{H(d)+l}\bigr\rangle}\,
\me_{(n)}(h)\\
&\times \chi_{l}(\det h)|\det h|_F^{s+(n-1)/2}\,dh.
\end{align*}
Decomposing $h=auk$ ($a\in A_{n}$, $u\in U_{n}$, $k\in K_{n}$) and 
applying Schur's orthogonality 
\cite[Corollary 1.10]{Knapp_002}
for the integration on $K_n$ with the equalities 
\begin{align*}
&\chi_{l}(\det h)
\rf_{d,\nu}(\zeta_{H(d)})(h)
=\eta_{\nu -\rho_{n}}(a)
\bigl\langle \tau_{d+l}(k)\zeta_{H(d)+l},
\zeta_{H(d)+l}\bigr\rangle,\\
&\tau_{d+l}(h)\zeta_{H(d)+l}
=\sum_{M\in \rG (d)}
\bigl\langle \tau_{d+l}(k)\zeta_{H(d)+l},\zeta_{M+l}\bigr\rangle \,
\tau_{d+l}(au)\zeta_{M+l}
\end{align*}
and $\dim V_{d+l}=\dim V_{d}$, we have 
\begin{align*}
c_\circ =&\frac{1}{\dim V_{d}}
\int_{A_{n}}\biggl(\int_{U_{n}}
\overline{\bigl\langle \tau_{d+l}(au)\zeta_{H(d)+l},
\zeta_{H(d)+l}\bigr\rangle }\me_{(n)}(au)\,du\biggr)\\
&\times \eta_{\nu -\rho_{n}}(a)
|\det a|_F^{s+(n-1)/2}\,da.
\end{align*}
By Lemma \ref{lem:Pint} and (\ref{eqn:Gamma_F_integral}), we have 
\begin{align*}
c_\circ =&\frac{1}{\dim V_{d}}\prod_{i=1}^{n}
\int_{0}^\infty 
\exp (-\pi \rc_F a_i^2)
a_i^{(s+\nu_i)\rc_F+d_i+l}\,\frac{2\rc_F\,da_i}{a_i}\\
=&\frac{1}{\dim V_{d}}
\prod_{i=1}^{n}\Gamma_{F}(s+\nu_i;\, d_i+l).
\end{align*}
Hence, the equality (\ref{eqn:god0_minKtype1}) follows from 
$\rg_\circ =c_\circ \rf_{d,\nu}$. 
\end{proof}

\begin{cor}
\label{cor:god_explicit_1}
We use the notation in Lemma \ref{lem:god_explicit} (1). 
If $d\in \Lambda_{n,F}$, we have  
\begin{align*}
&\rW_{\varepsilon}(\rf_{d,\nu}(\zeta_M))(g)
=\frac{\bigl(\dim V_{\widehat{d}}\,\bigr)
\chi_{d_n}(\det g)|\det g|_F^{\nu_n+(n-1)/2}}
{\prod_{i=1}^{n-1}\Gamma_{F}(\nu_n-\nu_i+1;\, d_i-d_n)}\\
&\phantom{==}\times 
\int_{G_{n-1}}\left(\int_{\rM_{n-1,1}(F)}
\Phi^+_{\widehat{d}-d_n}(\overline{\zeta_{H(\widehat{d})-d_n}}
\boxtimes \zeta_{M-d_n})\left(\left(h,hz\right)g\right)
\psi_{-\varepsilon }(e_{n-1}z)\,dz\right)\\
&\phantom{==}\times 
\rW_{\varepsilon}(\rf_{\widehat{d},\widehat{\nu}}(\zeta_{H(\widehat{d})}))
(h^{-1})\chi_{d_n}(\det h)|\det h|_F^{\nu_n+n/2}
\,dh
\end{align*}
for $M\in \rG (d)$ and $g\in G_n$. 
If $-d\in \Lambda_{n,F}$, we have  
\begin{align*}
&\rW_{\varepsilon}
(\bar{\rf}_{d,\nu}(\overline{\zeta_M}))(g)=
\frac{\bigl(\dim V_{-\widehat{d}}\,\bigr)
\chi_{d_n}(\det g)|\det g|_F^{\nu_n+(n-1)/2}}
{\prod_{i=1}^{n-1}\Gamma_{F}(\nu_n-\nu_i+1;\, d_n-d_i)}\\
&\hspace{2mm}\times 
\int_{G_{n-1}}\left(\int_{\rM_{n-1,1}(F)}
\overline{\Phi^+_{-\widehat{d}+d_n}}(\zeta_{H(-\widehat{d})+d_n}
\boxtimes \overline{\zeta_{M+d_n}})\left(\left(h,hz\right)g\right)
\psi_{-\varepsilon }(e_{n-1}z)\,dz\right)\\
&\hspace{2mm}\times 
\rW_{\varepsilon}(\bar{\rf}_{\widehat{d},\widehat{\nu}}
(\overline{\zeta_{H(-\widehat{d})}}))(h^{-1})
\chi_{d_n}(\det h)|\det h|_F^{\nu_n+n/2}
\,dh,
\end{align*}
for $M\in \rG (-d)$ and $g\in G_n$. 
\end{cor}
\begin{proof}
The assertion follows immediately from (\ref{eqn:W_god+}) and 
Lemma \ref{lem:god_explicit} (1). 
\end{proof}

\begin{cor}
\label{cor:god_explicit_2}
Let $\nu =(\nu_1,\nu_2,\cdots ,\nu_n)\in \bC^n$ and 
$d=(d_1,d_2,\cdots ,d_n)\in \bZ^n$ such that 
either $d\in \Lambda_{n,F}$ or $-d\in \Lambda_{n,F}$ holds. 
We take $\boldsymbol{\Gamma}_F(\nu;d)$ as in 
\S \ref{subsec:minKtype}. 
Then $1/{\boldsymbol{\Gamma}_F(\nu;d)}$ is nonzero 
if $\Pi_{d,\nu}$ is irreducible. 
\end{cor}
\begin{proof}
Let us give the proof for the case of $d\in \Lambda_{n,F}$. 
Let $g\in G_n$. 
By Corollary \ref{cor:god_explicit_1} and 
the entireness of the right hand side of (\ref{eqn:W_god+}), 
it is easy to show that $\rW_{\varepsilon}(\rf_{d,\nu}(\zeta_{H(d)}))(g)$
is the product of $1/\boldsymbol{\Gamma}_F(\nu;d)$ and 
some entire function of $\nu$ by induction. 
The irreducibility of $\Pi_{d,\nu}$ implies that 
there is some $g\in G_n$ such that 
$\rW_{\varepsilon}(\rf_{d,\nu}(\zeta_{H(d)}))(g)\neq 0$. 
Hence, we have 
$1/\boldsymbol{\Gamma}_F(\nu;d)\neq 0$. 
The proof for the case of $-d\in \Lambda_{n,F}$ is similar. 
\end{proof}

\subsection{Explicit recurrence relations}

Let $\varepsilon \in \{\pm 1\}$. 
Let $(\Pi_{d,\nu},I(d,\nu ))$ and 
$(\Pi_{d',\nu'},I(d',\nu'))$ be principal series representations 
of $G_n$ and $G_{n'}$, respectively, with parameters 
\begin{equation*}
\begin{aligned}
&d=(d_1,d_2,\cdots ,d_n)\in \bZ^n, \qquad &
&\nu =(\nu_1,\nu_2,\cdots ,\nu_n)\in \bC^n, \\
&d'=(d_1',d_2',\cdots ,d_{n'}')\in \bZ^{n'},& 
&\nu'=(\nu_1',\nu_2',\cdots ,\nu_{n'}')\in \bC^{n'}.
\end{aligned}
\end{equation*}

\begin{prop}
\label{prop:rec_0to+_explicit1}
Retain the notation. Assume $n'=n>1$, 
$-d'\in \Lambda_{n,F}$ and 
$d\in \Xi^\circ (-d')\cap \Lambda_{n,F}$. 
Let $l=\ell (d+d')$. 
$s\in \bC$ such that $\mathrm{Re}(s)$ is sufficiently large. 
Then we have 
\begin{align*}
&Z\bigl(s,\rW_\varepsilon (\rf_{d,\nu}(\zeta_{H(d)})),
\rW_{-\varepsilon }(\bar{\rf}_{d',\nu'}(\overline{\zeta_{H(-d')}})),
\overline{\varphi}_{1,n}^{(l)}(\overline{\zeta_{Q(d+d')}})\bigr)\\
&={\rC^{H(-d')+d_n',
Q((\mathbf{0}_{n-1},l))}_{H(-\widehat{d'})[d]+d_n'}}\,
\rC^{H(-d')+d_n',Q(d+d')}_{H(d)+d_n'}\\
&\phantom{=}\times 
\frac{\dim V_{-\widehat{d'}}}{\dim V_d}
\frac{\prod_{i=1}^{n}\Gamma_F(s+\nu_i+\nu_n';\, d_i+d_n')}
{\prod_{i=1}^{n-1}\Gamma_F(\nu_n'-\nu_i'+1;\,d_n'-d_i')}
\\
&\phantom{=}\times 
Z\bigl(s,\rW_\varepsilon (\rf_{d,\nu}(\zeta_{H(-\widehat{d'})[d]})),
\rW_{-\varepsilon }(\bar{\rf}_{\widehat{d'},\widehat{\nu'}}
(\overline{\zeta_{H(-\widehat{d'})}}))\bigr).
\end{align*}
\end{prop}
\begin{proof}
Let 
$\phi_1=\overline{\Phi^+_{-\widehat{d'}+d_n'}}(\zeta_{H(-\widehat{d'})+d_n'}
\boxtimes \overline{\zeta_{H(-d')+d_n'}})$ and 
$\phi_2=\overline{\varphi}_{1,n}^{(l)}(\overline{\zeta_{Q(d+d')}})$. 
By Proposition \ref{prop:zeta_god_recurrence1},  we have 
\begin{equation*}
\begin{aligned}
&Z\bigl(s,\,\rW_{\varepsilon} (\rf_{d,\nu}(\zeta_{H(d)})),\,
\rW_{-\varepsilon}\bigl(\rg^+_{d_n',\nu_n'}
(\bar{\rf}_{\widehat{d'},\widehat{\nu'}}
(\overline{\zeta_{H(-\widehat{d'})}}),\phi_1)\bigr),\,
\phi_2\bigr)\\
&=Z\bigl(s,\,\rW_{\varepsilon}
\bigl(\rg^\circ_{d_n',s+\nu_n'}(\rf_{d,\nu}(\zeta_{H(d)}),\phi_0 )\bigr),\,
\rW_{-\varepsilon}(\bar{\rf}_{\widehat{d'},\widehat{\nu'}}
(\overline{\zeta_{H(-\widehat{d'})}}))\bigr),
\end{aligned}
\end{equation*}
where $\phi_0(z)=\phi_1((1_{n-1},O_{n-1,1})z)
\phi_2(e_nz)$ \ ($z\in \rM_{n}(F)$). Since we have 
\begin{align*}
&\rg^+_{d_n',\nu_n'}
(\bar{\rf}_{\widehat{d'},\widehat{\nu'}}
(\overline{\zeta_{H(-\widehat{d'})}}),\phi_1)
=\frac{\prod_{i=1}^{n-1}\Gamma_{F}(
\nu_n'-\nu_i'+1;\, d_n'-d_i')}
{\dim V_{-\widehat{d'}}}\,
\bar{\rf}_{d',\nu'}(\overline{\zeta_{H(-d')}})
\end{align*}
by (\ref{eqn:god+_minKtype2}), 
it suffices to prove the equality 
\begin{equation}
\label{eqn:pf_rec_0to+_001}
\begin{aligned}
&Z\bigl(s,\,\rW_{\varepsilon}
\bigl(\rg^\circ_{d_n',s+\nu_n'}(\rf_{d,\nu}(\zeta_{H(d)}),\phi_0 )\bigr),\,
\rW_{-\varepsilon}(\bar{\rf}_{\widehat{d'},\widehat{\nu'}}
(\overline{\zeta_{H(-\widehat{d'})}}))\bigr)\\
&={\rC^{H(-d')+d_n',
Q((\mathbf{0}_{n-1},l))}_{H(-\widehat{d'})[d]+d_n'}}\,
\rC^{H(-d')+d_n',Q(d+d')}_{H(d)+d_n'}\\
&\quad \times 
\frac{\prod_{i=1}^{n}\Gamma_F(s+\nu_i+\nu_n';\, d_i+d_n')}{\dim V_d}\\
&\quad \times 
Z\bigl(s,\rW_\varepsilon (\rf_{d,\nu}(\zeta_{H(-\widehat{d'})[d]})),
\rW_{-\varepsilon }(\bar{\rf}_{\widehat{d'},\widehat{\nu'}}
(\overline{\zeta_{H(-\widehat{d'})}}))\bigr).
\end{aligned}
\end{equation}
By Lemma \ref{lem:poly_+_1n_decomp}, we have 
\begin{equation}
\label{eqn:pf_rec_0to+_003}
\begin{aligned}
&\rg^\circ_{d_n',s+\nu_n'}(\rf_{d,\nu}(\zeta_{H(d)}),\phi_0 )\\
&=
\sum_{\lambda'\in \Xi^\circ (-d'+d_n';l)}\,
\sum_{N,N'\in \rG (\lambda' )}
{\rC^{H(-d')+d_n',Q((\mathbf{0}_{n-1},l))}_{N}}\,
\rC^{H(-d')+d_n',Q(d+d')}_{N'}\\
&\phantom{=}\times 
\rg^\circ_{d_n',s+\nu_n'}(\rf_{d,\nu}(\zeta_{H(d)}),
\overline{\Phi^\circ_{\lambda'}}(\zeta_{N}\boxtimes \overline{\zeta_{N'}})).
\end{aligned}
\end{equation}
By (\ref{eqn:god0_Kact2}), we note that 
\begin{equation*}
v\otimes \overline{\zeta_M}
\mapsto \rg^\circ_{d_n',s+\nu_n'}\bigl(\rf_{d,\nu}(v),\,
\overline{\Phi^\circ_{\lambda'}}(v_1\boxtimes 
\overline{\zeta_{M+d_n'}})\bigr)(g)
\end{equation*}
defines an element of $\Hom_{K_n}(V_d\otimes_\bC \overline{V_{\lambda'-d_n'}},
\bC_{\mathrm{triv}})$ for $\lambda'\in \Xi^\circ (-d'+d_n';l)$, 
$v_1\in V_{\lambda'}$ and $g\in G_n$.  
Hence, by Lemma \ref{lem:schur_pairing}, 
for $\lambda'\in \Xi^\circ (-d'+d_n';l)$ and $N,N'\in \rG (\lambda' )$, 
we have 
\begin{align*}
&\rg^\circ_{d_n',s+\nu_n'}\bigl(\rf_{d,\nu}(\zeta_{H(d)}),\,
\overline{\Phi^\circ_{\lambda'}}(\zeta_{N}\boxtimes 
\overline{\zeta_{N'}})\bigr)(g)=0
\end{align*}
unless $\lambda'=d+d_n'$ and $N'=H(d)+d_n'$. 
By (\ref{eqn:pf_rec_0to+_003}) and this equality, we have   
\begin{equation}
\label{eqn:pf_rec_0to+_002}
\begin{aligned}
&\rg^\circ_{d_n',s+\nu_n'}(\rf_{d,\nu}(\zeta_{H(d)}),\phi_0 )\\
&=
\sum_{N\in \rG (d+d_n')}
{\rC^{H(-d')+d_n',Q((\mathbf{0}_{n-1},l))}_{N}}\,
\rC^{H(-d')+d_n',Q(d+d')}_{H(d)+d_n'}\\
&\phantom{=}\times 
\rg^\circ_{d_n',s+\nu_n'}(\rf_{d,\nu}(\zeta_{H(d)}),
\overline{\Phi^\circ_{d+d_n'}}(\zeta_{N}\boxtimes 
\overline{\zeta_{H(d)+d_n'}})).
\end{aligned}
\end{equation}
By (\ref{eqn:zeta+_Kact}) and (\ref{eqn:god0_Kact1}), 
we note that  
\begin{align*}
&\zeta_{M}\otimes \overline{v}\mapsto 
Z\bigl(s,\rW_{\varepsilon}
\bigl(\rg^\circ_{d_n',s+\nu_n'}(\rf_{d,\nu}(v_1),
\overline{\Phi^\circ_{d+d_n'}}(\zeta_{M+d_n'}\boxtimes \overline{v_2}))
\bigr),\rW_{-\varepsilon}(\bar{\rf}_{\widehat{d'},\widehat{\nu'}}
(\overline{v}))\bigr) 
\end{align*}
defines an element of 
$\Hom_{K_{n-1}}(V_d\otimes_\bC \overline{V_{-\widehat{d'}}},
\bC_{\mathrm{triv}})$ for $v_1\in V_d$ and $v_2\in V_{d+d_n'}$. 
Hence, by Lemma \ref{lem:schur_pairing+}, 
for $N\in \rG (d+d_n')$, we have 
\begin{equation*}
Z\bigl(s,\rW_{\varepsilon}
\bigl(\rg^\circ_{d_n',s+\nu_n'}(\rf_{d,\nu}(\zeta_{H(d)}),
\overline{\Phi^\circ_{d+d_n'}}(\zeta_{N}\boxtimes \overline{\zeta_{H(d)+d_n'}})
)\bigr),\rW_{-\varepsilon}(\bar{\rf}_{\widehat{d'},\widehat{\nu'}}
(\overline{\zeta_{H(-\widehat{d'})}}))\bigr)=0
\end{equation*}
unless $\widehat{N}=H(-\widehat{d'})+d_n'$. 
By (\ref{eqn:god0_minKtype1}), (\ref{eqn:pf_rec_0to+_002}) 
and this equality, 
we obtain (\ref{eqn:pf_rec_0to+_001}). 
\end{proof}

\begin{prop}
\label{prop:rec_+to0_explicit1}
Retain the notation. Assume $n'=n-1$, $d\in \Lambda_{n,F}$ and 
$-d'\in \Xi^+(d)\cap \Lambda_{n-1,F}$. 
Let $l=\ell (\widehat{d}+d')$. Then we have 
\begin{align*}
&Z\bigl(s,\rW_\varepsilon (\rf_{d,\nu}(\zeta_{H(-d')[d]})),
\rW_{-\varepsilon}(\bar{\rf}_{d',\nu'}
(\overline{\zeta_{H(-d')}}))\bigr)\\
&=(-\varepsilon \sI )^{l}\,
\Bigl(\rC^{H((-d'-d_n,0)),\,Q((\widehat{d}+d',0))}_{H(d)-d_n}\,
{\rC^{H((-d'-d_n,0)),\,Q((\mathbf{0}_{n-1},l))}_{H(-d')[d]-d_n}}\,
\Bigr)^{-1}\\
&\hspace{5mm}
\times \frac{\dim V_{\widehat{d}}}{\dim V_{-d'}}
\frac{\prod_{i=1}^{n-1}\Gamma_F(
s+\nu_n+\nu_i';\,-d_n-d_i')}
{\prod_{i=1}^{n-1}\Gamma_F(
\nu_n-\nu_i+1;\, -d_n+d_i)}\\
&\hspace{5mm}
\times Z\bigl(s,\rW_\varepsilon (\rf_{\widehat{d},\widehat{\nu}}
(\zeta_{H(\widehat{d})})),
\rW_{-\varepsilon }(\bar{\rf}_{d',\nu'}
(\overline{\zeta_{H(-d')}})),
\overline{\varphi}_{1,n-1}^{(l)}(\overline{\zeta_{Q(\widehat{d}+d')}})\bigr).
\end{align*}
\end{prop}
\begin{proof}
Let 
$\phi_1=\Phi^\circ_{-d'-d_n}\bigl(\overline{\zeta_{H(-d')-d_n}}
\boxtimes \zeta_{H(-d')-d_n}\bigr)$ and 
$\phi_2=\varphi_{n-1,1}^{(l)}(\overline{\zeta_{Q(\widehat{d}+d')}})$. 
By Proposition \ref{prop:zeta_god_recurrence2},  we have 
\begin{align*}
&Z\bigl(s,\,
\rW_{\varepsilon}\bigl(\rg^+_{d_n,\nu_n}(
\rf_{\widehat{d},\widehat{\nu}}(\zeta_{H(\widehat{d})})
,\phi_0)\bigr),\,
\rW_{-\varepsilon}(
\bar{\rf}_{d',\nu'}(\overline{\zeta_{H(-d')}})
)\bigr)\\
&=Z\bigl(s,\,
\rW_{\varepsilon}(\rf_{\widehat{d},\widehat{\nu}}(\zeta_{H(\widehat{d})})),\,
\rW_{-\varepsilon}(\rg^\circ_{d_n,s+\nu_n}
(\bar{\rf}_{d',\nu'}(\overline{\zeta_{H(-d')}}),\phi_1)),\,
\cF_\varepsilon (\phi_2)\bigr),
\end{align*}
where $\phi_0(z)=\phi_1(z\,{}^t\!(1_{n-1},O_{n-1,1}))\phi_2(z\,{}^t\!e_{n})$ \ 
($z\in \rM_{n-1,n}(F)$). 
We have 
\begin{equation*}
\rg^\circ_{d_n,s+\nu_n}(\bar{\rf}_{d',\nu'}
(\overline{\zeta_{H(-d')}}),\,\phi_1)
=\frac{\prod_{i=1}^{n}
\Gamma_{F}(s+\nu_n+\nu_i';\,-d_n-d_i')}{\dim V_{-d'}}
\bar{\rf}_{d',\nu'}(\overline{\zeta_{H(-d')}})
\end{equation*}
by (\ref{eqn:god0_minKtype2}). 
Because of these equalities and Lemma \ref{lem:fourier_explicit}, 
it suffices to prove 
\begin{equation}
\label{eqn:pf_rec_+to0_001}
\begin{aligned}
&Z\bigl(s,\,
\rW_{\varepsilon}\bigl(\rg^+_{d_n,\nu_n}(
\rf_{\widehat{d},\widehat{\nu}}(\zeta_{H(\widehat{d})})
,\phi_0)\bigr),\,
\rW_{-\varepsilon}(
\bar{\rf}_{d',\nu'}(\overline{\zeta_{H(-d')}})
)\bigr)\\
&=
\rC^{H((-d'-d_n,0)),\,Q((\widehat{d}+d',0))}_{H(d)-d_n}\,
{\rC^{H((-d'-d_n,0)),\,Q((\mathbf{0}_{n-1},l))}_{H(-d')[d]-d_n}}\\
&\quad \times 
\frac{\prod_{i=1}^{n-1}\Gamma_F(
\nu_n-\nu_i+1;\, -d_n+d_i)}{\dim V_{\widehat{d}}}\\
&\quad 
\times Z\bigl(s,\rW_\varepsilon (\rf_{d,\nu}(\zeta_{H(-d')[d]})),
\rW_{-\varepsilon}(\bar{\rf}_{d',\nu'}
(\overline{\zeta_{H(-d')}}))\bigr).
\end{aligned}
\end{equation}
By Lemma \ref{lem:poly_0_n1_decomp}, 
we have 
\begin{equation}
\label{eqn:pf_rec_+to0_003}
\begin{aligned}
&\rg^+_{d_n,\nu_n}(
\rf_{\widehat{d},\widehat{\nu}}(\zeta_{H(\widehat{d})})
,\phi_0)\\
&=
\sum_{\lambda' \in \Xi^\circ (-d'-d_n;l)}
\,\underset{N'\in \rG ((\lambda',0))}
{\sum_{N\in \rG ((\lambda',0);\lambda')}}
\rC^{H((-d'-d_n,0)),Q((\widehat{d}+d',0))}_{N}\\
&\quad 
\times {\rC^{H((-d'-d_n,0)),Q((\mathbf{0}_{n-1},l))}_{N'}}\,
\rg^+_{d_n,\nu_n}\bigl(
\rf_{\widehat{d},\widehat{\nu}}(\zeta_{H(\widehat{d})})
,\Phi^+_{\lambda'}\bigl(\overline{\zeta_{\widehat{N}}}\boxtimes \zeta_{N'})
\bigr).
\end{aligned}
\end{equation}
By (\ref{eqn:god+_Kact2}), we note that 
\begin{align*}
&v\otimes \overline{\zeta_M}\mapsto 
\rg^+_{d_n,\nu_n}\bigl(\rf_{\widehat{d},\widehat{\nu}}(v),
\Phi^+_{\lambda'}\bigl(\overline{\zeta_{M-d_n}}\boxtimes v_1)
\bigr)(g)
\end{align*}
defines an element of $\Hom_{K_{n-1}}
(V_{\widehat{d}}\otimes_\bC \overline{V_{\lambda'+d_n}},\bC_{\mathrm{triv}})$ 
for $\lambda'\in \Xi^\circ (-d'-d_n;l)$, 
$v_1\in V_{(\lambda',0)}$ and $g\in G_n$.  
Hence, by Lemma \ref{lem:schur_pairing}, 
for $N\in \rG ((\lambda',0);\lambda')$, $N'\in \rG ((\lambda',0))$ 
and $\lambda'\in \Xi^\circ (-d'-d_n;l)$, we have 
\begin{align*}
&\rg^+_{d_n,\nu_n}\bigl(
\rf_{\widehat{d},\widehat{\nu}}(\zeta_{H(\widehat{d})})
,\Phi^+_{\lambda'}\bigl(\overline{\zeta_{\widehat{N}}}\boxtimes \zeta_{N'})
\bigr)=0
\end{align*}
unless $\lambda'=\widehat{d}-d_n$ and $\widehat{N}=H(\widehat{d})-d_n$. 
By (\ref{eqn:pf_rec_+to0_003}) and this equality, we have   
\begin{equation}
\label{eqn:pf_rec_+to0_002}
\begin{aligned}
&\rg^+_{d_n,\nu_n}(
\rf_{\widehat{d},\widehat{\nu}}(\zeta_{H(\widehat{d})})
,\phi_0)\\
&=
\sum_{N'\in \rG (d-d_n)}
\rC^{H((-d'-d_n,0)),Q((\widehat{d}+d',0))}_{H(d)-d_n}
\,{\rC^{H((-d'-d_n,0)),Q((\mathbf{0}_{n-1},l))}_{N'}}\\
&\quad \times 
\rg^+_{d_n,\nu_n}\bigl(
\rf_{\widehat{d},\widehat{\nu}}(\zeta_{H(\widehat{d})})
,\Phi^+_{\widehat{d}-d_n}\bigl(\overline{\zeta_{H(\widehat{d})-d_n}}
\boxtimes \zeta_{N'})\bigr).
\end{aligned}
\end{equation}
By (\ref{eqn:zeta+_Kact}) and (\ref{eqn:god+_Kact1}), 
we note that  
\begin{equation*}
\zeta_{M}\otimes \overline{v}\mapsto 
Z\bigl(s,\,
\rW_{\varepsilon}\bigl(\rg^+_{d_n,\nu_n}(
\rf_{\widehat{d},\widehat{\nu}}(v_1)
,\Phi^+_{\widehat{d}-d_n}\bigl(\overline{v_2}
\boxtimes \zeta_{M-d_n}))\bigr),\,
\rW_{-\varepsilon}(
\bar{\rf}_{d',\nu'}(\overline{v})
)\bigr)
\end{equation*}
defines an element of 
$\Hom_{K_{n-1}}(V_d\otimes_\bC \overline{V_{-d'}},
\bC_{\mathrm{triv}})$ for $v_1\in V_{\widehat{d}}$ 
and $v_2\in V_{\widehat{d}-d_n}$. 
Hence, by Lemma \ref{lem:schur_pairing+}, 
for $N'\in \rG (d-d_n)$, we have 
\begin{equation*}
Z\bigl(s,\,
\rW_{\varepsilon}\bigl(\rg^+_{d_n,\nu_n}(
\rf_{\widehat{d},\widehat{\nu}}(\zeta_{H(\widehat{d})})
,\Phi^+_{\widehat{d}-d_n}\bigl(\overline{\zeta_{H(\widehat{d})-d_n}}
\boxtimes \zeta_{N'}))\bigr),\,
\rW_{-\varepsilon}(
\bar{\rf}_{d',\nu'}(\overline{\zeta_{H(-d')}})
)\bigr)=0
\end{equation*}
unless $\widehat{N'}=H(-\widehat{d'})-d_n$. 
By (\ref{eqn:god+_minKtype1}), (\ref{eqn:pf_rec_+to0_002}) 
and this equality, 
we obtain (\ref{eqn:pf_rec_+to0_001}). 
\end{proof}

\begin{thm}
\label{thm:sub1}
Retain the notation. 
Assume $d\in \Lambda_{n,F}$ and $-d'\in \Lambda_{n',F}$. 
We take $\boldsymbol{\Gamma}_F(\nu;d)$ and 
$\boldsymbol{\Gamma}_F(\nu';d')$ as in 
\S \ref{subsec:minKtype}. 
\vspace{1mm}

\noindent 
(1) Assume $n'=n$ and $d\in \Xi^\circ (-d')$. 
Let $l=\ell (d+d')$.  
Then we have 
\begin{align*}
&Z\bigl(s,\rW_\varepsilon (\rf_{d,\nu}(\zeta_{H(d)})),
\rW_{-\varepsilon }(\bar{\rf}_{d',\nu'}(\overline{\zeta_{H(-d')}})),
\overline{\varphi}_{1,n}^{(l)}(\overline{\zeta_{Q(d+d')}})
\bigr)\\
&=
\frac{(-\varepsilon \sI )^{\sum_{i=1}^{n-1}(n-i)(d_i+d_i')}
\sqrt{\rb (d+d')}\rC^\circ (d;-d')
L(s,\Pi_{d,\nu}\times \Pi_{d',\nu'})}
{(\dim V_{d})\boldsymbol{\Gamma}_F(\nu;d)
\boldsymbol{\Gamma}_F(\nu';d')}.
\end{align*}

\noindent 
(2) Assume $n'=n-1$ and $-d'\in \Xi^+(d)$. 
Then we have 
\begin{align*}
&Z\bigl(s,\rW_\varepsilon (\rf_{d,\nu}(\zeta_{H(-d')[d]})),
\rW_{-\varepsilon }(\bar{\rf}_{d',\nu'}(\overline{\zeta_{H(-d')}}))
\bigr)\\
&=
\frac{(-\varepsilon \sI )^{\sum_{i=1}^{n-1}(n-i)(d_i+d_i')}
L(s,\Pi_{d,\nu}\times \Pi_{d',\nu'})}
{(\dim V_{-d'})\sqrt{\rr (H(-d')[d])}
\boldsymbol{\Gamma}_F(\nu;d)
\boldsymbol{\Gamma}_F(\nu';d')}.
\end{align*}
Here $\rr (H(-d')[d])$ is defined by (\ref{eqn:def_rM}). 
\end{thm}
\begin{proof}
Let us prove the statement (1) by induction with respect to $n$. 
First, we consider the case of $n=1$. 
Since 
\begin{align*}
&\rW_{\varepsilon}\bigl(\rf_{d_1,\nu_1}(\zeta_{d_1})\bigr)(ak)
=a^{\nu_1\rc_F}k^{d_1},\hspace{5mm}
\rW_{-\varepsilon}\bigl(\bar{\rf}_{d_1',\nu_1'}
(\overline{\zeta_{-d_1'}})\bigr)(ak)=a^{\nu_1'\rc_F}k^{d_1'},\\
&\overline{\varphi}_{1,1}^{(d_1+d_1')}(\overline{\zeta_{d_1+d_1'}})(ak)
=(a\overline{k})^{d_1+d_1'}\exp (-\pi \rc_F a^2)
\end{align*}
for $a\in A_1=\bR_{+}^\times $ and $k\in K_1$, we have 
\begin{align*}
&Z\bigl(s,
\rW_{\varepsilon}\bigl(\rf_{d_1,\nu_1}(\zeta_{d_1})\bigr), 
\rW_{-\varepsilon}\bigl(\bar{\rf}_{d_1',\nu_1'}
(\overline{\zeta_{-d_1'}})\bigr),
\overline{\varphi}_{1,1}^{(d_1+d_1')}(\overline{\zeta_{d_1+d_1'}})\bigr)\\
&=\left(\int_0^\infty \exp (-\pi \rc_F a^2)a^{(s+\nu_1+\nu_1')\rc_F+d_1+d_1'}
\,\frac{2\rc_F\,da}{a}\right)\left(\int_{K_1}dk\right)\\
&=\Gamma_F (s+\nu_1+\nu_1';\, d_1+d_1')
=L(s,\Pi_{d_1,\nu_1}\times \Pi_{d_1',\nu_1'}).
\end{align*}
Here the second equality follows from (\ref{eqn:Gamma_F_integral}). 
Next, we consider the case of $n\geq 2$. 
Let $q=\ell (\widehat{d}+\widehat{d'})$. 
By Propositions \ref{prop:rec_0to+_explicit1} and 
\ref{prop:rec_+to0_explicit1}, we have 
\begin{align*}
&Z\bigl(s,\rW_\varepsilon (\rf_{d,\nu}(\zeta_{H(d)})),
\rW_{-\varepsilon }(\bar{\rf}_{d',\nu'}(\overline{\zeta_{H(-d')}})),
\overline{\varphi}_{1,n}^{(l)}(\overline{\zeta_{Q(d+d')}})\bigr)\\
&=(-\varepsilon \sI )^{q}\,
\frac{\rC^{H(-d')+d_n',Q(d+d')}_{H(d)+d_n'}}
{\rC^{H((-\widehat{d'}-d_n,0)),\,
Q((\widehat{d}+\widehat{d'},0))}_{H(d)-d_n}}
\,
\frac{{\rC^{H(-d')+d_n',
Q((\mathbf{0}_{n-1},l))}_{H(-\widehat{d'})[d]+d_n'}}}
{{\rC^{H((-\widehat{d'}-d_n,0)),\,
Q((\mathbf{0}_{n-1},q))}_{H(-\widehat{d'})[d]-d_n}}}\\
&\quad \times 
\frac{\dim V_{\widehat{d}}}{\dim V_d}
\frac{\prod_{i=1}^{n}\Gamma_F(s+\nu_i+\nu_n';\, d_i+d_n')}
{\prod_{i=1}^{n-1}\Gamma_F(\nu_n'-\nu_i'+1;\,d_n'-d_i')}
\frac{\prod_{i=1}^{n-1}\Gamma_F(
s+\nu_n+\nu_i';\,-d_n-d_i')}
{\prod_{i=1}^{n-1}\Gamma_F(
\nu_n-\nu_i+1;\, -d_n+d_i)}\\
&\quad 
\times Z\bigl(s,\rW_\varepsilon (\rf_{\widehat{d},\widehat{\nu}})
(\zeta_{H(\widehat{d})})),
\rW_{-\varepsilon }(\bar{\rf}_{\widehat{d'},\widehat{\nu'}}
(\overline{\zeta_{H(-\widehat{d'})}})),
\overline{\varphi}_{1,n-1}^{(q)}
(\overline{\zeta_{Q(\widehat{d}+\widehat{d'})}})\bigr).
\end{align*}
Moreover, by (\ref{eqn:CG_explicit}), we have 
\begin{align*}
&\frac{\rC^{H(-d')+d_n',Q(d+d')}_{H(d)+d_n'}}
{\rC^{H((-\widehat{d'}-d_n,0)),\,
Q((\widehat{d}+\widehat{d'},0))}_{H(d)-d_n}}
\,
\frac{{\rC^{H(-d')+d_n',
Q((\mathbf{0}_{n-1},l))}_{H(-\widehat{d'})[d]+d_n'}}}
{{\rC^{H((-\widehat{d'}-d_n,0)),\,
Q((\mathbf{0}_{n-1},q))}_{H(-\widehat{d'})[d]-d_n}}}\\
&=
\sqrt{\frac{l!}{q!(d_n+d_n')!}}
\prod_{h=1}^{n-1}
\frac{(d_{h}-d_{n}-h+n)!(-d_{h}'+d_{n}'-h+n-1)!}
{(d_{h}+d_{n}'-h+n)!(-d_{h}'-d_{n}-h+n-1)!}.
\end{align*}
By the above equalities and the induction hypothesis, 
we obtain the formula in the statement (1). 

The statement (2) follows 
from (\ref{eqn:CG_explicit}), Proposition \ref{prop:rec_+to0_explicit1}, 
the statement (1) and 
\begin{align*}
&\frac{1}{\sqrt{\rr (H(\mu )[\lambda ])}} 
=\sqrt{\prod_{1\leq i\leq  j\leq n-1}
\frac{(\mu_i-\mu_j-i+j)!(\lambda_{i}-\lambda_{j+1}-i+j)!}
{(\lambda_i-\mu_j-i+j)!(\mu_{i}-\lambda_{j+1}-i+j)!}}
\end{align*}
for $\lambda =(\lambda_1,\lambda_2,\cdots ,\lambda_n)\in \Lambda_n$ 
and $\mu =(\mu_1,\mu_2,\cdots ,\mu_{n-1})\in \Xi^+(\lambda )$. 
\end{proof}

\begin{proof}[Proof of Theorem \ref{thm:main1}]
The equality (\ref{eqn:thm_main1}) follows from 
Theorem \ref{thm:sub1} (2). 
Since (\ref{eqn:zeta_hom_1}) is an element of 
$\Hom_{K_{n-1}}(V_{d}\otimes_\bC \overline{V_{-d'}},\bC_{\mathrm{triv}})$, 
we completes the proof by Lemma \ref{lem:schur_pairing+}. 
\end{proof}

\begin{proof}[Proof of Theorem \ref{thm:main2} (2)]
The equality (\ref{eqn:main2_12}) follows from 
Theorem \ref{thm:sub1} (1). Since 
(\ref{eqn:zeta_hom_22}) is an element of 
$\Hom_{K_{n}}(V_{d}\otimes_\bC \overline{V_{-d'}}
\otimes_\bC \overline{V_{(-l,\mathbf{0}_{n-1})}},\bC_{\mathrm{triv}})$, 
we completes the proof by Lemma \ref{lem:schur_pairing0}. 
\end{proof}

Similar to Propositions \ref{prop:rec_0to+_explicit1} 
and \ref{prop:rec_+to0_explicit1}, 
we obtain the following propositions. 

\begin{prop}
\label{prop:rec_0to+_explicit2}
Retain the notation. Assume $n'=n>1$, $d'\in \Lambda_{n,F}$ and 
$-d\in \Xi^\circ (d')\cap \Lambda_{n,F}$. 
Let $l=\ell (-d-d')$. Then we have 
\begin{align*}
&Z\bigl(s,\rW_\varepsilon (\bar{\rf}_{d,\nu}(\overline{\zeta_{H(-d)}})),
\rW_{-\varepsilon }(\rf_{d',\nu'}(\zeta_{H(d')})),
\varphi_{1,n}^{(l)}(\zeta_{Q(-d-d')})\bigr)\\
&=\rC^{H(d')-d_n',Q((\mathbf{0}_{n-1},l))}_{H(\widehat{d'})[-d]-d_n'}\,
{\rC^{H(d')-d_n',Q(d+d')}_{H(-d)-d_n'}}\\
&\quad 
\times \frac{\dim V_{\widehat{d'}}}{\dim V_{-d}}
\frac{\prod_{i=1}^{n}\Gamma_F(s+\nu_i+\nu_n';\, -d_i-d_n')}
{\prod_{i=1}^{n-1}\Gamma_F(\nu_n'-\nu_i'+1;\,-d_n'+d_i')}\\
&\quad \times 
Z\bigl(s,\rW_\varepsilon (\bar{\rf}_{d,\nu}
(\overline{\zeta_{H(\widehat{d'})[-d]}})),
\rW_{-\varepsilon }(\rf_{\widehat{d'},\widehat{\nu'}}
(\zeta_{H(\widehat{d'})}))\bigr).
\end{align*}
\end{prop}

\begin{prop}
\label{prop:rec_+to0_explicit2}
Retain the notation. Assume $n'=n-1$, $-d\in \Lambda_{n,F}$ and 
$d'\in \Xi^+(-d)\cap \Lambda_{n-1,F}$. 
Let $l=\ell (-\widehat{d}-d')$. Then we have 
\begin{align*}
&Z\bigl(s,\rW_\varepsilon (\bar{\rf}_{d,\nu}
(\overline{\zeta_{H(d')[-d]}})),
\rW_{-\varepsilon}(\rf_{d',\nu'}(\zeta_{H(d')}))\bigr)\\
&=(-\varepsilon \sI )^{l}\,
\Bigl({\rC^{H((d'+d_n,0)),\,Q((-\widehat{d}-d',0))}_{H(-d)+d_n}}\,
\rC^{H((d'+d_n,0)),\,Q((\mathbf{0}_{n-1},l))}_{H(d')[-d]+d_n}\,
\Bigr)^{-1}\\
&\quad \times \frac{\dim V_{-\widehat{d}}}{\dim V_{d'}}
\frac{\prod_{i=1}^{n-1}\Gamma_F(
s+\nu_n+\nu_i';\,d_n+d_i')}
{\prod_{i=1}^{n-1}\Gamma_F(
\nu_n-\nu_i+1;\, d_n-d_i)}\\
&\quad \times Z\bigl(s,\rW_\varepsilon (\bar{\rf}_{\widehat{d},\widehat{\nu}}
(\overline{\zeta_{H(-\widehat{d})}})),
\rW_{-\varepsilon }(\rf_{d',\nu'}(\zeta_{H(d')})),
\varphi_{1,n-1}^{(l)}(\zeta_{Q(-\widehat{d}-d')})\bigr).
\end{align*}
\end{prop}

Similar to Theorem \ref{thm:sub1}, 
we obtain the following theorem using Propositions 
\ref{prop:rec_0to+_explicit2} and \ref{prop:rec_+to0_explicit2}.

\begin{thm}
\label{thm:sub2}
Retain the notation. 
Assume $-d\in \Lambda_{n,F}$ and $d'\in \Lambda_{n',F}$. 
We take $\boldsymbol{\Gamma}_F(\nu;d)$ and 
$\boldsymbol{\Gamma}_F(\nu';d')$ as in 
\S \ref{subsec:minKtype}. 
\vspace{1mm}

\noindent 
(1) Assume $n'=n$ and $-d\in \Xi^\circ (d')$. 
Let $l=\ell (-d-d')$. 
Then we have 
\begin{align*}
&Z\bigl(s,\rW_\varepsilon 
(\bar{\rf}_{d,\nu}(\overline{\zeta_{H(-d)}})),
\rW_{-\varepsilon }(\rf_{d',\nu'}(\zeta_{H(d')})),
\varphi_{1,n}^{(l)}(\zeta_{Q(-d-d')})\bigr)\\
&=
\frac{(\varepsilon \sI )^{\sum_{i=1}^{n-1}(n-i)(d_i+d_i')}
\sqrt{\rb (-d-d')}\rC^\circ (-d;d')
L(s,\Pi_{d,\nu}\times \Pi_{d',\nu'})}
{(\dim V_{-d})\boldsymbol{\Gamma}_F(\nu;d)
\boldsymbol{\Gamma}_F(\nu';d')}.
\end{align*}

\noindent 
(2) Assume $n'=n-1$ and $d'\in \Xi^+(-d)$. 
Then we have 
\begin{align*}
&Z\bigl(s,\rW_\varepsilon (\bar{\rf}_{d,\nu}(\overline{\zeta_{H(d')[-d]}})),
\rW_{-\varepsilon }(\rf_{d',\nu'}(\zeta_{H(d')}))
\bigr)\\
&=
\frac{(\varepsilon \sI )^{\sum_{i=1}^{n-1}(n-i)(d_i+d_i')}
L(s,\Pi_{d,\nu}\times \Pi_{d',\nu'})}
{(\dim V_{d'})\sqrt{\rr (H(d')[-d])}
\boldsymbol{\Gamma}_F(\nu;d)
\boldsymbol{\Gamma}_F(\nu';d')}.
\end{align*}
\end{thm}

\begin{proof}[Proof of Theorem \ref{thm:main2} (1)]
The equality (\ref{eqn:main2_22}) follows from 
Theorem \ref{thm:sub2} (1) and (\ref{eqn:zeta0_symmetry}). Since 
(\ref{eqn:zeta_hom_21}) is an element of 
$\Hom_{K_{n}}(V_{d}\otimes_\bC \overline{V_{-d'}}
\otimes_\bC V_{(l,\mathbf{0}_{n-1})},\bC_{\mathrm{triv}})$, 
we completes the proof by Lemma \ref{lem:schur_pairing0} and 
the properties of complex conjugate representations 
in \S \ref{subsec:com_conj_rep}. 
\end{proof}

\appendix 

\section{Explicit formulas of Whittaker functions} 

In this appendix, we consider the explicit formulas 
of the radial parts of Whittaker functions on $G_n$. 
Let $\varepsilon \in \{\pm 1\}$, 
$d=(d_1,d_2,\cdots,d_n)\in \bZ^n$, and 
$\nu =(\nu_1,\nu_2,\cdots,\nu_n)\in \bC^n$. Assume that 
either $d\in \Lambda_{n,F}$ or $-d\in \Lambda_{n,F}$ holds. 
We set 
\begin{align*}
&\widetilde{W}_{d,\nu}^{(\varepsilon )}(a)
=\left\{\begin{array}{ll}
\eta_{-\rho_n}(a)\rW_{\varepsilon}(\rf_{d,\nu}(\xi_{H(d)}))(a)
&\text{if $d\in \Lambda_{n,F}$},\\[1mm]
\eta_{-\rho_n}(a)\rW_{\varepsilon}
(\bar{\rf}_{d,\nu}(\overline{\xi_{H(-d)}}))(a)
&\text{if $-d\in \Lambda_{n,F}$}
\end{array}\right.&
&(a\in A_n). 
\end{align*}
Then we have the following theorem, which is 
the generalization of 
the explicit formulas \cite[Theorem 14]{Ishii_Stade_001} 
of spherical Whittaker functions on $\mathrm{GL}(n,\bR)$.

\begin{thm}
Retain the notation, and we assume $n>1$. 
We take $\widehat{d}$ and $\widehat{\nu}$ as in \S \ref{subsec:god_Jacquet}. 
Let $a=\diag (a_1,\,a_2,\,\cdots ,a_n )\in A_n$. Then we have 
\begin{align*}
&\widetilde{W}_{d,\nu}^{(\varepsilon )}(a)
=\frac{\prod_{i=1}^{n}a_i^{\nu_n\rc_F+|d_i-d_n|}}
{\prod_{i=1}^{n-1}\Gamma_{F}(\nu_n-\nu_i+1;\, |d_i-d_n|)}
\\
&\hspace{5mm}\times 
\int_{(\bR_{+}^\times )^{n-1}}
\widetilde{W}^{(\varepsilon )}_{\widehat{d},\widehat{\nu}}(t)
\prod_{i=1}^{n-1}
\exp \biggl(-\pi \rc_F
\biggl(\frac{t_i^2}{a_{i+1}^{2}}+\frac{a_i^2}{t_i^2}\biggr)\biggr)
t_i^{-\nu_n\rc_F-|d_i-d_n|}
\,\frac{2\rc_F\,dt_i}{t_i}
\end{align*}
with $t=\diag (t_1,\,t_2,\,\cdots ,t_{n-1} )\in A_{n-1}$. 
\end{thm}
\begin{proof}
We will prove here only the case of $d\in \Lambda_{n,F}$, 
since the proof for the case of $-d\in \Lambda_{n,F}$ is similar. 
By Lemmas \ref{lem:K_GLn_poly1} and \ref{lem:K_GLn_poly2}, 
we have 
\begin{align*}
&\Phi^+_{\widehat{d}-d_n}
(\overline{\zeta_{H(\widehat{d})-d_n}}\boxtimes \zeta_{H(d)-d_n})((h,hz)a)\\
&=\left(\prod_{i=1}^{n-1}a_i^{d_i-d_n}\right)
\bigl\langle \tau_{\widehat{d}-d_n}(h)\zeta_{H(\widehat{d})-d_n},\,
\zeta_{H(\widehat{d})-d_n}\bigr\rangle 
\me_{(n-1,n)}((h,hz)a)
\end{align*}
for $h\in G_{n-1}$ and $z\in \mathrm{M}_{n-1,1}(F)$. 
Hence, by Corollary \ref{cor:god_explicit_1}, we have  
\begin{align*}
&\widetilde{W}_{d,\nu}^{(\varepsilon )}(a)=
\eta_{-\rho_n}(a)\rW_{\varepsilon}(\rf_{d,\nu}(\zeta_{H(d)}))(a)\\
&=\frac{\bigl(\dim V_{\widehat{d}}\,\bigr)
\prod_{i=1}^{n}a_i^{(\nu_n+i-1)\rc_F+d_i-d_n}}
{\prod_{i=1}^{n-1}\Gamma_{F}(\nu_n-\nu_i+1;\, d_i-d_n)}\\
&\phantom{=}
\times 
\int_{G_{n-1}}\left(\int_{\rM_{n-1,1}(F)}
\me_{(n-1,n)}((h,hz)a)
\psi_{-\varepsilon }(e_{n-1}z)\,dz\right)
\rW_{\varepsilon}(\rf_{\widehat{d},\widehat{\nu}}(\zeta_{H(\widehat{d})}))
(h^{-1})\\
&\phantom{=}
\times 
\bigl\langle \tau_{\widehat{d}-d_n}(h)\zeta_{H(\widehat{d})-d_n},\,
\zeta_{H(\widehat{d})-d_n}\bigr\rangle 
\chi_{d_n}(\det h)|\det h|_F^{\nu_n+n/2}
\,dh.
\end{align*}
Decomposing $h^{-1}=xtk$ 
($x\in N_{n-1}$, $t=\diag (t_1,t_2,\cdots ,t_{n-1})\in A_{n-1}$, 
$k\in K_{n-1}$) and applying Schur's orthogonality 
\cite[Corollary 1.10]{Knapp_002} 
for the integration on $K_{n-1}$ 
together with the equalities 
\begin{align*}
&\bigl\langle \tau_{\widehat{d}-d_n}(h)\zeta_{H(\widehat{d})-d_n},\,
\zeta_{H(\widehat{d})-d_n}\bigr\rangle 
\chi_{d_n}(\det h)\\
&=\left(\prod_{i=1}^{n-1}t_i^{-d_i+d_n}\right)
\overline{\bigl\langle \tau_{\widehat{d}}(k)
\zeta_{H(\widehat{d})},\,\zeta_{H(\widehat{d})}
\bigr\rangle}&
&(\text{by Lemma \ref{lem:Pint}})
\end{align*}
and 
\begin{align*}
\rW_{\varepsilon}(\rf_{\widehat{d},\widehat{\nu}}(\zeta_{H(\widehat{d})}))
(h^{-1})
&=\psi_{\varepsilon ,n-1}(x)
\rW_{\varepsilon}(\rf_{\widehat{d},\widehat{\nu}}
(\tau_{\widehat{d}}(k)\zeta_{H(\widehat{d})}))(t)\\
&=\sum_{M\in \rG (\widehat{d})}
\langle \tau_{\widehat{d}}(k)\zeta_{H(\widehat{d})},\zeta_{M}\rangle \ 
\psi_{\varepsilon ,n-1}(x)
\rW_{\varepsilon}(\rf_{\widehat{d},\widehat{\nu}}(\zeta_{M}))(t),
\end{align*}
we have  
\begin{align*}
\widetilde{W}_{d,\nu}^{(\varepsilon )}(a)
=&\frac{\prod_{i=1}^{n}a_i^{(\nu_n+i-1)\rc_F+d_i-d_n}}
{\prod_{i=1}^{n-1}\Gamma_{F}(\nu_n-\nu_i+1;\, d_i-d_n)}
\\
&\times 
\int_{(\bR_{+}^\times )^{n-1}}\biggl(\int_{N_{n-1}}\int_{\rM_{n-1,1}(F)}
\me_{(n-1,n)}((t^{-1}x^{-1},t^{-1}x^{-1}z)a)\psi_{\varepsilon ,n-1}(x)\\
&\times 
\psi_{-\varepsilon }(e_{n-1}z)\,dz\,dx\biggr)
\widetilde{W}^{(\varepsilon )}_{\widehat{d},\widehat{\nu}}(t)
\prod_{i=1}^{n-1}t_i^{-(\nu_n+n-i)\rc_F-d_i+d_n}
\,\frac{2\rc_F\,dt_i}{t_i}.
\end{align*}
Let us consider the integral 
\begin{align*}
&\int_{N_{n-1}}\int_{\rM_{n-1,1}(F)}
\me_{(n-1,n)}((t^{-1}x^{-1},t^{-1}x^{-1}z)a)\psi_{\varepsilon ,n-1}(x)
\psi_{-\varepsilon }(e_{n-1}z)\,dz\,dx.
\end{align*}
Substituting $\displaystyle 
\left(\begin{array}{cc}
x^{-1}&x^{-1}z\\
O_{1,n-1}&1
\end{array}\right)\to x$, this integral becomes 
\begin{align}
\label{eqn:pf_app_002}
&\int_{N_{n}}
\me_{(n-1,n)}((1_{n-1},O_{n-1,1})\iota_n(t^{-1})xa)
\psi_{-\varepsilon ,n}(x)\,dx.
\end{align}
By the elementary formula (\ref{eqn:el_fmla_int}) and the equality 
\begin{align*}
&\me_{(n-1,n)} (\left(1_{n-1},O_{n-1,1}\right)\iota_n(t^{-1})xa)\\
&=\prod_{i=1}^{n-1}\exp (-\pi \rc_Ft_i^{-2}a_i^2)
\prod_{j=i+1}^n\exp (-\pi \rc_Ft_i^{-2}a_j^2x_{i,j}\overline{x_{i,j}}),
\end{align*}
we know that (\ref{eqn:pf_app_002}) is equal to  
\begin{align*}
&\prod_{i=1}^{n-1}\exp \biggl(-\pi \rc_F
\biggl(\frac{t_i^2}{a_{i+1}^{2}}+\frac{a_i^2}{t_i^2}\biggr)\biggr)
t_i^{(n-i)\rc_F}a_{i+1}^{-i\rc_F}.
\end{align*}
Therefore, we obtain the assertion.  
\end{proof}

\def\cprime{$'$}


\begin{thebibliography}{10}

\bibitem{Alisauskas_Jucys_001}
S.~J. Ali\v{s}auskas, A.-A.~A. Jucys, and A.~P. Jucys.
\newblock On the symmetric tensor operators of the unitary groups.
\newblock {\em J. Mathematical Phys.}, 13:1329--1333, 1972.

\bibitem{Dong_Xue_001}
Chao-Ping Dong and Huajian Xue.
\newblock On the nonvanishing hypothesis for {R}ankin-{S}elberg convolutions
  for {${\rm GL}_n(\Bbb C)\times{\rm GL}_n(\Bbb C)$}.
\newblock {\em Represent. Theory}, 21:151--171, 2017.

\bibitem{Gelfand_Tsetlin_001}
I.~M. Gel\cprime~fand and M.~L. Cetlin.
\newblock Finite-dimensional representations of the group of unimodular
  matrices.
\newblock {\em Doklady Akad. Nauk SSSR (N.S.)}, 71:825--828, 1950.

\bibitem{Goodman_Wallach_001}
Roe Goodman and Nolan~R. Wallach.
\newblock {\em Symmetry, representations, and invariants}, volume 255 of {\em
  Graduate Texts in Mathematics}.
\newblock Springer, Dordrecht, 2009.

\bibitem{Grenie_001}
Lo\"{\i}c Greni\'{e}.
\newblock Critical values of automorphic {$L$}-functions for {${\rm
  GL}(r)\times{\rm GL}(r)$}.
\newblock {\em Manuscripta Math.}, 110(3):283--311, 2003.

\bibitem{Grobner_Harris_001}
Harald Grobner and Michael Harris.
\newblock Whittaker periods, motivic periods, and special values of tensor
  product {$L$}-functions.
\newblock {\em J. Inst. Math. Jussieu}, 15(4):711--769, 2016.

\bibitem{Hecke_001}
Erich Hecke.
\newblock {\em Mathematische {W}erke}.
\newblock Herausgegeben im Auftrage der Akademie der Wissenschaften zu
  G\"{o}ttingen. Vandenhoeck \& Ruprecht, G\"{o}ttingen, 1959.

\bibitem{Hirano_Ishii_Miyazaki_pre}
Miki Hirano, Taku Ishii, and Tadashi Miyazaki.
\newblock Archimedean zeta integrals for {$GL(3)\times GL(2)$}.
\newblock {\em To appear in Memoir of the AMS.}

\bibitem{Ishii_Stade_001}
Taku Ishii and Eric Stade.
\newblock New formulas for {W}hittaker functions on {${\rm GL}(n,\Bbb R)$}.
\newblock {\em J. Funct. Anal.}, 244(1):289--314, 2007.

\bibitem{Ishii_Stade_002}
Taku Ishii and Eric Stade.
\newblock Archimedean zeta integrals on {$\mathrm{GL}_n\times \mathrm{GL}_m$}
  and {$\mathrm{SO}_{2n+1}\times \mathrm{GL}_m$}.
\newblock {\em Manuscripta Math.}, 141(3-4):485--536, 2013.

\bibitem{Jacquet_Langlands_001}
H.~Jacquet and R.~P. Langlands.
\newblock {\em Automorphic forms on {${\rm GL}(2)$}}.
\newblock Lecture Notes in Mathematics, Vol. 114. Springer-Verlag, Berlin-New
  York, 1970.

\bibitem{J_P_S_S_001}
H.~Jacquet, I.~I. Piatetskii-Shapiro, and J.~A. Shalika.
\newblock Rankin-{S}elberg convolutions.
\newblock {\em Amer. J. Math.}, 105(2):367--464, 1983.

\bibitem{Jacquet_003}
Herv{\'e} Jacquet.
\newblock {\em Automorphic forms on {${\rm GL}(2)$}. {P}art {II}}.
\newblock Lecture Notes in Mathematics, Vol. 278. Springer-Verlag, Berlin-New
  York, 1972.

\bibitem{Jacquet_001}
Herv{\'e} Jacquet.
\newblock Archimedean {R}ankin-{S}elberg integrals.
\newblock In {\em Automorphic forms and {$L$}-functions {II}. {L}ocal aspects},
  volume 489 of {\em Contemp. Math.}, pages 57--172. Amer. Math. Soc.,
  Providence, RI, 2009.

\bibitem{Jacuqet_Shalika_001}
Herv{\'e} Jacquet and Joseph Shalika.
\newblock Rankin-{S}elberg convolutions: {A}rchimedean theory.
\newblock In {\em Festschrift in honor of I. I. Piatetski-Shapiro on the
  occasion of his sixtieth birthday, Part I (Ramat Aviv, 1989)}, volume~2 of
  {\em Israel Math. Conf. Proc.}, pages 125--207. Weizmann, Jerusalem, 1990.

\bibitem{Jucis_001}
A.~A. Jucis.
\newblock The isoscalar factors of the {C}lebsch-{G}ordan coefficients of
  unitary groups.
\newblock {\em Litovsk. Fiz. Sb.}, 10:5--12, 1970.

\bibitem{Knapp_002}
Anthony~W. Knapp.
\newblock {\em Representation theory of semisimple groups}.
\newblock Princeton Landmarks in Mathematics. Princeton University Press,
  Princeton, NJ, 2001.
\newblock An overview based on examples, Reprint of the 1986 original.

\bibitem{Kostant_001}
Bertram Kostant.
\newblock On {W}hittaker vectors and representation theory.
\newblock {\em Invent. Math.}, 48(2):101--184, 1978.

\bibitem{Miyazaki_003}
Tadashi Miyazaki.
\newblock The local zeta integrals for {$GL(2,\bold{C})\times GL(2,\bold{C})$}.
\newblock {\em Proc. Japan Acad. Ser. A Math. Sci.}, 94(1):1--6, 2018.

\bibitem{Popa_001}
Alexandru~A. Popa.
\newblock Whittaker newforms for {A}rchimedean representations.
\newblock {\em J. Number Theory}, 128(6):1637--1645, 2008.

\bibitem{Raghuram_001}
A.~Raghuram.
\newblock Critical values of {R}ankin-{S}elberg {$L$}-functions for
  {$\text{GL}_n\times\text{GL}_{n-1}$} and the symmetric cube {$L$}-functions
  for {$\text{GL}_2$}.
\newblock {\em Forum Math.}, 28(3):457--489, 2016.

\bibitem{Rankin_001}
R.~A. Rankin.
\newblock Contributions to the theory of {R}amanujan's function {$\tau(n)$} and
  similar arithmetical functions. {I}. {T}he zeros of the function
  {$\sum^\infty_{n=1}\tau(n)/n^s$} on the line {$\Re s=13/2$}. {II}. {T}he
  order of the {F}ourier coefficients of integral modular forms.
\newblock {\em Proc. Cambridge Philos. Soc.}, 35:351--372, 1939.

\bibitem{Selberg001}
Atle Selberg.
\newblock Bemerkungen \"{u}ber eine {D}irichletsche {R}eihe, die mit der
  {T}heorie der {M}odulformen nahe verbunden ist.
\newblock {\em Arch. Math. Naturvid.}, 43:47--50, 1940.

\bibitem{Speh_Vogan_001}
Birgit Speh and David~A. Vogan, Jr.
\newblock Reducibility of generalized principal series representations.
\newblock {\em Acta Math.}, 145(3-4):227--299, 1980.

\bibitem{Stade_002}
Eric Stade.
\newblock Mellin transforms of {${\rm GL}(n,\Bbb R)$} {W}hittaker functions.
\newblock {\em Amer. J. Math.}, 123(1):121--161, 2001.

\bibitem{Stade_001}
Eric Stade.
\newblock Archimedean {$L$}-factors on {${\rm GL}(n)\times{\rm GL}(n)$} and
  generalized {B}arnes integrals.
\newblock {\em Israel J. Math.}, 127:201--219, 2002.

\bibitem{Sun_001}
Binyong Sun.
\newblock The nonvanishing hypothesis at infinity for {R}ankin-{S}elberg
  convolutions.
\newblock {\em J. Amer. Math. Soc.}, 30(1):1--25, 2017.

\bibitem{Vilenkin_Klimyk_001}
N.~Ja. Vilenkin and A.~U. Klimyk.
\newblock {\em Representation of {L}ie groups and special functions. {V}ol. 3},
  volume~75 of {\em Mathematics and its Applications (Soviet Series)}.
\newblock Kluwer Academic Publishers Group, Dordrecht, 1992.
\newblock Classical and quantum groups and special functions, Translated from
  the Russian by V. A. Groza and A. A. Groza.

\bibitem{Wallach_003}
Nolan~R. Wallach.
\newblock {\em Real reductive groups. {II}}, volume 132 of {\em Pure and
  Applied Mathematics}.
\newblock Academic Press Inc., Boston, MA, 1992.

\bibitem{Watanabe_001}
Takao Watanabe.
\newblock Global theta liftings of general linear groups.
\newblock {\em J. Math. Sci. Univ. Tokyo}, 3(3):699--711, 1996.

\bibitem{Zhang_001}
Shou-Wu Zhang.
\newblock Gross-{Z}agier formula for {${\rm GL}_2$}.
\newblock {\em Asian J. Math.}, 5(2):183--290, 2001.

\bibitem{Zhelobenko_001}
D.~P. Zhelobenko.
\newblock On {G}el\cprime fand-{Z}etlin bases for classical {L}ie algebras.
\newblock In {\em Representations of {L}ie groups and {L}ie algebras
  ({B}udapest, 1971)}, pages 79--106. Akad. Kiad\'o, Budapest, 1985.

\end{thebibliography}

\end{document}